# A Conjecture Regarding the Riemann Hypothesis as Visualized by t Strings

by

Ronald F. Fox
Smyrna, Georgia
May 7, 2020

**Abstract**

The introduction of strings into the study of the Riemann Hypothesis provides a visualization of the genesis of zeros of the Zeta function. The method is heuristic and when originally introduced suggested strong visual evidence for the truth of the Riemann Hypothesis. Some sort of organizing principle for strings with similar $t$ values is observed and points towards a region between (1, 0) and (0, 0) on the abscissa, and within order unity along the ordinate. Progress in understanding these observations has been made by expanding the domain of $\sigma$ from the critical strip, [0, 1], to the half-line [0, ∞]. The nature of the organizing principle is explained. A generic structure for strings over the expanded domain is proffered. New perspective is gained regarding the truth of the Riemann Hypothesis.

# I. Introduction

The Riemann hypothesis remains one of the outstanding problems in analytic number theory after more than 160 years [1]. It is concerned with the distribution of prime numbers and is encapsulated by the Riemann Zeta function and an equivalent expression written exclusively in terms of primes found by Euler [2] much earlier (for real s).
(1)

$$\zeta(s) = \sum_{n=1}^{\infty} \frac{1}{n^s} = \prod_{p\ prime} \frac{1}{1-p^{-s}}$$

in which $s = \sigma + ti$, a complex number with real part $\sigma$ and imaginary part coefficient $t$. This is the notation of Riemann. The expressions above are valid for $\sigma > 1$ and are absolutely convergent there. In [1] Riemann used analytic continuation to extend the Zeta function into the regime $0 < \sigma < 1$ (and into the regime $\sigma < 0$) and found the Dirichlet Eta function equation (conjectured by Euler in 1749)
(2)

$$\eta(s) = \sum_{n=1}^{\infty} (-1)^{n-1} \frac{1}{n^s} = (1 - 2^{1-s})\zeta(s)$$

This result involves a conditionally convergent series and must be handled with care. In the form
(3)

$$\zeta(s) = \frac{1}{(1-2^{1-s})} \sum_{n=1}^{\infty} (-1)^{n-1} \frac{1}{n^s} = \frac{1}{(1-2^{1-s})} \eta(s)$$

it provides a way of computing the Zeta function for $0 < \sigma < 1$ in terms of the alternating Eta series.

Also, in [1] Riemann established the *reflection formula for the Zeta function* for $0 < \sigma < 1$,
(4)

$$\zeta(s) = 2^s \pi^{s-1} \sin\left(\frac{\pi s}{2}\right) \Gamma(1-s) \zeta(1-s)$$

Combining Eq.(3) with Eq.(4) yields the *reflection formula for the Eta function*

(5)

$$\eta(s) = 2^s \pi^{s-1} \sin\left(\frac{\pi s}{2}\right) \Gamma(1-s) \frac{(1-2^{1-s})}{(1-2^s)} \eta(1-s)$$

These reflection formulae imply that if $s$ is a zero of Zeta (or Eta) then so is $1-s$. There are zeros of Zeta and Eta of three kinds. There are some so-called trivial zeros for $s = \sigma = -2, -4, \ldots$ created by the sine term. There are also trivial zeros associated with the factor $\frac{(1-2^{1-s})}{(1-2^s)}$. These are of the form $1 - i\frac{k 2\pi}{\ln(2)}$ for integer $k$. The remaining *non-trivial* zeros of Zeta are also non-trivial zeros for Eta and *vice versa*. The *critical strip* defined by $0 < \sigma < 1$ and $-\infty < t < \infty$ contains these non-trivial zeros. The Riemann hypothesis is: **$\sigma = 0.5$ for all non-trivial zeros.** In 2004, Gourdon [3] verified the hypothesis for the first $10^{13}$ zeros using an algorithm invented by Odlyzko [4].

There is also a *modified* reflection formula. Both the Zeta and the Eta functions have the special property [5] that if $s_0$ is a zero then so is its complex conjugate, $s_0^*$. Together with the reflection formula associated with Eqs.(4 - 5) this implies that $1 - s_0^*$ is also a zero. This modified reflection formula can be expressed:

*if $\sigma_0 + t_0 i$ is a zero of Eta then so is $1 - \sigma_0 + t_0 i$.*

These two expressions have the *same* imaginary parts and that is a key to what follows.

We now introduce the idea of a $t$ string. In the critical strip, $0 < \sigma < 1$ and $-\infty < t < \infty$, fix any value of $t$ and plot the value of Eta for that $t$ and all $\sigma \in [0, 1]$. If instead, we use Zeta to compute the string values we get a $t$ string for Zeta. In the following we will use Eta. As far as the nontrivial zeros are concerned there is no difference between the two choices. The extension of the results described up to now make the choice of Eta natural as will be seen.

As an example, the strings in figure {1}

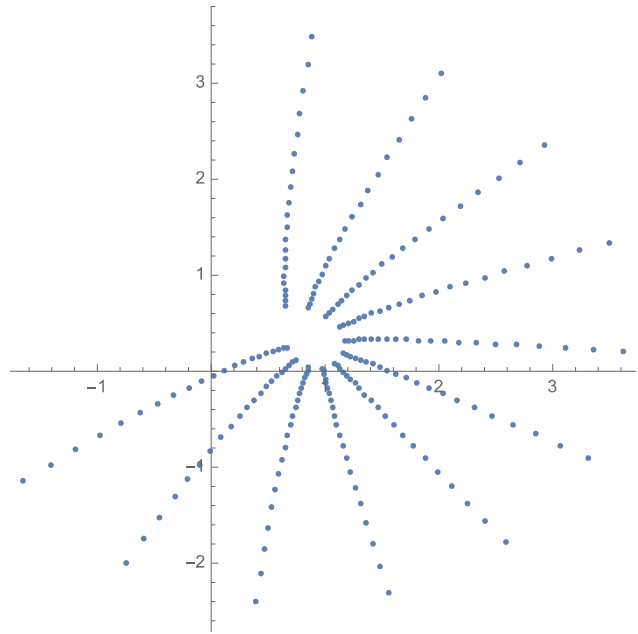

{1}

were generated using Mathematica 12 to execute the command

$$Table[DirichletEta[\sigma + ti], \{\sigma, 0, 1, .05\}, \{t, 19, 21, .2\}]$$

This command means calculate the value of the Dirichlet Eta function for the $t$ value 19 and all the $\sigma$'s in $\{\sigma, 0, 1, .05\}$ which means a discrete set of $\sigma$ values from 0 to 1 in steps of 0.05, and then do it again for the $t$ value 19.2, and again for 19.4 and so on up to the $t$ value 21. This creates 11 strings each of which is composed of 21 $\sigma$ values. The choice of discrete values of $\sigma$ has an advantage relative to a using a continuous line of $\sigma's$. The string created from a uniformly distributed set of $\sigma$ values is not uniform. The denser end of the string is always for the larger values of $\sigma$ and the more rarified end is for the smaller values of $\sigma$. The strings are laid down clockwise as the value of $t$ increases. In this example the string with $t$ value 19 is the most vertical one pointing at just past 12 o'clock.

Two more observations about this example lead to the results to be presented below. There is a zero with $t = 21.022039639...$ Very close to this value is the string for $t = 21$, the one furthest around clockwise from the first string at $t = 19$. By counting dots, from the outer most to the inner most, there are 21 dots (the last two are fused at this magnification) and the one corresponding with $\sigma = 0.5$ is the 11th one. It is almost on the origin (0, 0) the condition required by the Riemann Hypothesis (RH). I can show this explicitly by plotting the string for $t = 21.022039639…$ in figure {2}

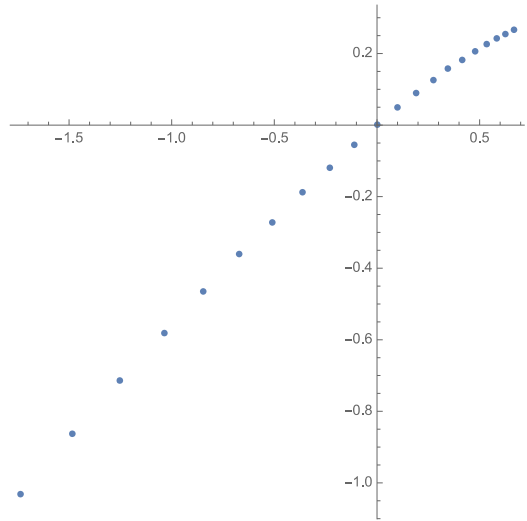

{2}

The aspect ratios in the two figures are not identical so there is an apparent difference in the slope of the two strings for $t = 21$ and for $t = 21.022039639…$, beyond their numerical difference that is only 0.1%. The second observation has to do with how I get these plots from the Mathematica Table command. The output from the Table command organizes the data in accord with the instructions intended by $\{\sigma, 0, 1, .05\}$, $\{t, 19, 21, .2\}$. In doing so, it introduces internal parentheses reflecting these instructions. Removal of these parentheses reduces the data to a collection of coordinates that can be plotted automatically by clicking on the plot option. When the sizes of the sets of discrete entities in the two instructions is different, as in this example, placing the larger set to the right of the smaller set produces fewer parentheses than the other way around. In this example a reversal of the order changes the number of pairs of internal parentheses from 21 to 11. This means the more efficient way to execute this example involves 11 sets of internal parentheses, each containing 21 entities.

As we compute the $t$ strings for increasing values of $t$ some of the time a string comes close to the origin (0, 0). By adjusting the value of $t$ we can make the string and the origin coincide. As an example of how strings increase our ability to visualize the genesis of zeros, we look at the behavior leading up to the first zero. Displaying the behavior of Eta visually greatly enhances our ability to understand how and why the non-trivial zeros of Eta (Zeta) occur for $\sigma = 0.5$. We begin by restricting the argument of Eta to the real axis of the critical strip ($t = 0$). For $\sigma \in (0, 1)$ we find $\eta(\sigma) \in \big(0.5,\ ln(2)\big)$. The unit interval for $\sigma$ is compressed by $\eta$ into an interval of length $0.193\ldots$. Because $\eta$ is continuous and infinitely differentiable in both $\sigma$ and $t$, we expect that small variations in $t$ will create small variations in $\eta$. We know that the first zero of eta has an imaginary part slightly bigger than 14. Therefore, we will begin by looking at what happens to a set of points for $\{\sigma,\ 0.02, 0.98, 0.02\}$ and for a *fixed* value of $t$. The choice of discrete values of $\sigma$ is imposed in order to keep the magnitude of the computation reasonable and to illustrate non-uniformity in the output. For figure $\{3\}$, the choice of spacing works pretty well and the sets of points for a fixed $t$ represent a *string* very well. Some discreteness does appear in the longer strings (larger $t$ values) and shows a non-uniform density.

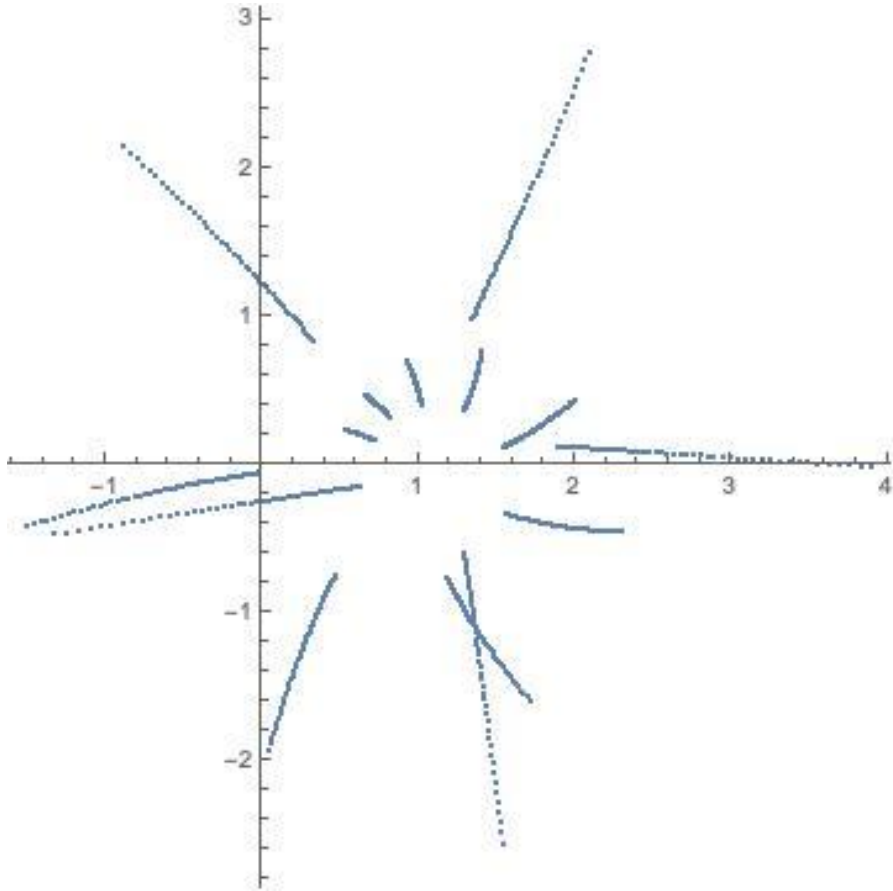

{3}

There is a lot of information in this figure. Each string is labeled by $t$ and is made up of the values of Eta for fixed $t$ and for 49 evenly spaced values of $\sigma$ with the 25th value equal to 0.5. There are 14 strings in the figure corresponding to 14 values of $t$ in $\{t, 1, 14, 1\}$. The $t = 1$ string is located near the 1 on the abscissa at about 10 o'clock. It is the shortest string in the plot, slightly larger than 0.193 (the arc length for the case $t = 0$). As we go clockwise, we see the strings for $t = 2, 3, 4, 5,$ and 6. Clearly, they are getting longer in arc length and are moving away from the origin. String 6 appears below the abscissa at about 4 o'clock. The points on the string are arranged by Eta with the small values of $\sigma$ distal to the 1 on the abscissa and the larger values of $\sigma$ proximal to the 1 on the abscissa. Since the Eta values for uniformly spaced $\sigma$'s are not uniformly spaced the $\sigma = 0.5$ value of Eta is not halfway along the length of the string. By enlarging the figure, the location of this special $\sigma$ can be found. The 9th and the 14th strings are close to each other just before 9 o'clock. The 9th string is the shorter of the two. It's $\sigma = 1$ value is closest to the origin. It looks like a slight adjustment in the $t$ value would put the $\sigma = 1$ right on the origin. Indeed, $t = 9.0647\ldots$ does do so. This corresponds precisely with a *trivial* zero given above by $1 - i\frac{k2\pi}{\ln(2)}$ for the case of $k = 1$. Similarly, the 14th string looks like an adjustment of its $t$ value would possibly put $\sigma = 0.5$ on the origin. Let us try $t = 14.134725\ldots$. This is plotted below in figure {4}.

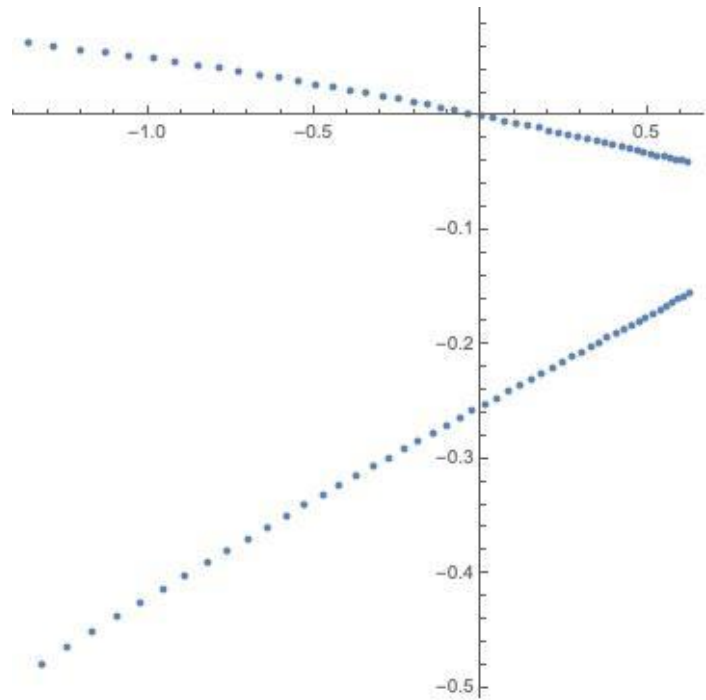

{4}

The lower string is the $t = 14$ string from the previous plot. The difference in apparent slope is the result of different aspect ratios in the two plots. The upper curve is for $t = 14.134725$, the imaginary part of first non-trivial zero for Eta. Moreover, $\sigma = 0.5$ corresponds with the 25th dot in the 49-dot representation of the $\sigma$ interval (it can be located by counting from the left). As you can see the typical value for a dot is a real part in the tenths and an imaginary part in the tenths as well. However, the 25th $\sigma$ dot has the Eta value $1.62123 \times 10^{-6} - 2.6635 \times 10^{-7} i$. Since we have expressed the $t$ for the first zero of Eta to only one part in $10^6$ we cannot expect to get "zero" to any better precision. The smaller values of $\sigma$ produce the points to the left in the figure and the larger values of $\sigma$ produce the points to the right in the figure. After it became clear that the trivial zeros on the edge of the critical strip did not cause problems the sigma range was changed to $\{\sigma, 0, 1, .05\}$.

The question that shouts out from this account is why does the string intersect the origin at $\sigma = 0.5$ ?! This is where the modified reflection formula comes into play.

*A line and a point can intersect only once (?).*

Let them intersect for $\sigma_0 + t_0 i$, i.e. $\eta(\sigma_0 + t_0 i) = 0$. The reflection formula implies that they (*same $t_0$ string*) also intersect at $1 - \sigma_0 + t_0 i$. If they can intersect only once, then these two expressions must describe the same point. That is possible only if $\sigma_0 = 0.5$ (RH). After executing several more special cases I happened to discover that a string could loop around and cross itself. If such a crossing point *coincided with the origin* then it is possible for a string to intersect a point, the origin, at two different points of the string, corresponding to the crossing values, $\sigma_0'$ and $\sigma_0''$. While I did see self-crossings of strings, I never saw one at the origin. Nevertheless, I did see a case that was remarkably tantalizing.

From Odlyzko's data [6] there is a zero with $t = $ 267653395649.3623669687. There is also a self-crossing close to the origin. I have observed before that for larger $t$ values the length of the string gets longer at the $\sigma < 0.5$ end. This particular string has a length of order 15000. Most of it cannot be shown in figure {5} below because the loop takes up a square area of about 8 ($2 \times 4$). The crossing is "close" to the origin, at (0.04067, 0.0936). For comparison the values of the coordinates closest to the origin where the zero is manifested are (-0.003229, -0.001320).

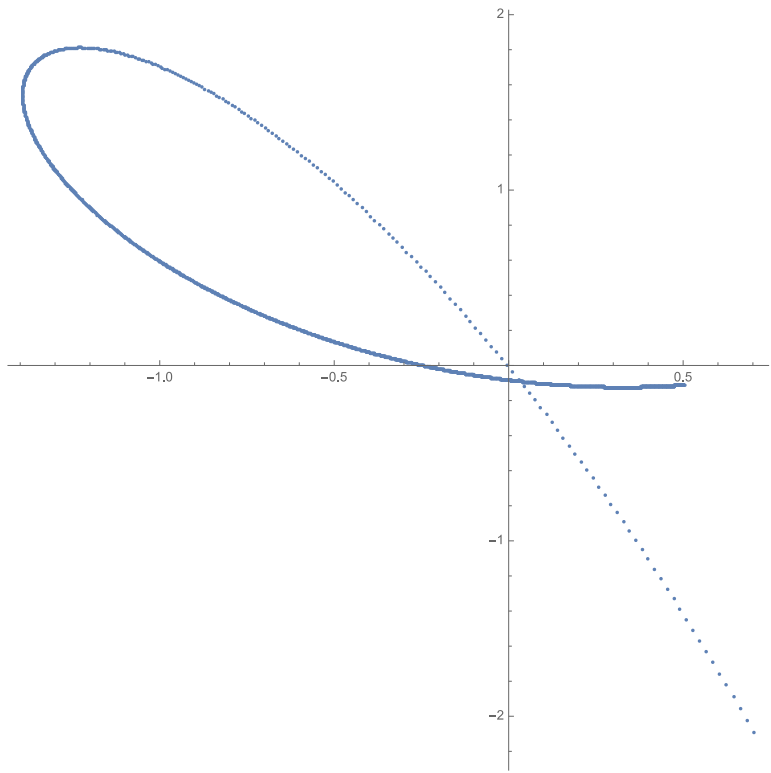

{5}

The zeros are discrete. Thus, it is not possible to make a small adjustment in the $t$ value and create a small adjustment in the crossing point of a desired amount. Nearby $t's$ for other zeros either produce no loop or produce one crossing further away from the origin, at least for the cases I explored.

A study of many strings for many values of $t$ up to one in the trillions suggested that there was some sort of organizing principle not centered on the origin, (0, 0). For example, the case defined by (figure {6})

$$Table[DirichletEta[\sigma + ti], \{t, 21, 23, .2\}, \{\sigma, 0, 1, .05\}]$$

produces 11 strings that appear to be organized around the point (1.5, 0.5), at least approximately. In the other examples given below (and in the appendix) I will simply write down the instructions for $\sigma$ and $t$ and not the front part stating Table and DirichletEta.

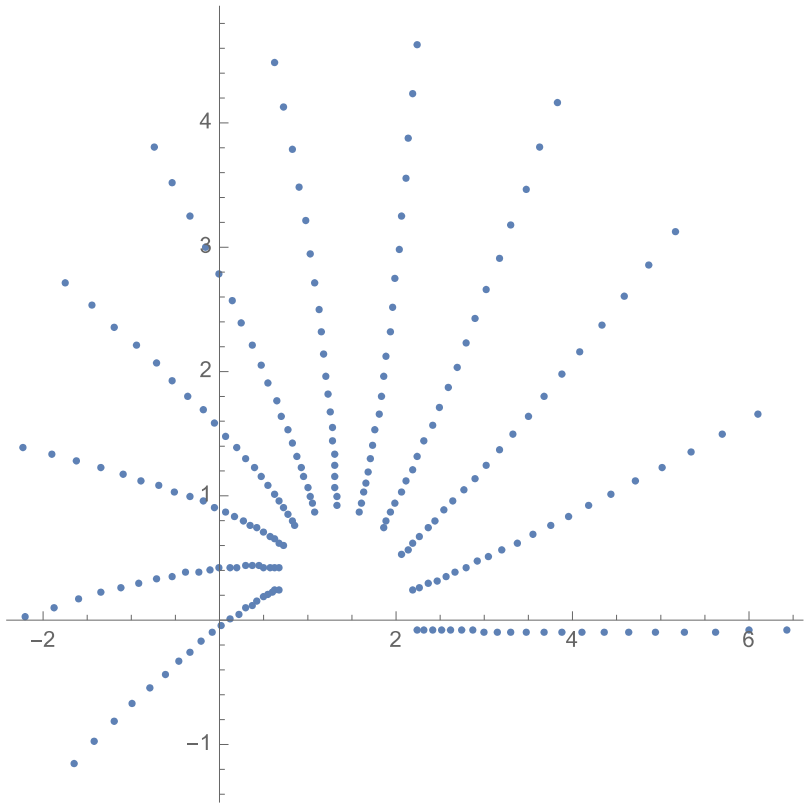

{6}

A more striking example is given by $\{t, 24, 26, .2\}, \{\sigma, 0, 1, .05\}$ in figure {7}. Clearly, there is a zero very close to $t = 25$.

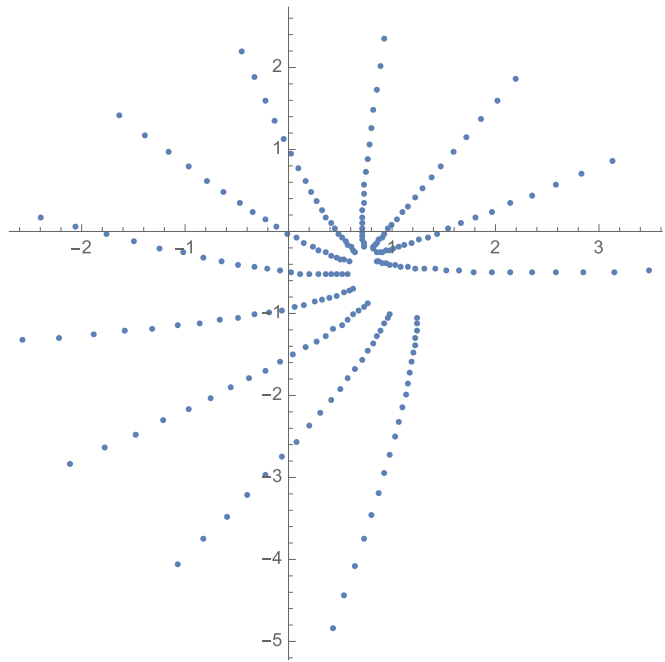

{7}

Half of the strings in figure {7} appear to be organized around one center and the other half are organized around a different center. A more extreme version of this type of behavior is given in figure {8}.

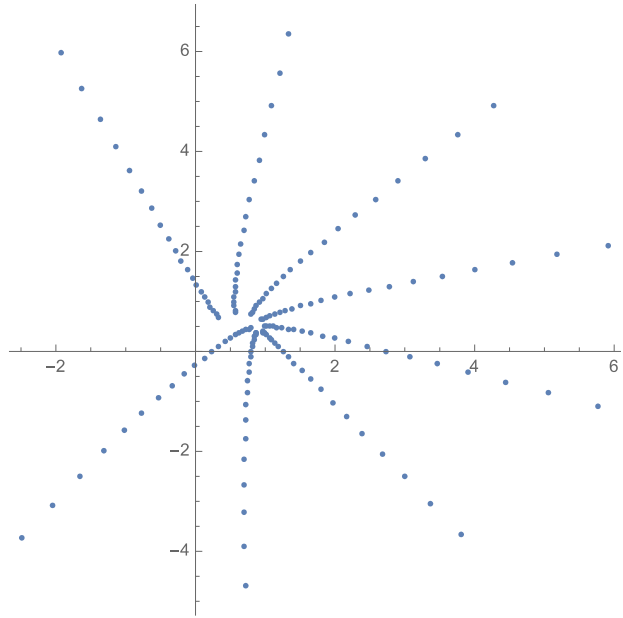

{8}

Figure {8} is defined by $\{t, 55, 56.4, .2\}, \{\sigma, 0, 1, .05\}$. How do we account for these structures? Why are the strings nearly straight rather than much more highly curved? Note that in figures {7-8} the axes are very close to the same scale so that there is no apparent suppression of curvature created by different scales for the two axes. In the **appendix**, the details of how the zeros occur as $t$ is increased from 1 to 67.6 are exhibited.

To answer these questions and related ones, we expand the domain for $\sigma$ from [0, 1] to [0, ∞]. This will make it possible to see the structures produced for $\sigma$ in [0, 1] as having their genesis in the [1, ∞] portion of the $\sigma$ domain. All of this will be presented in Section II of the paper. In Section III the remainder of the justification for RH is explored. Can self-crossings of $t$ strings occur at the origin? Can there be more than one self-crossing per string?

## II. Expansion of the σ domain

Expand the σ domain for Eta from [0, 1] to [0, ∞]. This means that for σ in [0,1] we are thinking about the Eta function. Eta is alternating on [0, ∞] and is absolutely convergent on [1, ∞]. Zeta, by contrast, is absolutely convergent for σ in [1, ∞] and converted by analytic continuation [1] to a simple function times Eta in [0, 1], Eq. (3). Moreover, Zeta is a positive series on [1, ∞] whereas Eta is an alternating series there. Thus, we have come back to where Zeta began but continue with the alternating series instead. What we learn about Eta for σ in [1, ∞] that can be used to understand Eta for σ in [0, 1] will also apply to Zeta for σ in [0, 1].

For σ in [0, ∞], the definition of Eta is
(6)
$$Eta[\sigma + ti] = \sum_{n=1}^{\infty} (-1)^{n-1} \frac{1}{n^{\sigma+ti}}$$

$$= 1 + \sum_{n=2}^{\infty} (-1)^{n-1} \frac{1}{n^{\sigma+ti}}$$

Eta, for every choice of $\sigma + ti$, begins with 1. Each summand has a natural expression in polar coordinates of the form
(7)
$$\frac{1}{n^{\sigma+ti}} = n^{-\sigma} Exp[-itLog[n]]$$

wherein Log denotes the Napierian logarithm. These polar forms are defined relative to the origin, (0, 0). Of course, the 1 is trivially also in polar form. The minus signs in the even $n$ terms can be absorbed into the angle factor as $e^{-i\pi}$. The addition of polar forms with respect to a fixed origin, in this case, (0, 0), is closed, i.e. the result is always expressible as a polar form with the same origin. For fixed σ the length of the polar form ($n^{-\sigma}$) is decreasing as $n$ increases, and the angle of the polar form ($tLog[n]$) is increasing in the clockwise direction (a negative angle) as $n$ increases.

A "large" value for σ is 10. Let's see how 7 strings are arrayed for { σ, 9, 10, 0.1},{t, 22, 28, 1}. We run

Table[DirichletEta[σ + ti], { σ, 9, 10, 0.1},{t, 22, 28, 1}]

and get

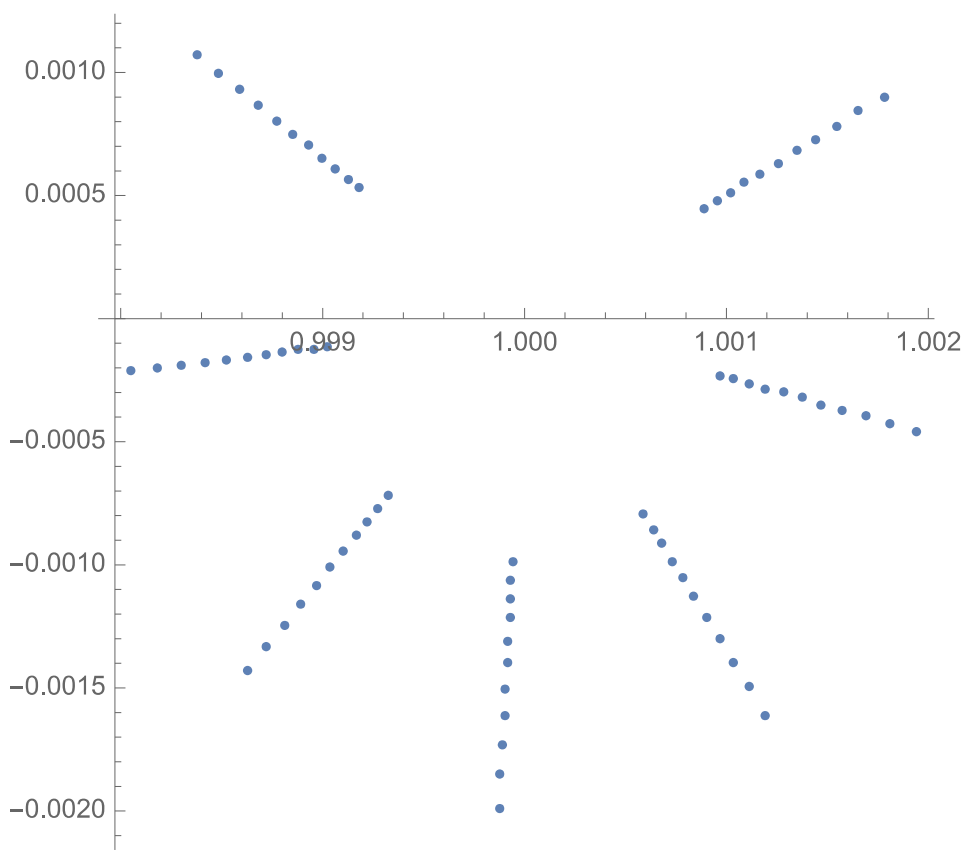

{9}

The strings are straight, they point towards (1, 0), the string at 2 o'clock is the $t = 22$ string, and they are of length ~ 0.001. Can we understand all of this?

The initial terms in the Eta summation are the biggest terms overall. For example, for $\sigma = 9$ we have

(8)

$$\eta(9 + 22i) \cong 1 - 2^{-9} Exp[-i\, 22 Log[2]] + 3^{-9} Exp[-i\, 22 Log[3]]\ldots$$

in which the first term is 1, the second has a modulus of order 1/512 (0.001953125) and the third is of order 1/19683 (0.0000508053). Clearly the value of Eta for $\sigma \geq 9$ is dominated by the first two terms. Since the 1 appears in every evaluation of Eta, we can think of it as a shift of the origin, from (0, 0) to (1, 0), and drop further consideration of the 1. Figure {9} exhibits the dominance of the $n = 2$ term, at least to a few parts per 1000. Moreover, this single term works for all the values of $\sigma$ in $\{\sigma, 9, 10, 0.1\}$ since they are larger than 9, making the contributions to Eta smaller. Each string has an angle exponential that factors from the magnitude term in this limit. That makes all the $\sigma$ terms for a given $t$ have the same angle and appear as straight lines pointing to (1, 0). This translates in the present case into

(9)

$$-2^{-\sigma} Exp[-i\, 22 Log[2]], \{\sigma, 9, 10, 0.1\}$$

The angle in radians is $-22 \times Log[2] = -22 \times .693... = -15.249...$ This is more negative than $-4\pi$ by $-2.6828...$ radians or $-153.71...$ degrees. The effect of the minus sign in front of the expression above is to add $\pi$ radians to the angle. This means the angle for the $t = 22$ string is 26.29… degrees. In figure {9}, the two axes are not the same scale. They aren't off by much but to the eye it looks like the angle is close to 50 degrees. From the raw data, the coordinates of the points at the ends of the $t = 22$ string are (1.00088, 0.000445673) at the $\sigma = 10$ end and (1.00178, 0.000904055) at the $\sigma = 9$ end. For a straight line in between, the angle is ArcTan [0.000458382/0.00090] = 26.99… degrees. The discrepancy is 70/2674 = 2.6% which is the same size as the contribution of the $n = 3$ term to the $n = 2$ term. We include the $n = 3$ term in Eq. (8) and get small corrections to the approximation for the DirichletEta of figure {9}. The corresponding coordinates of the ends of the string are (1.00089, 0.000446332) and (1.00178, 0.000906560). The resulting angle satisfies ArcTan [0.000460228/0.00089] = 27.34… degrees. This analysis can be repeated for each of the remaining 6 strings in figure {9}.

Another check on the accuracy of the approximation involves the straightness assumption. Since the strings point to (1, 0) we can compare the coordinates of the end points by looking at the angle made by each end relative to (1, 0). The raw data for the DirichletEta plot is given above for $t = 22$ and is (1.00088, 0.000445673) and (1.00178, 0.000904055) for $\sigma = 10$ and $\sigma = 9$ respectively. Relative to (1, 0) the two end point angles are respectively ArcTan [0.000445673/0.00088] = 26.85… degrees and ArcTan [0.000904055/0.00178] = 26.92… degrees. Clearly the string is not perfectly straight because of the small admixture of the $n = 3$ term with the dominate $n = 2$ term. In addition, the string may not be pointing at $(1, 0)$ perfectly for these $\sigma$ values. For $\sigma = 20$, the relative importance of the $n = 3$ term to the $n = 2$ term is $3^{-20}$ to $2^{-20}$ or 0.0003 (0.03%). Note that all of the strings for $\sigma$ values between 20 and $\infty$ plot inside the $\sigma = 10$ ends of the strings in figure {9}. Expanding the $\sigma$ domain into $[1, \infty]$ only increases the extent of the string plots by order 1 corresponding to the organizing center at precisely $(1, 0)$.

To make the last point completely clear we plot in figure {10} the results for

Table[DirichletEta[$\sigma + ti$], { $\sigma$, 19, 20, 0.1},{t, 22, 28, 1}]

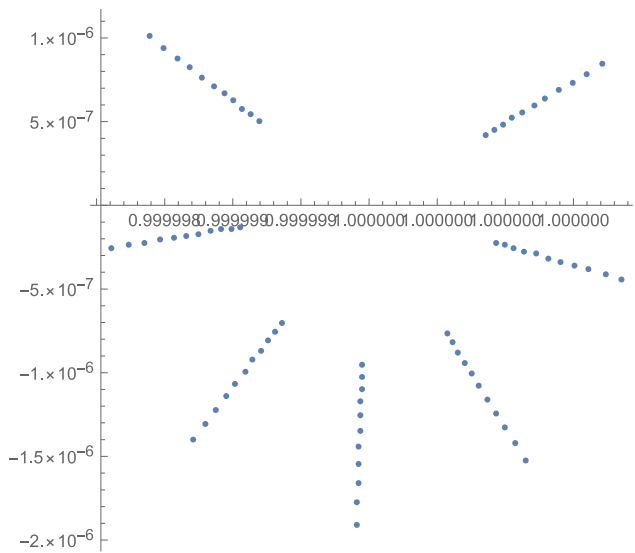

{10}

Notice that the figures {9-10} look the same except for the difference in scale. As $\sigma$ get larger the scale gets smaller, exponentially. On the other hand, as $\sigma$ gets smaller eventually the structure of the string pattern changes dramatically. Since many scales are still manifest it is necessary to look at different stretches of $\sigma$ to get the full picture. In anticipation of what is to come later let us look at three plots:

Table[DirichletEta[$\sigma + ti$], {t, 22, 28, 1},{ $\sigma$, 1.5, 4, 0.1}]

Table[DirichletEta[$\sigma + ti$], {t, 22, 28, 1},{ $\sigma$, 0.5, 1.5, 0.1}]

Table[DirichletEta[$\sigma + ti$], {t, 22, 28, 1},{ $\sigma$, 1, 2, 0.1}]

The first plot in figure {11} is for σ greater than 1 and, therefore, still outside the critical strip.

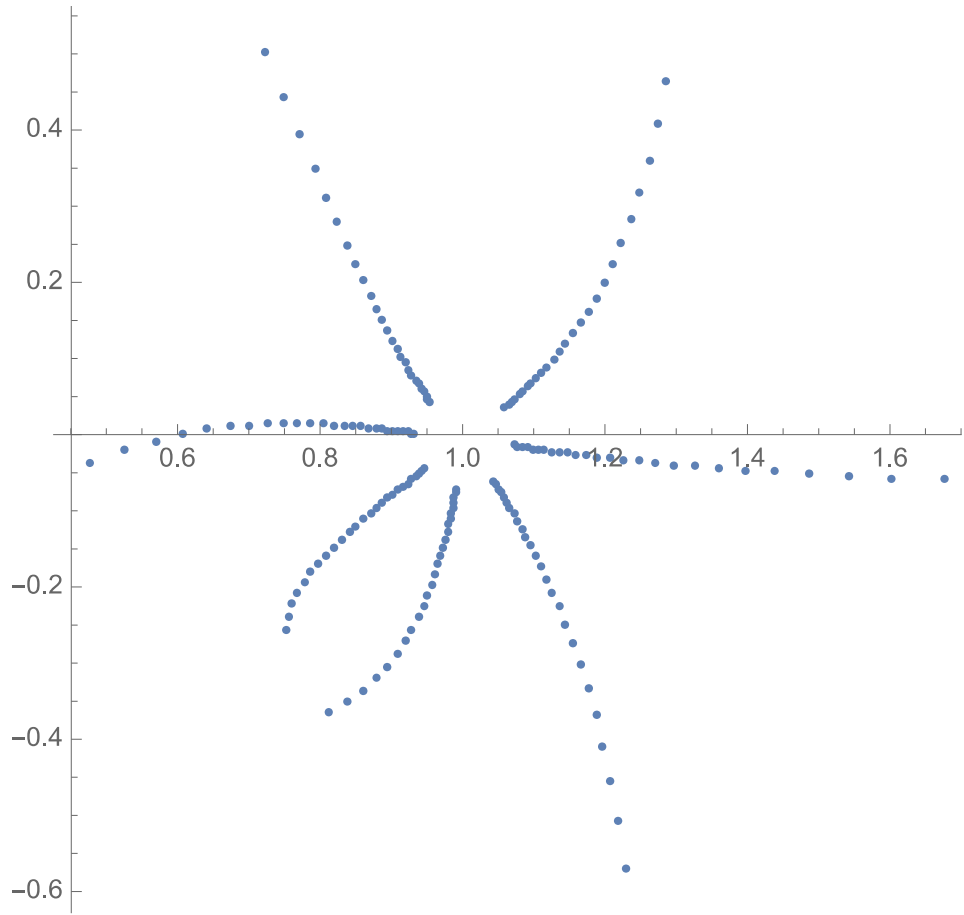

{11}

One can see that the outer dots in figure {9} match up with the inner dots (a smear) in figure {11}. However, the outer dots in figure {11} are showing some curvature that represents the influence of higher values of $n$ than 2 and 3. The second plot in figure {12} shows more structure and straddles the critical strip boundary at σ = 1.

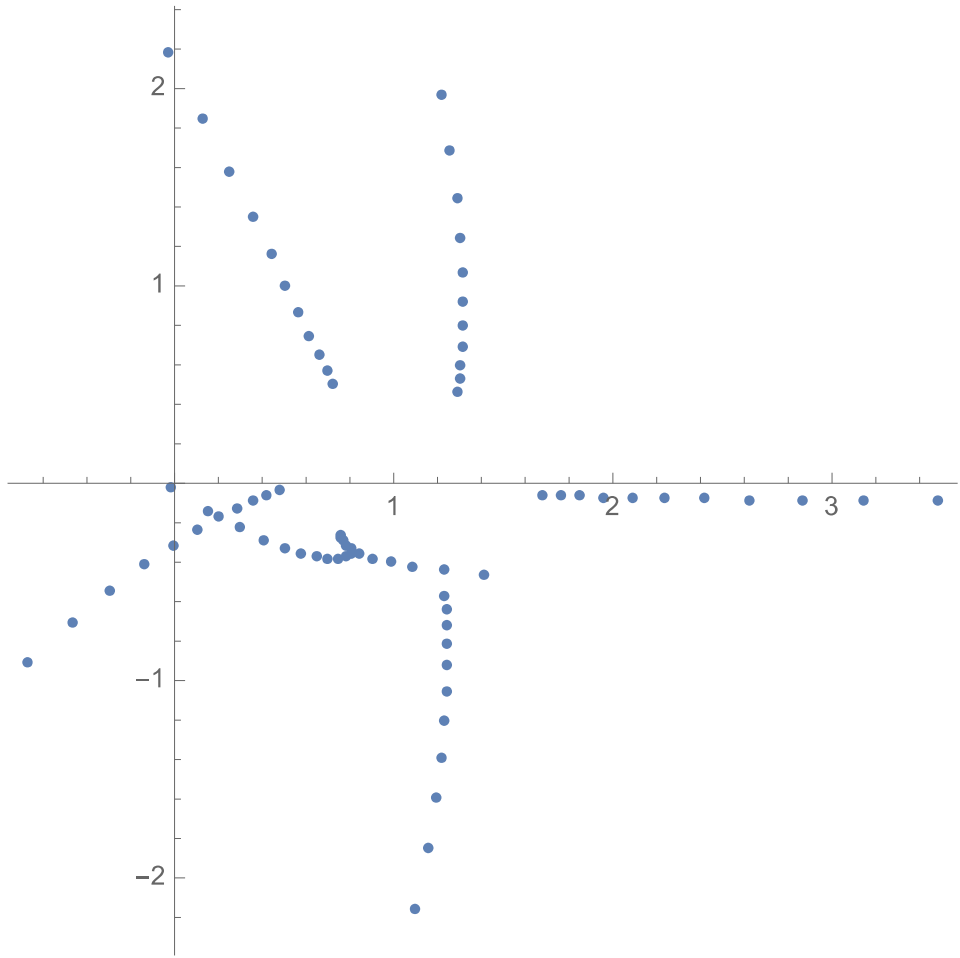

{12}

Figures {11-12} are on different scales. The outer dots of figure {11} smoothly connect to the inner dots of figure {12}. It is not clear from these plots where the two horizonal strings near (1, 0) in figure {12} connect. The third plot in figure {13} makes it clear. In figure {11} the two strings (lower left) curving towards each other are for $t = 25$ and $t = 26$. In figure {12} the $t = 25$ string is curving to the left while the $t = 26$ string is hooking to the right. The "direction" of these "motions" is from dense end to rarified end. Thus, in figure {12} the short horizontal string just below the 1 on the abscissa is the $t = 26$ string and the short string to its left is the $t = 25$ string. These two strings have switched places and in figure {12} it is the $\sigma = 0.5$ end of the $t = 25$ string that is almost coinciding with the origin (0, 0). This means that there must be a zero of Eta (Zeta) nearby. In fact, there is a zero with $t = 25.010857580...$

The third plot is given by figure {13} and its σ domain overlaps those of the first two plots

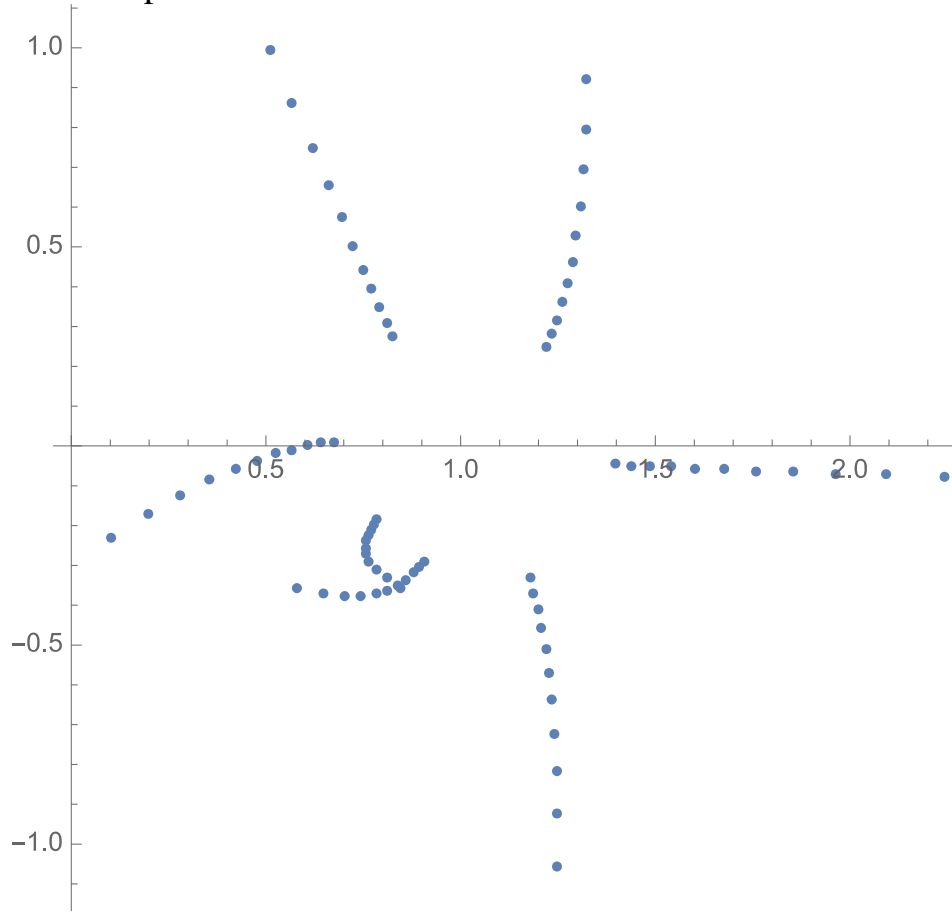

{13}

This figure makes the behavior of strings for $t = 25$ and $t = 26$ clearer. We can fuse these three views into one composite figure. If we look at

$$\text{Table}[\text{DirichletEta}[\sigma + ti], \{t, 22, 28, 1\}, \{\sigma, 0, 4, 0.1\}]$$

we get figure {14} which smoothly joins the fragments previously exhibited. Some strings have changed very little and some have changed greatly.

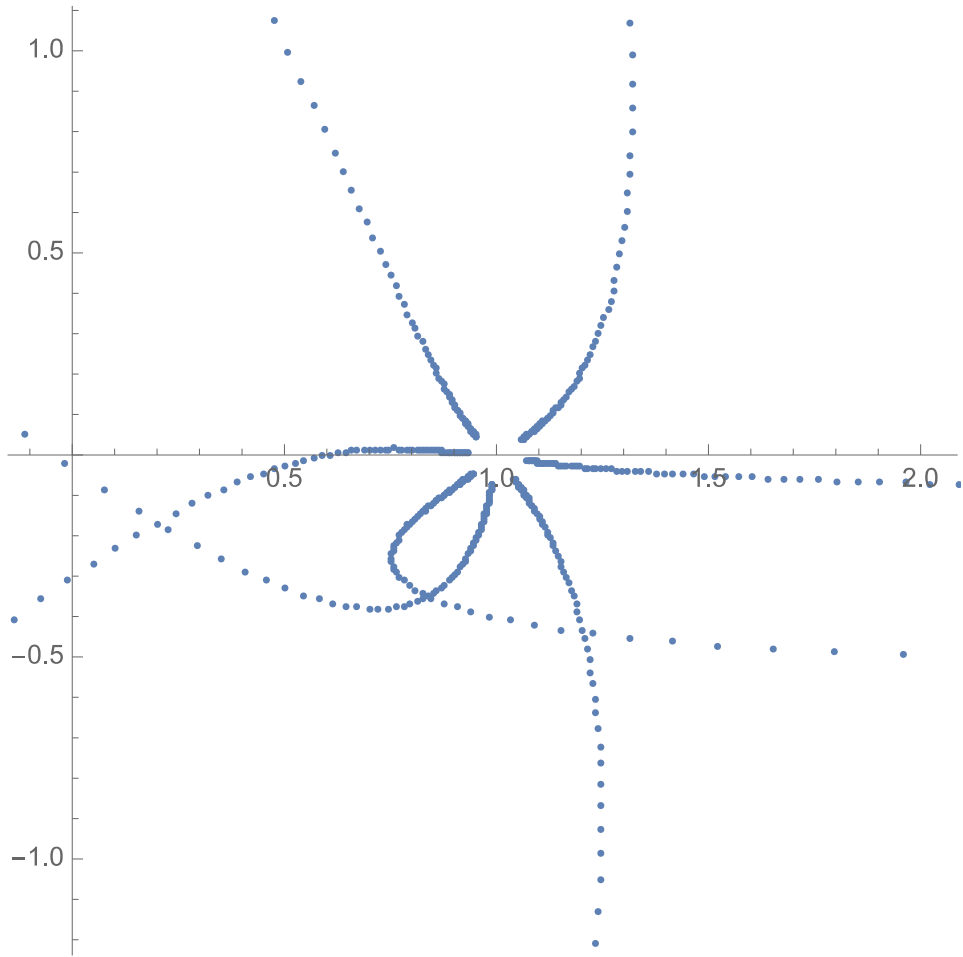

{14}

The cores of this figure are given in figures {9-10} but on a very much smaller scale.

    Two concepts with very different foundations will be used to gain more understanding of the string figures. We want a quantitative manner to explain the emergence and transformation of strings when they no longer remain essentially straight or only slightly curved. We also will use the subjective notion that the human eye is not able to do better than a linear perception of 0.1%. This means that in one dimension, on a page of the size we are looking at now, we do not perceive variations as small or smaller than 0.1%. In two dimensions the precision is the square of this, or 0.0001%. I will demonstrate the truth of the perception issue as we go along using examples from the figures.

To determine how many terms in the sum that constitutes the Eta function we need to keep for a specified precision we will compare the moduli of each term to that of the $n = 2$ term. We will have precision p for $n$ summands if

(10)
$$\frac{\frac{1}{(n+1)^\sigma}}{\frac{1}{2^\sigma}} < 10^{-p}$$

and

$$\frac{\frac{1}{n^\sigma}}{\frac{1}{2^\sigma}} > 10^{-p}$$

This means that $n$ is the last integer lacking the desired precision and $n + 1$ is the first that does have the desired precision. Simple algebra renders these expressions in the useful form: $n = 2 \times 10^{p/\sigma}$. As long as $\sigma > p$, $n$ will remain relatively small, but when the inequality turns around there will be a rapid exponential (exponential of reciprocal $\sigma$) increase in $n$.

A few examples of how well a truncated series represents the Eta function are given below. The first example uses a $t$ range that covers two consecutive zeros, 111.0295... and 111.8746... 10 uniformly spaced values of $t$ are used. Start with "large" $t$.
Table[DirichletEta[$\sigma + ti$],{t, 111.0295, 111.8746, 0.0939},{ $\sigma$, 4, 7, 0.02}]

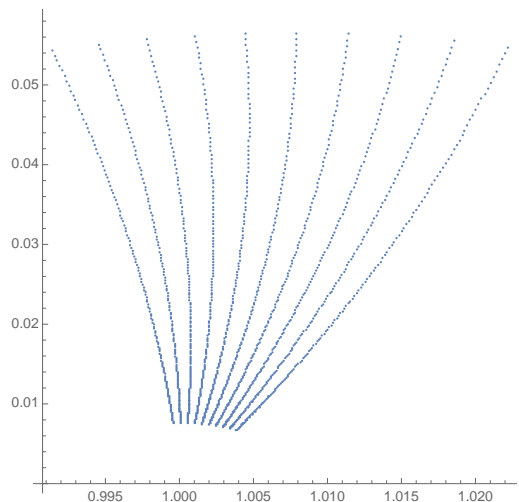

{15}

The result is plotted in figure {15}. Above, I said the values of $t$ were uniformly distributed between the values for the zeros. When we get to the $\sigma$ values around 0.5 it will be clear that the distribution of strings is ~ uniform around a circle determined by two consecutive zeros, the strings for which must go through (0, 0). For values of $\sigma$ near to 4 or 7 the strings are arranged as a parallel flare. In figure {15} the flare is seen to originate near (1, 0). Decreasing the $\sigma$ values produces the continuation of figure {15}:

Table[DirichletEta[$\sigma + ti$],{t, 111.0295, 111.8746, 0.0939},{ $\sigma$, 1.5, 4, .02}]

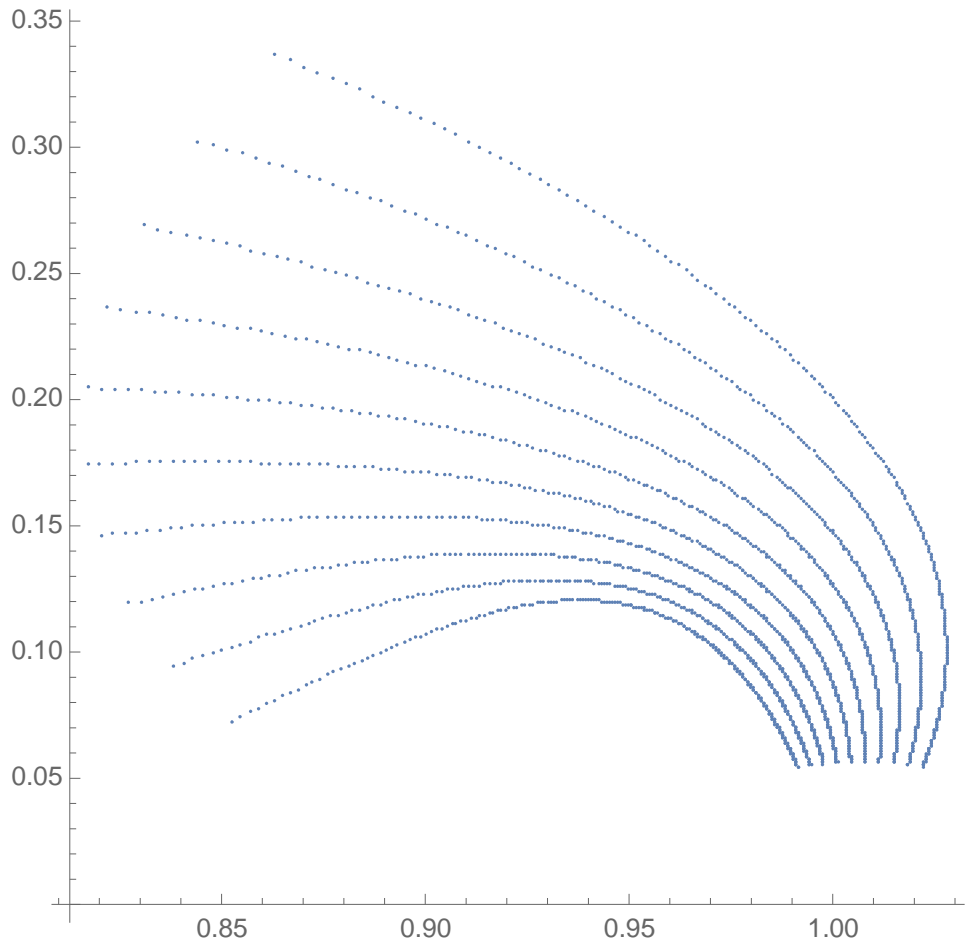

{16}

Aside from the difference in scale, figure {16} smoothly joins figure {15}. Moreover, the strings have taken a turn towards the critical strip that is still far away on this scale. After the turn, the strings resume their parallel flare. Finally, we cross into the critical strip and witness dramatic changes.
Table[DirichletEta[$\sigma + ti$],{t, 111.0295, 111.8746, 0.0939},{ $\sigma$, .4, 1.5, .01}]

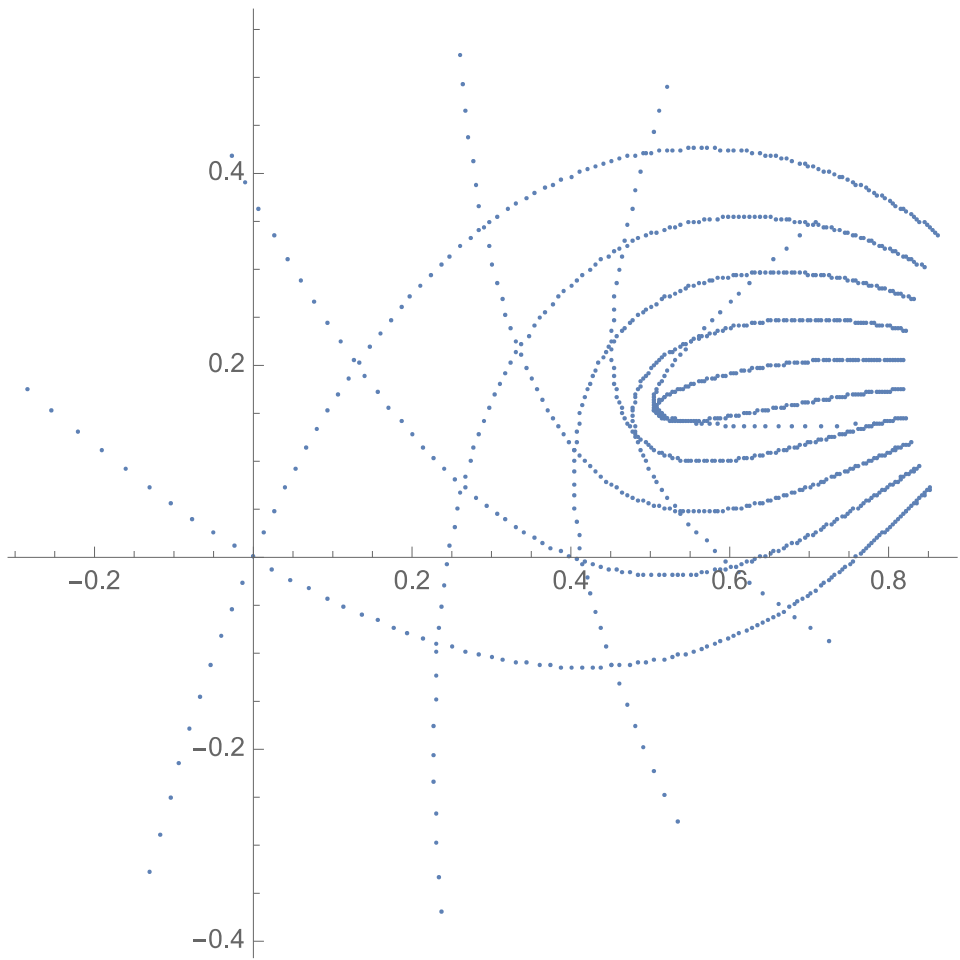

{17}

Once again, if you compensate for the change in scale, figures {16-17} smoothly connect. The parallel flare continues into the critical strip where it explodes into roughly uniformly distributed strings around a circle that is not centered on (0, 0) nor on (1, 0). Remember that the outer most two strings correspond with $t$ values for zeros and must go through (0, 0) as they plainly do. The $\sigma$ dots are spaced by 0.1 units apart and the dots for 0.4 to 0.5 are the outer most 11 dots of each string. If you just look at these outer 11 dots you will see the apparent circulation around a center (a radial flare) somewhere near (0.4, 0.1). Also notice that two of the strings have had to make almost complete U-turns in order to achieve the equal spacing around the circle. The transformation from parallel flare to radial flare seems to be universal, based on a small number of explored cases.

The preceding analysis used the Eta function provided by Mathematica 12. Let's see how well a truncated series of specified precision works in producing images that are indistinguishable to the human eye.

Table[Sum[$(-1)^{n-1}\left(\frac{1}{n^{\sigma+ti}}\right)$, {n, 1, 12, 1}],
 {t, 111.0295, 111.8746, 0.0939}, {σ, 4, 7, 0.01}]

Note that the limit on the Sum over $n$ is $n = 12$. Since the smallest value of $\sigma$ in this computation is 4, the precision formula gives: $n = 2 \times 10^{3/4} = 11.24$, which is rounded up to the integer 12. Figure {18} shows the result.

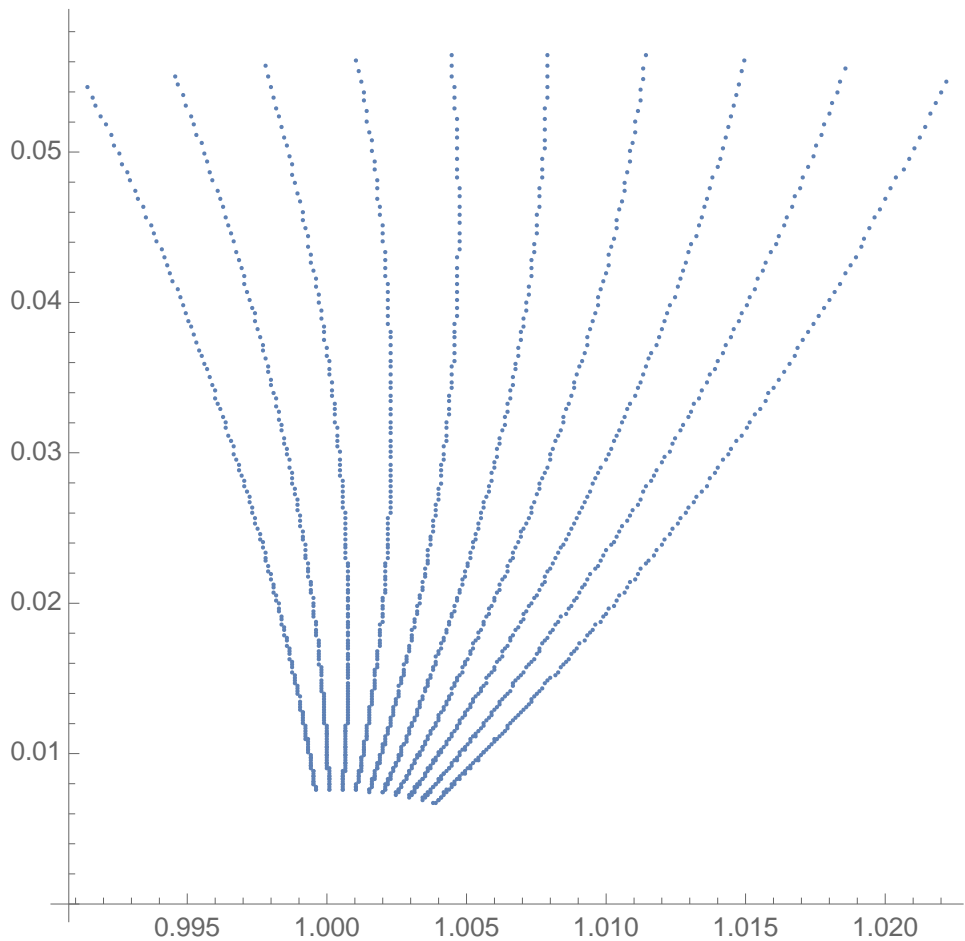

{18}

Next we do:

Table[Sum[$(-1)^{n-1}\left(\frac{1}{n^{\sigma+ti}}\right)$, {n, 1, 200, 1}],
 {t, 111.0295, 111.8746, 0.0939}, {σ, 1.5, 4, 0.01}]

For this segment the smallest value of $\sigma$ in this computation is 1.5, so that the precision formula gives $n = 2 \times 10^{3/1.5} = 200$. Note that the turn towards the critical strip is reproduced in figure {19}.

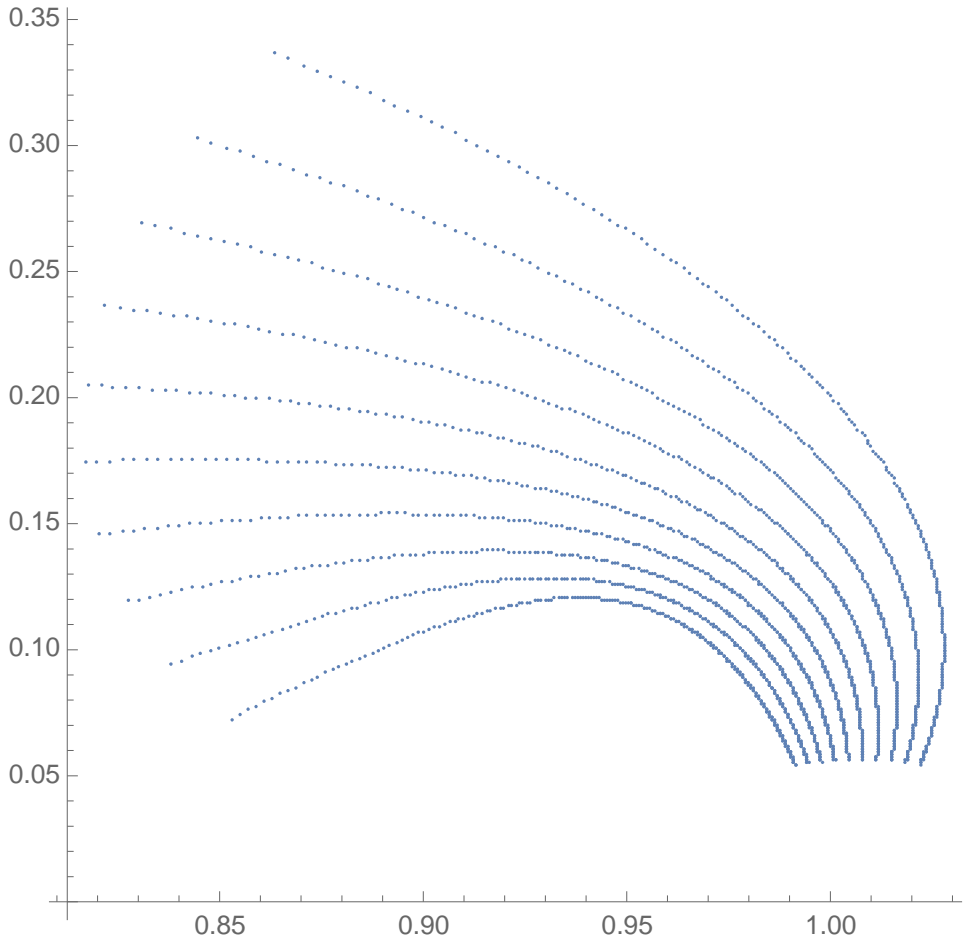

{19}

These first two figures can be smoothly joined if they are put on the same scale. If one generates the raw numbers, as I have done, then one can see the agreement with the Eta function calculations up to the desired precision for each segment. Is there going to be agreement for the next segment wherein for Eta dramatic changes occurred in the string pattern? We need:

Table[Sum[$(-1)^{n-1} \left(\frac{1}{n^{\sigma+ti}}\right)$, {n, 1, 2000000, 1}],
{t, 111.0295, 111.8746, 0.0939}, {$\sigma$, 0.4, 1.5, 0.01}]

Even though the smallest $\sigma$ value is 0.4, I have used 0.5 instead. For this choice we get $n = 2 \times 10^{3/0.5} = 2000000$. When you use 0.4 as would be expected the time to run the computation takes 160 hours. To do 0.5 takes only 5 hours. We expect to see deviations for $\sigma$ less than 0.5.

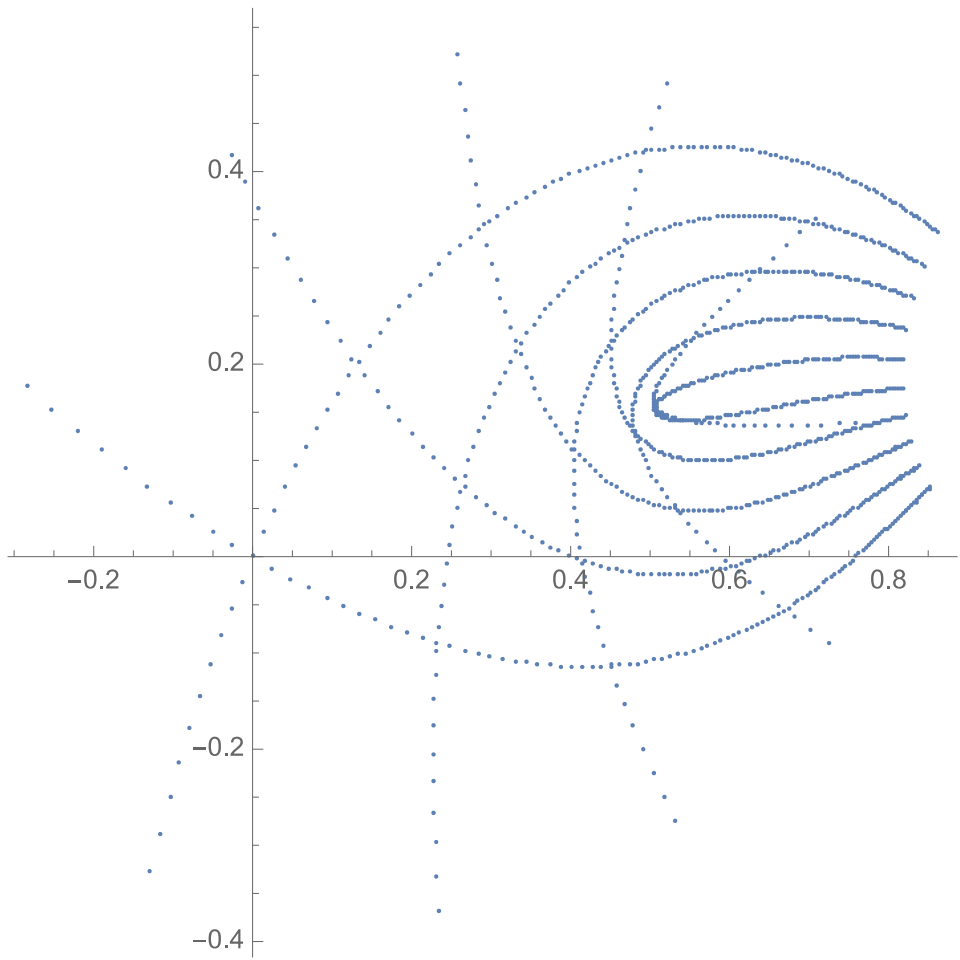

{20}

The agreement between the full Dirichlet Eta function calculations and the truncated approximations is remarkably good.

    I have mentioned that the structure generated for $\sigma \in [0, 0.5]$ is a radial flare. In figures {17, 20} one can see the beginnings of this since $\sigma$ is as small as 0.4. However, if one looks at the full domain $\sigma \in [0, 0.5]$ then the magnitudes of the strings at the 0 end of this domain are so much larger than for the rest of the larger values of $\sigma$ that the structure in figure {17} is indiscernible. (While the analog to figure {20} can also be contemplated, we do not consider it explicitly because the minimum $n$ needed to maintain the precision for p = 3 approaches ∞ as $\sigma$ approaches 0.)
Table[DirichletEta[$\sigma + ti$],{t, 111.0295, 111.8746, 0.0939},{ $\sigma$, 0, 0.5, .01}] generates figure {21}.

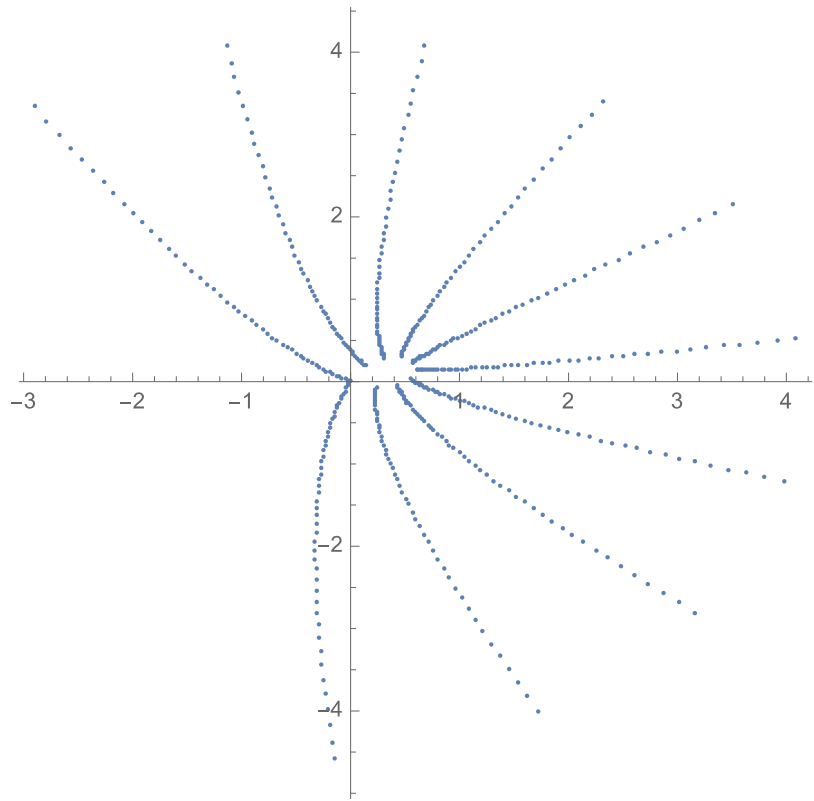

{21}

The scales for figures {20-21} differ by ten-fold. By overlapping the $\sigma$ domains of figures {20-21} the difference is clearer in figure {22}:
Table[DirichletEta[$\sigma + ti$],{t, 111.0295, 111.8746, 0.0939},{ $\sigma$, 0, 0.7, .01}]

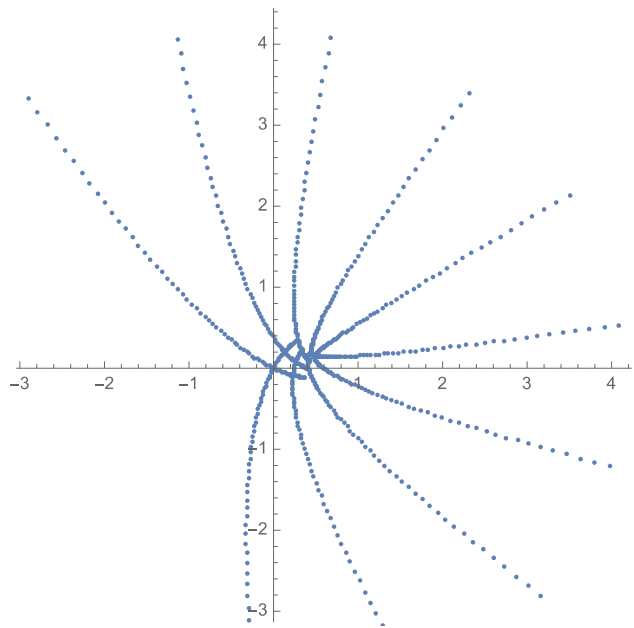

{22}

Also recall that these structures emanated from parallel flares that originated at (1, 0). This is very nearby on the scale of figure {22}.

Another example for two consecutive zeroes with $t = 357.151302252$ and $t = 357.952685102$ is given here by
Table[DirichletEta[$\sigma + ti$],{t, 357.151, 357.952, 0.089},{ $\sigma$, 4, 7, .02}]

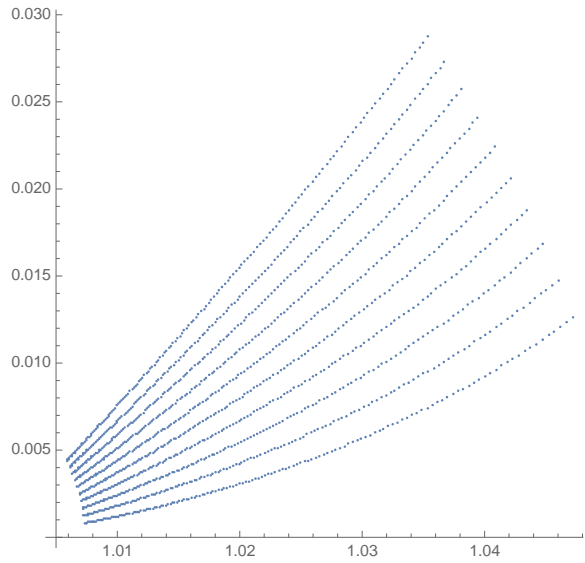

{23}

The approximation Sum is given by:
Table[Sum[$(-1)^{n-1} \left(\frac{1}{n^{\sigma+ti}}\right)$, {n, 1, 12, 1}],
{t, 357.151, 357.952, 0.089}, {$\sigma$, 4, 7, 0.01}]

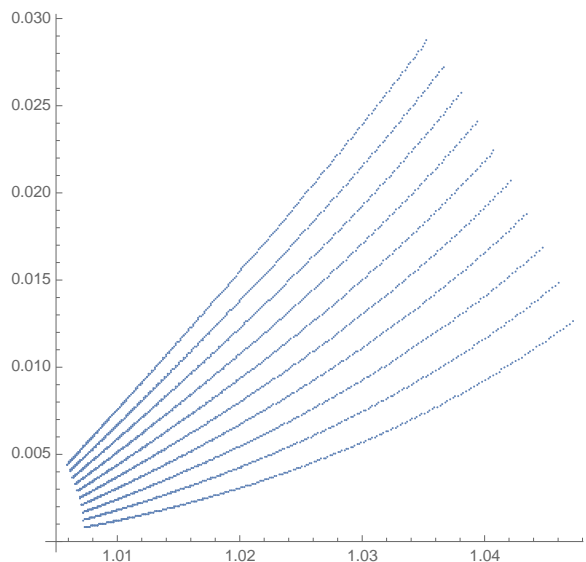

{24}

These strings are headed away from (0, 0). They smoothly connect to (DirichletEta, {σ, 1.5, 4, 0.01})

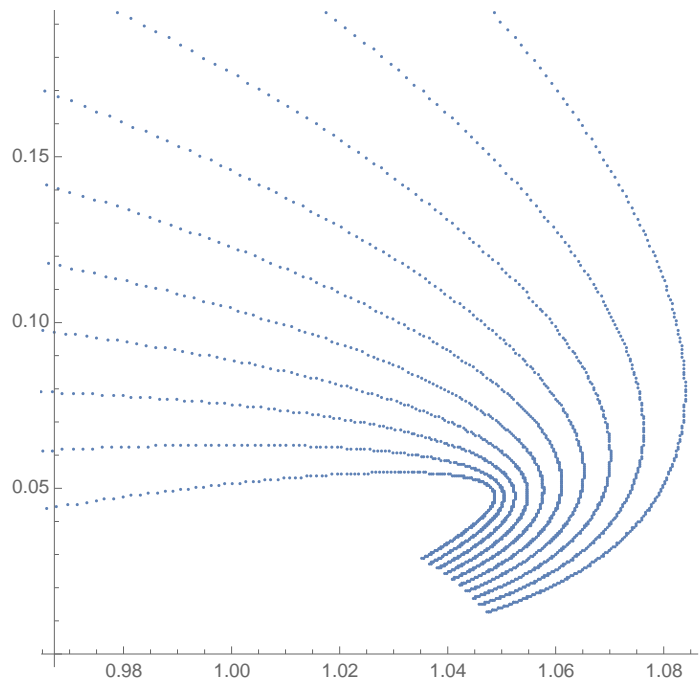

{25}

and to (Sum, {n, 1, 200, 1}, {σ, 4, 7, 0.01})

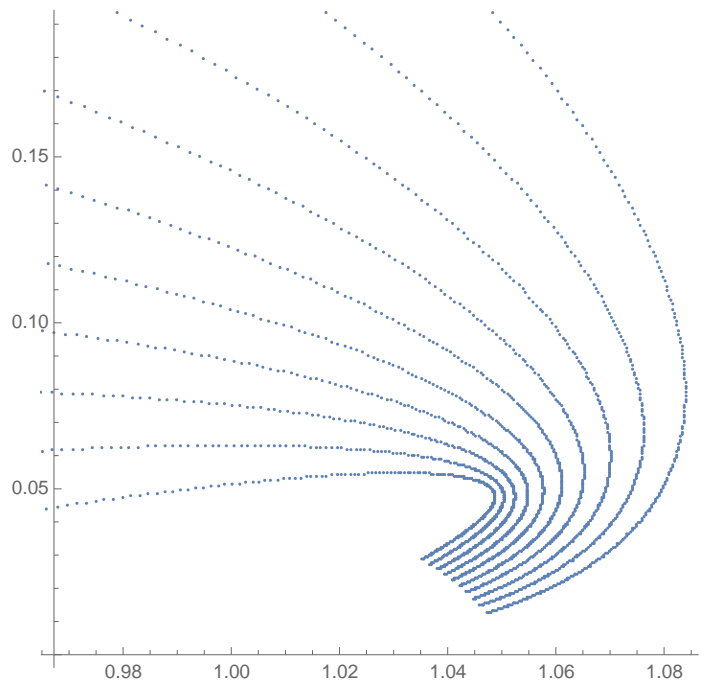

{26}

The short-hand notation for these calculations should be plain enough. The abrupt U-turn with the strings headed for the origin is captured very accurately by the 200 summands approximation (good for precision $p = 3$). The conversion from parallel flare to radial flare that was already begun by $\sigma = 0.5$ in the previous example is well on its way in this example in figure {27}. (DirichletEta, $\{\sigma, 0.4, 1.5, 0.01\}$)

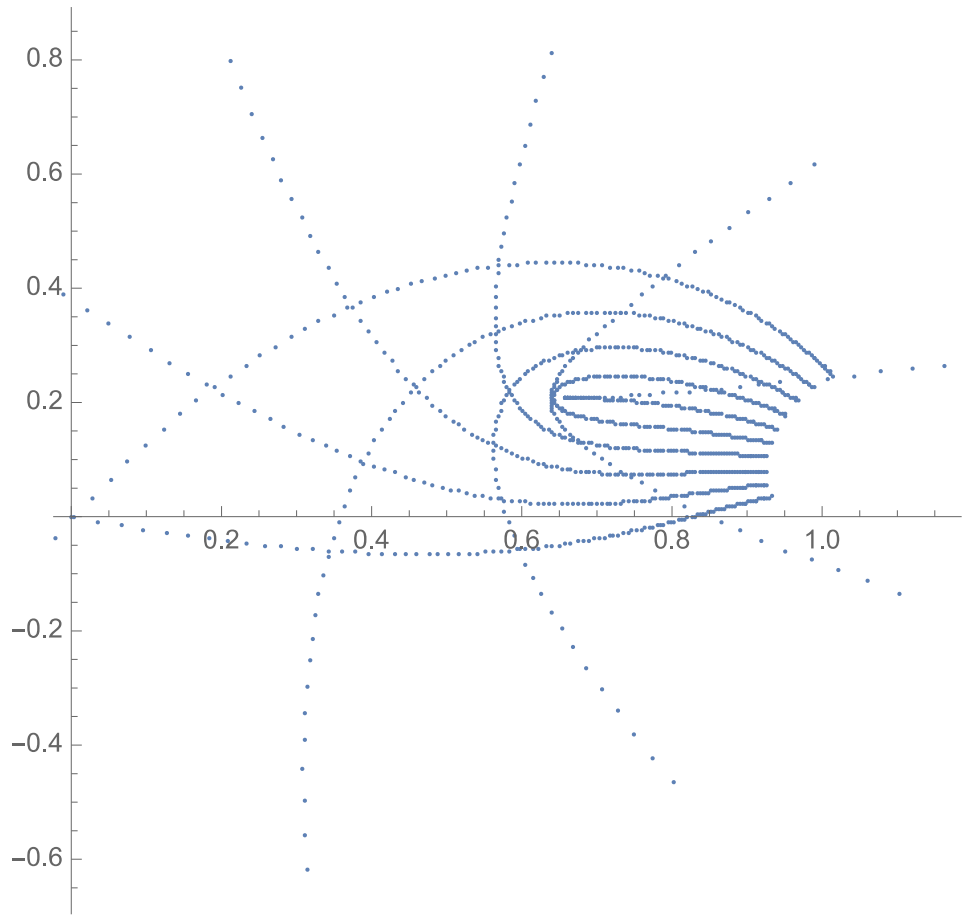

{27}

The corresponding Sum approximation is characterized in shorthand by (Sum, $\{n, 1, 2000000, 1\}, \{\sigma, 0.4, 1.5, 0.01\}$). As before, I chose the precision that goes with $\sigma = 0.5$ rather than $\sigma = 0.4$ because of the length of time needed for the more precise case. In both figures {27, 28} there is a fine detail that is barely visible. On the right-hand side you can still see the parallel flare formation. There are 10 strings in it. The sixth string from the bottom ends at a point at its left end. Coming back to the right from that end point is a rarified set of dots. This is an abrupt U-turn with a crossing. These details will be exhibited shortly. They are preserved in the approximation Sum shown in figure {28}.

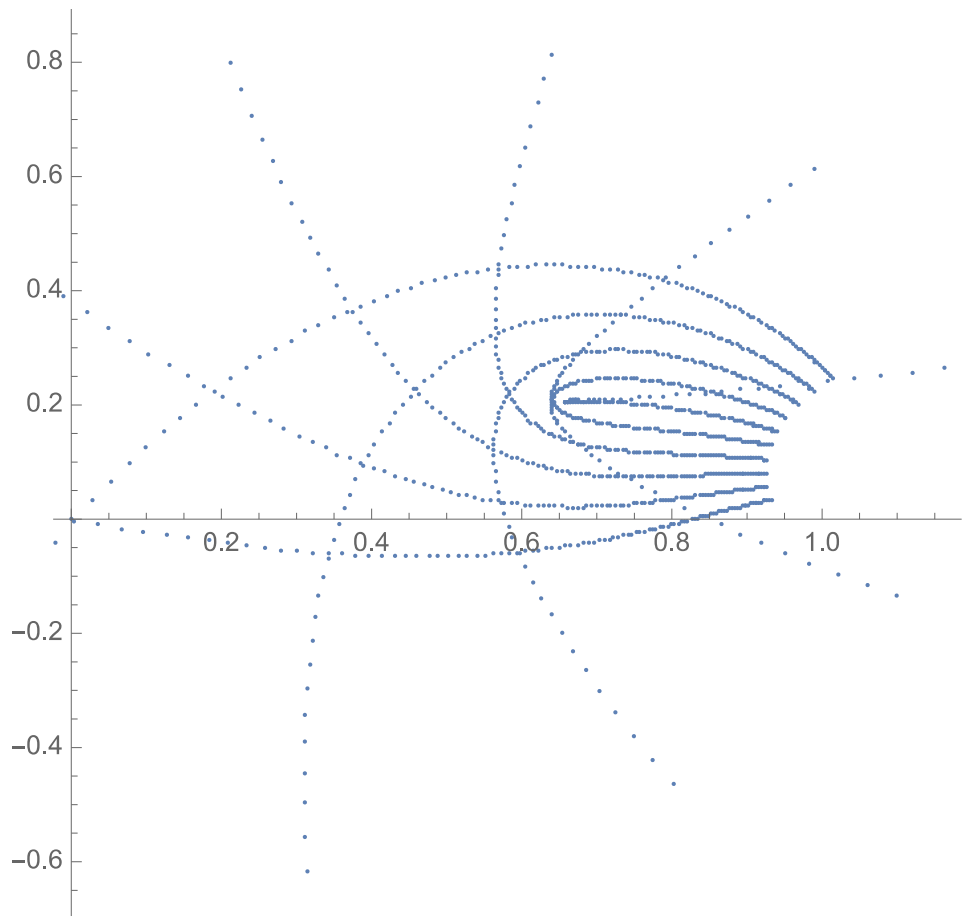

{28}

The U-turn in the preceding figure can be viewed by itself and also be magnified. The sixth string up from the bottom is at $t = 357.151 + 5 \times 0.089 = 357.596$. If I run Table[DirichletEta[$\sigma + ti$],{ $t$, 357.596, 357.596, .089},{ $\sigma$, .4, 1.5, .01}] then I get a magnified view of the U-turn with a much greater splay angle than in figure {28} because of the finer vertical scale used in figure {29}. While it may at first look like a cusp it is in fact a small loop. The ordinate scale in figure {28} is 5 times finer than is the abscissa scale.

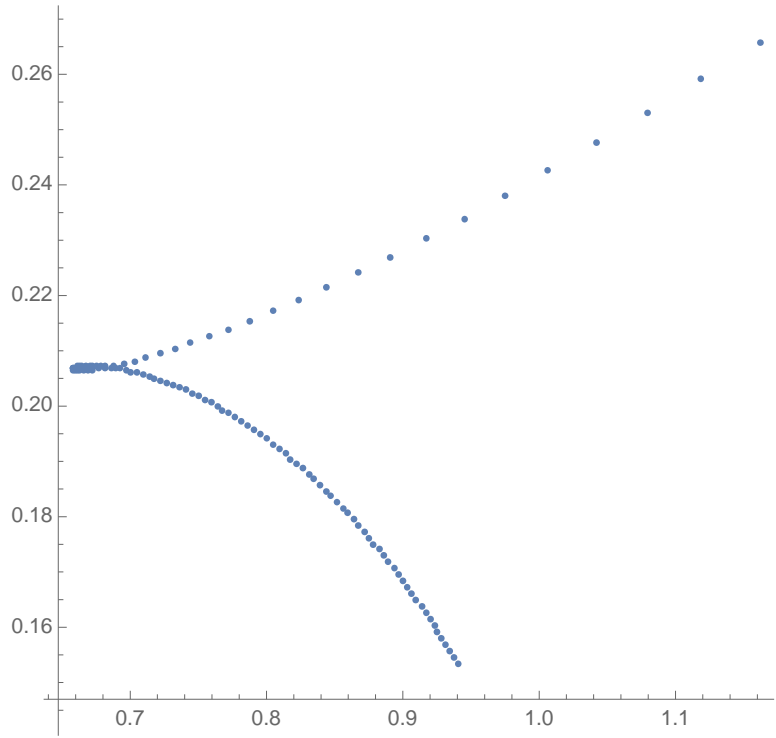

{29}

This means that there are nearby strings not sampled by the original 10 strings that have more pronounced loops. If I use a 10 times finer scale for a string with $t = 357.612$ I get figure {30}

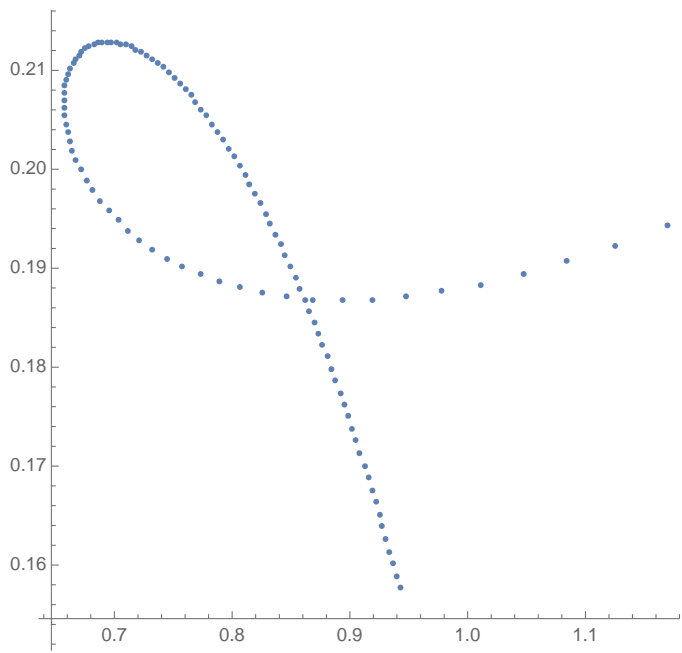

{30}

The apparently horizonal orientation in figure {28} is greatly distorted in figure {30} because of the big difference in scales for the two axes.

The analog of figure {22} for the previous example is figure {31} for this example.
Table[DirichletEta[σ + ti],{t, 357.151, 357.952, 0.089},{ σ, 0, 0.7, .01}]

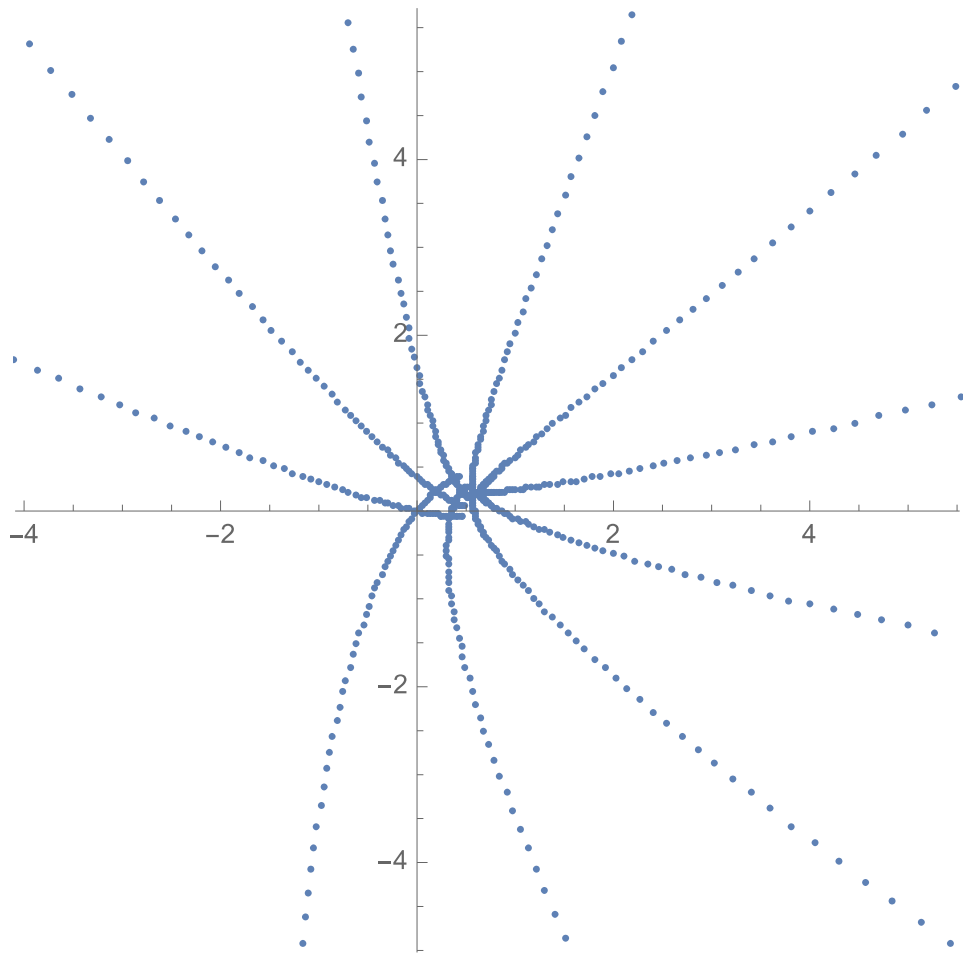

{31}

You can clearly see that the first and the last strings go through the origin as befits the *t*'s of zeros. With these examples, we can begin to see a general feature that is conversion from parallel flare to radial flare occurs between (0, 0) and (1, 0) on the abscissa. Events such as loop formation also occurs in the jumble of conversion somewhere between where strings begin (1, 0) and where they end, radially far from (0, 0). This suggests that only that portion of a string that begins at $\sigma = \infty$ and ends before 0.5 can have a self-crossing. That would mean that a string can intersect the origin, (0, 0), only once (RH).

The final example to be exhibited here is for two of Odlyzko's large $t$ values, 267653395747.6039956773 and 267653395747.8122426610. I can run 8 strings by implementing
Table[DirichletEta[$\sigma$ + $ti$],{t, 267653395747.6039956773, 267653395747.81224266, .0297495691},{ $\sigma$, 4, 7, .02}]. See figure {32}.

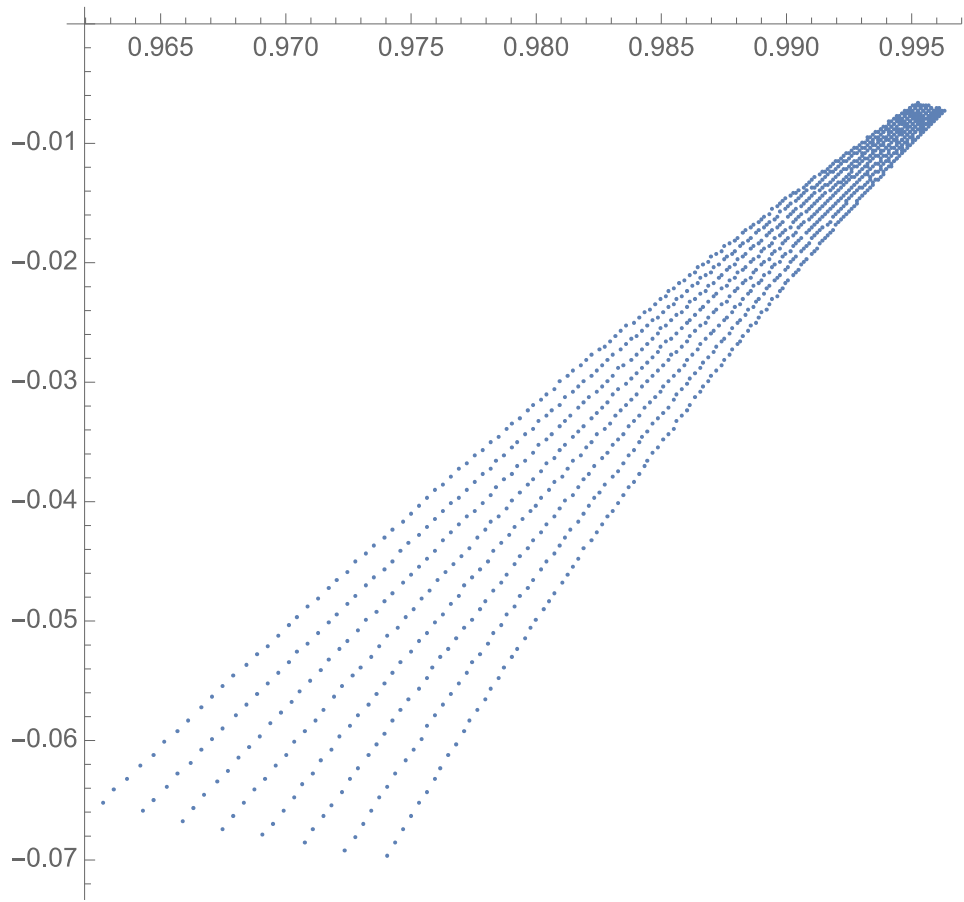

{32}

The next segment is given by
Table[DirichletEta[$\sigma$ + $ti$],{t, 267653395747.6039956773, 267653395747.81224266, .0297495691},{ $\sigma$, 1.5, 4, .01}].

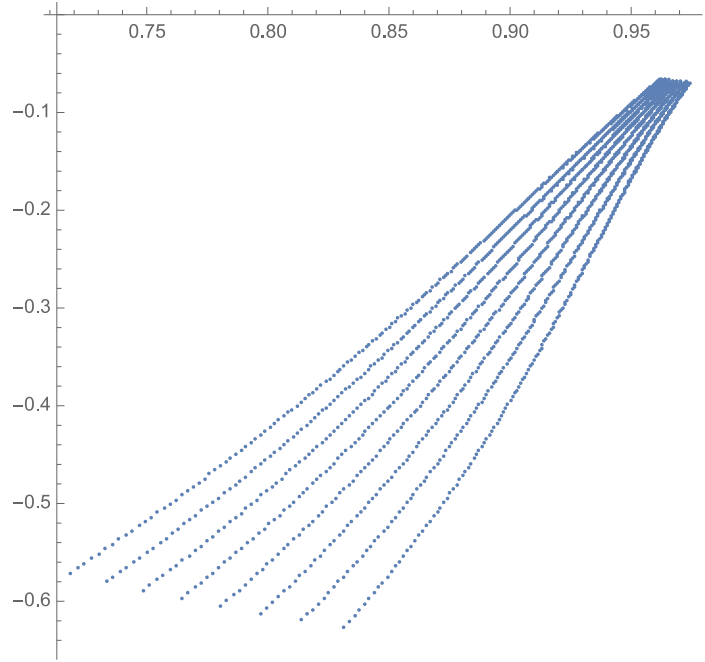

{33}

The last segment is given by
Table[DirichletEta[$\sigma + ti$],{t, 267653395747.6039956773, 267653395747.81224266, .0297495691},{ $\sigma$, .4, 1.5, .01}].

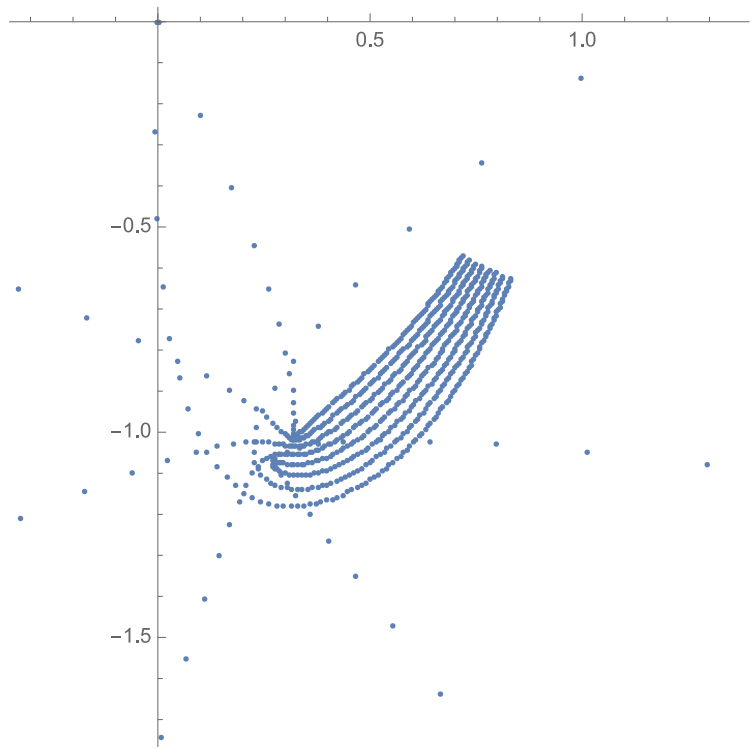

{34}

This last figure can use some refinement. The parameters are changed as in Table[DirichletEta[$\sigma + ti$],{t, 267653395747.6039956773, 267653395747.81224266, .0297495691},{ $\sigma$, .49, 1.2, .005}].

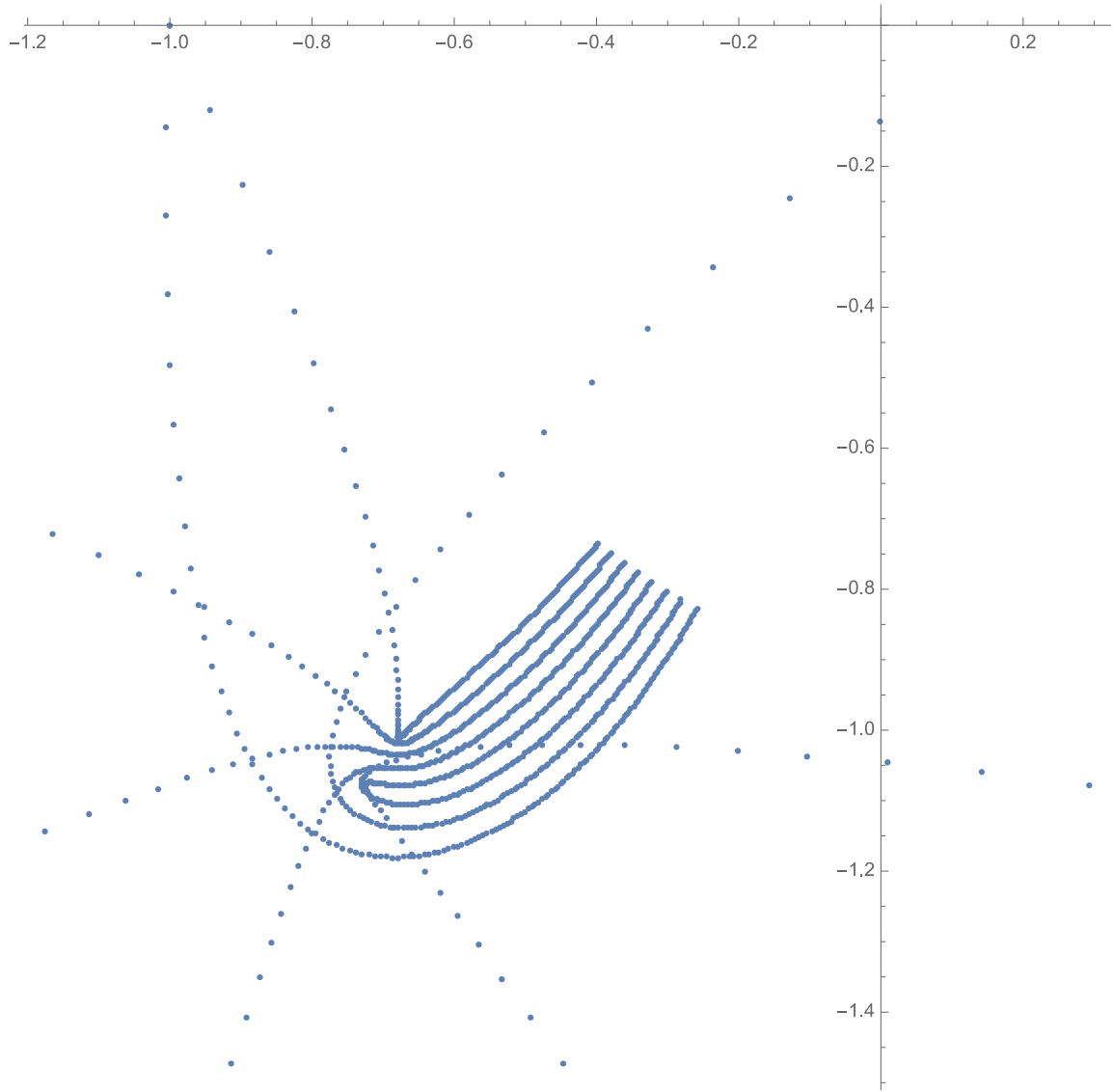

{35}

Note that strings 1 and 8 go through the origin (actually through (-1, 0) because in this case 1 was subtracted in accord with an earlier discussion). Note also that the first two figures were just parallel flares and that it looked like there might not be room enough for conversion to radial flares. Not only did conversion occur but it happened quite far removed from the origin. However, the details of the conversion are obscure and one further magnification is needed. In the figure below, we have used {σ, .57, .8, .002} (there is still a shift in the abscissa by 1.)

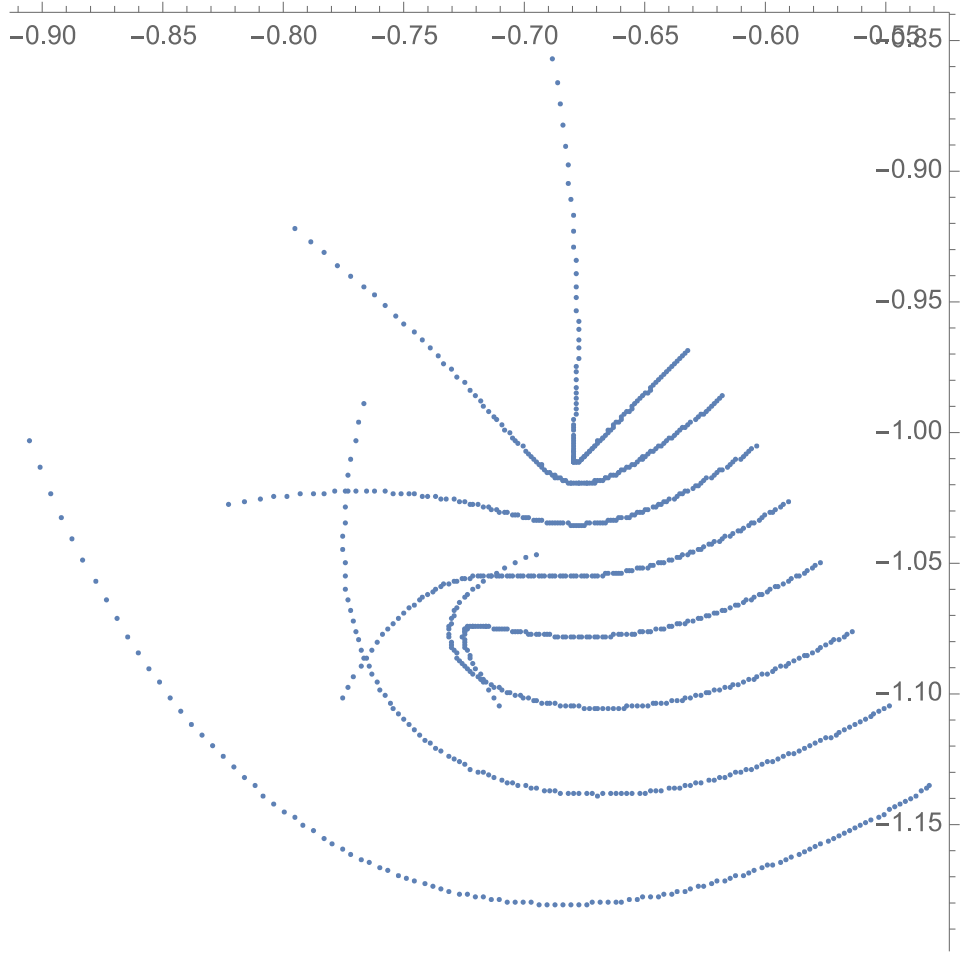

{36}

Comparing this with the figure {35} and numbering the strings from 1 to 8, starting with the topmost one as number 1, then the order of the long extensions starting with the vertical one as number 1 and going clockwise produces the sequence 1, 7, 6, 5, 4, 3, 2, 8. However, strings 1 and 8 meet at (-1, 0), which is (0, 0) without the shift of origin, and then cross over each other giving the reversed order 8, 7, 6, 5, 4, 3, 2, 1. This can be looked at from the perspective of very small σ by running { σ, 0, .7, 0.002} which produces figure {37}. Reading clockwise from one o'clock the order is 8, 7, 6, 5, 4, 3, 2, 1. Notice how long the strings are. In fact, the strings are all cut off short for this figure and one string is actually longer than 382118. Most curious is their straightness.

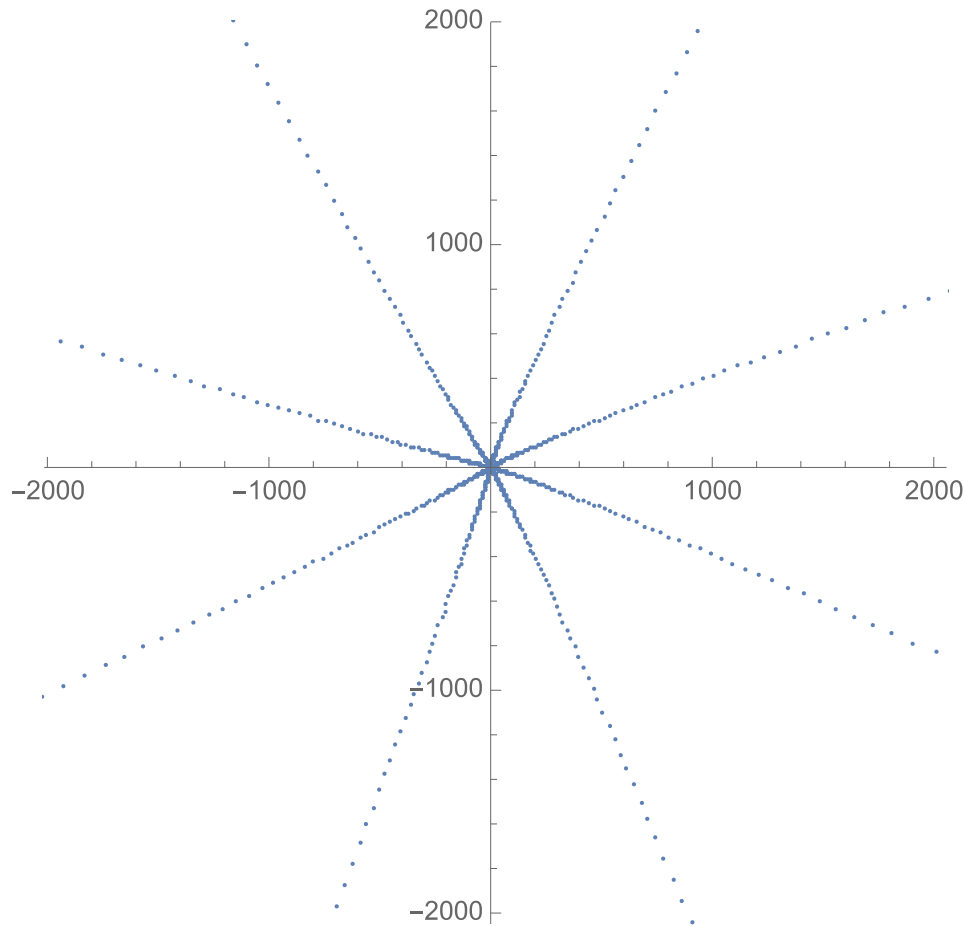

{37}

Why did the parallel flares turn towards the origin? At first, the strings for large $\sigma$ are dominated by the term 1 with small corrections arising from the $n = 2$ term. These are polar forms pointing at (1,0), As $\sigma$ decreases, larger values of $n$ contribute and their contributions are larger so that eventually the sum is dominated not by 1 but by all the other terms. These other terms, other than 1, are polar forms pointing at (0, 0). As $\sigma$ crosses into the critical strip the number of terms other than 1 becomes very large. As $t$ varies the relative amounts of these different polar forms changes, creating the shifts in apparent organizing centers. The approximating sums give some idea of how this develops. For precision $p = 3$, we found that for $\sigma = 12$, $n = 4$; for $\sigma = 3$, $n = 20$; for $\sigma = 1.5$, $n = 200$; and for $\sigma = .5$, $n = 2000000$.

# III. Discussion

Strings provide us with a visual appreciation of the behavior of the Eta and Zeta functions. Starting from small $t$ values (less than 14) strings cannot reach the origin. Thus, there are no zeros. For larger values (greater than 14.135), there can be strings that cross the origin. Every time this happens coincides with a $t$ value for a non-trivial zero. As $t$ grows so does the number of crossings at the origin. The integer parts of $t$ for the first several crossings are: 14, 21, 25, 30, 32, 37, 40, 43, 48, 49, 52…The spacing between consecutive values eventually settles down to the order of one. For example, 4080, 4081, 4081, 4082, 4083, 4084, 4086… For much larger values [6], such as 267653395648, the integer part of 267653395648.8475231278, the crossing values are 4 per integer. Instead of writing the entire number I will write 648 as shorthand for 267653395648. One gets: 648, 649,649, 649, 650, 650, 650, 650, 651, 651, 651, 651,.. Plots of strings show these crossings happening for the string position corresponding to $\sigma = 0.5$. The justification for this was presented earlier and is equivalent to the statement that *a line and a point can intersect only once*. If so, then the modified reflection formula forces $\sigma = 0.5$. The "if so" is the hitch. If a string can cross itself then such a crossing happening at the origin would invalidate this line of reasoning and instead *prove* not-RH. Even though we found a crossing very close to the origin we did not find one precisely on the origin. Instead, we found that strings change, as $\sigma$ decreases, in a generic way.

Having fixed $t$, we can view the string as starting with large $\sigma$ and developing as $\sigma$ decreases all the way down to zero. Large $\sigma$ means $\sigma > 20$. The value of Eta is a small correction to 1 by an amount given by modulus $1/2^{20}$ or smaller, with the next term's modulus of order $1/3^{20}$ or smaller. Figure {10} is an example. Even for several consecutive values of $\sigma$, the angle dependent factor really does factor producing a straight string. This holds for values of $t$ that take you completely around the operative center, $(1, 0)$. If values of $t$ are chosen that create a group of closely spaced strings, such as in figures {15, 23, 32}, then a parallel flare is created. For figure {15}, the initial direction of the parallel flare is towards one o'clock, for figure {23}, it is towards 3 o'clock and for figure {32} it is towards 7 o'clock. As $\sigma$ decreases more values of $n$ are important and the values of the points on the string get larger. Eventually, these additional terms are numerous enough and large enough to overwhelm the beginning value, 1.

This means that all the terms added up to make Eta are polar forms centered on (0, 0). In a second phase, the parallel flare turns towards the origin and continues in the new direction as $\sigma$ gets even smaller. This is seen in figures {16, 25} but in figure {33} the initial behavior is continued in the same direction during the second phase. At first the growth in the number of $n$ terms needed is slow but eventually it grows exponentially. The discussion of precision p captures the essence of these events. Once the sudden increase in the number of contributing $n$'s takes off the strings engage in a dramatic jumble (there is not enough space here to elaborate on this highly technical terminology) displayed in figures {17, 21, 27, 34, 35}. Out of this jumble a radial flare develops, usually ending, as $\sigma$ goes to zero, in a circular pattern with very straight and very long radii, as in figure {37}. Wherever the jumble occurs, which is clearly to the left of (1, 0) determines the particular location and nature of the organizing principle of a center. The center is usually also to the right of (0, 0). It may be displaced up or down the ordinate.

The most striking feature of the jumble is the occurrence of loops in some strings. These are self-crossings of a string. In figures {28, 29} there is loop formed during a string's U-turn. In figure {5} the loop is produced on a string that also crosses the origin. Although this loop occurs very close to the origin relative to the size of the entire string (0.002 versus 15000) it does not have its self-crossing at the origin. In fact, this loop's crossing point is just outside of the critical strip. If it can be proved that all crossings occur inside the jumbles and that the jumbles are over before the loops reach (0, 0) then it would *not* be possible for a string's self-crossing to coincide with the origin. This would make the assertion that a string's self-crossing *never* coincides with the origin true and, therefore, RH is true.

One final avenue of approach is directed at the structure for small $\sigma$. We see a radial flare of long straight strings. Very close to the origin these strings generally deviate from their overall straightness and can even turn through 270 degrees. That doesn't allow them to invalidate RH. Can a definitive statement be made and proved regarding this structural constraint? The same observations that have been made about how loops form in jumbles suggest that a string never has more than one loop. That may help in establishing the possibilities for string self-crossings that simultaneously occur at the origin.

Return to figure {5} for one last look at the possibility of not-RH. The figure depicts a string for a zero with $t = 267653395649.3623$. It clearly goes through the origin with a self-crossing nearby. In a greater context it is the 10th string in the set-up:

Table[DirichletEta[$\sigma + ti$], {$t$, 267653395648.8475, 267653395649.3623, .0572}, { $\sigma$, .48, 1.4, .002}]. The first string in this command is also for a zero. The result is figure {38}.

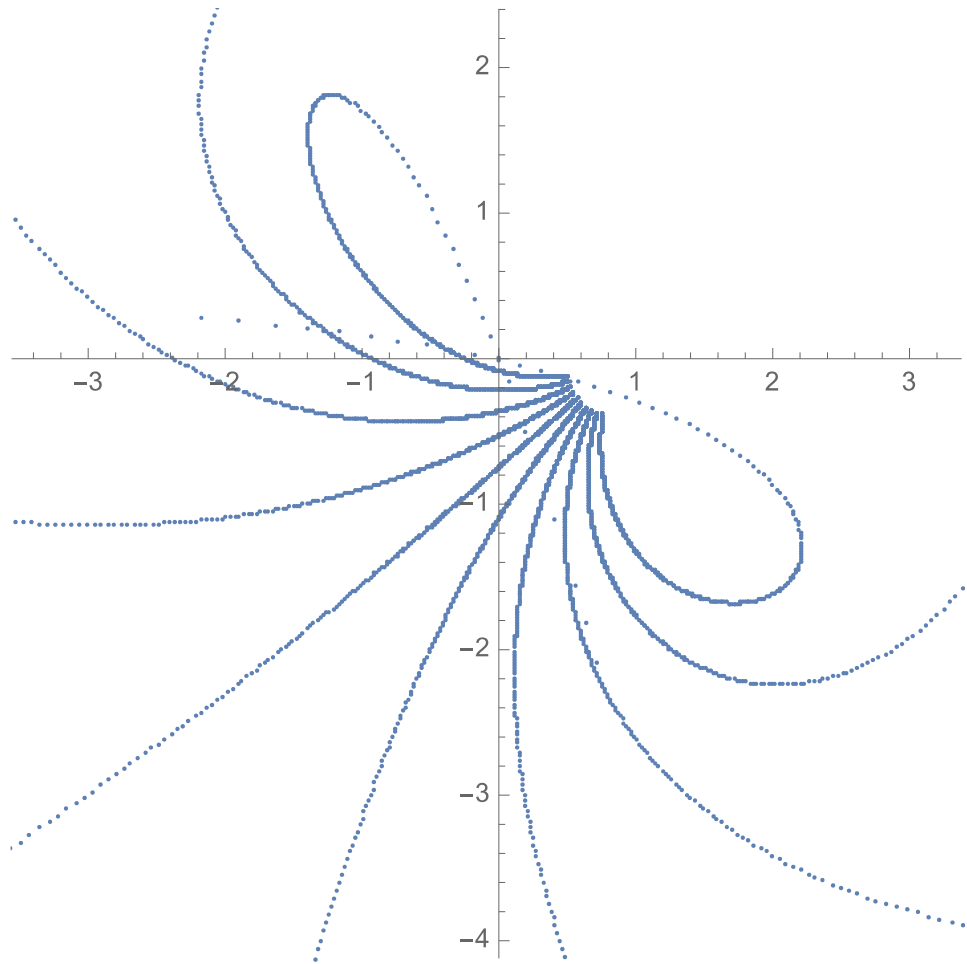

{38}

The loop in figure {5} is the upper loop here, and opposite it is another loop for the first string. The upper loop clearly crosses itself. One can see that this figure would attach to one closer to (1, 0), i.e. one with $\sigma$ greater than 1.4. Then the first loop would also cross itself, but not as close to (0, 0) as the upper loop. The zeros are discrete so you cannot simply continuously tune this figure until the crossing is at the origin.

# Appendix

In this appendix a variety of $t$ strings formed for $t$ in the range 1 to 67.6, with various spacing's from 1 to 0.1 are exhibited in order to justify the rules presented above. The $\sigma - t$ ranges label the figures. This is followed by a discussion of the figure contents. Then the figure is presented.

[{$\sigma$, 0, 1, .05}, {$t$, 1, 11, 1}]

There exists an imaginary part of a *trivial* zero at $t = 9.0647\ldots = \frac{2\pi}{\ln(2)}$.

This is the $k = 1$ case. The rightmost sigma dot ($\sigma = 1$) of the ninth string counting clockwise from the shortest string almost crosses the origin. Increasing $t$ from 9 to 9.0647… will make that happen. The maximum $t$ is 11 which is not large enough to produce a *nontrivial* zero whose imaginary part must be larger than 14. Proceeding clockwise from $t = 1$ to $t = 11$ illustrates the points made above about string location, length, and dot density. The seventh string is the first to regain at least unit arc length.

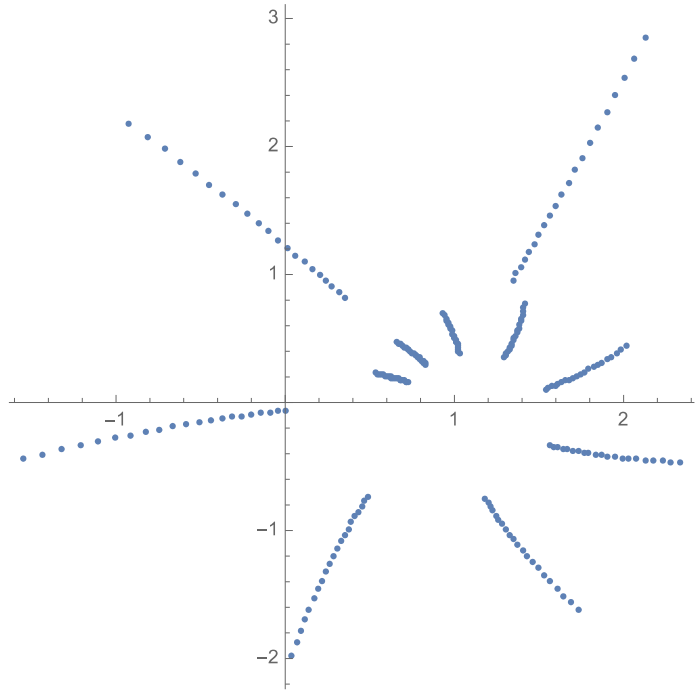

[{σ, 0, 1, .05}, {t, 11, 15, .5}]

There is an imaginary part of a zero for $t$ between 14 and 14.5. The eleventh sigma ($\sigma = 0.5$) almost crosses the origin. Compared to the figure above, the spacings of the $t$ values have changed from 1 to 0.5. The $t = 11$ string was the last string in the figure above and is the first string in this figure. The strings with $t = 14$ and $t = 14.5$ are near 9 o'clock. A string between them will cross the origin for $\sigma = 0.5$. The value for this string is $t = 14.134725$ ...

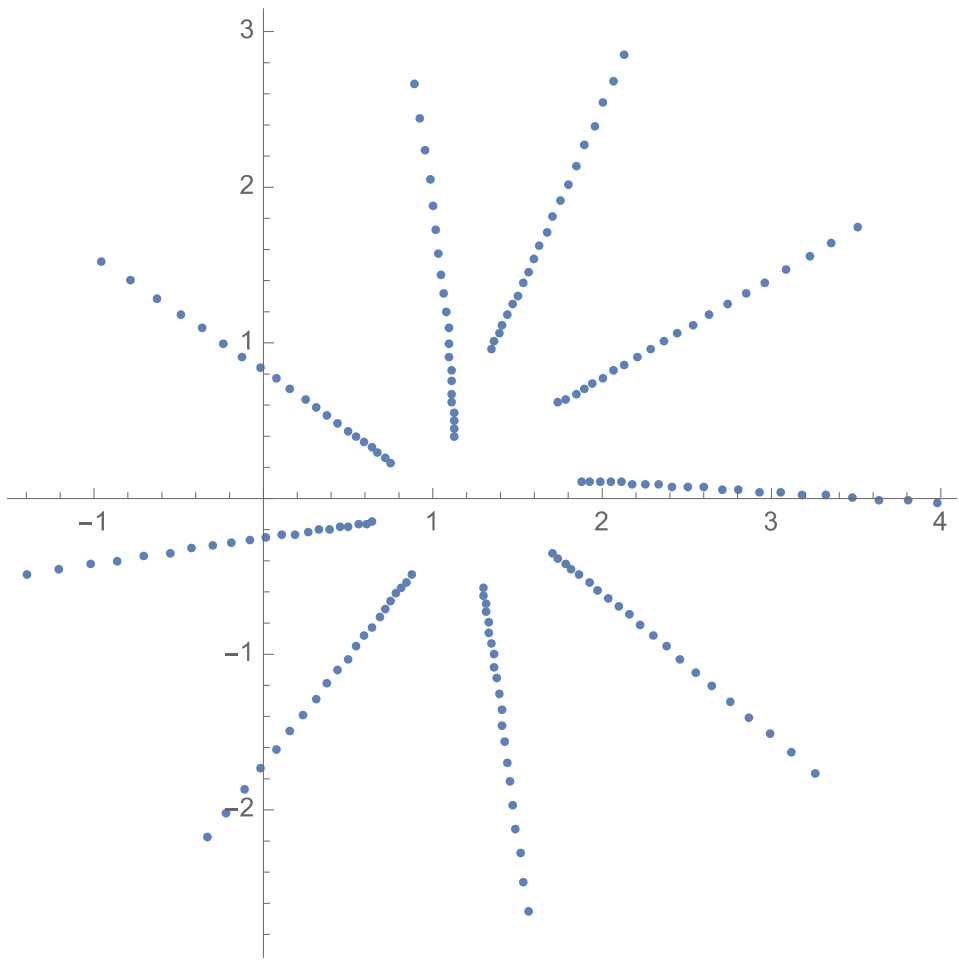

[{σ, 0, 1, .05}, {t, 15, 18, .5}]
Another trivial zero occurs for $t = 18.1294 \ldots$ ($\sigma = 1$). This is the $k = 2$ case. The last string in this figure is for $t = 18$. If this value of $t$ is increased by 0.1294… then the $\sigma = 1$ end will sit on the origin.

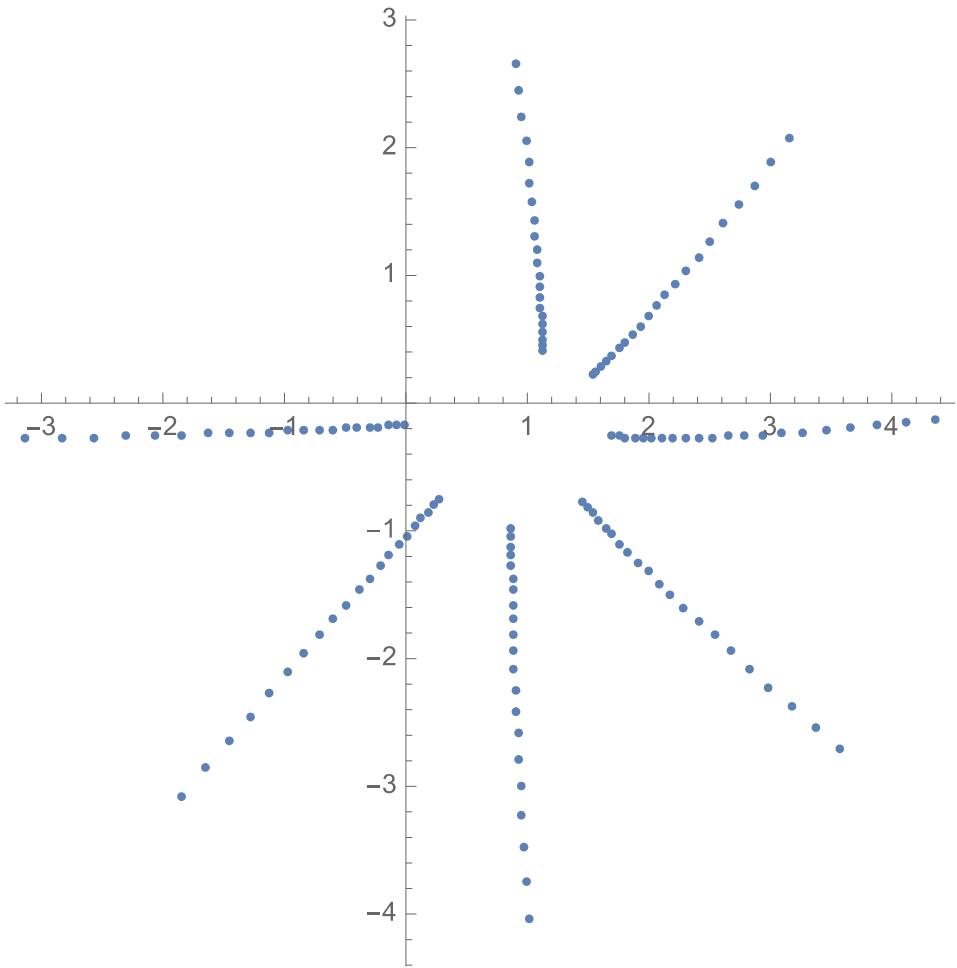

[{σ, 0, 1, .05}, {t, 19, 21, .2}]
There is a non-trivial zero for t = 21.0220 ... The eleventh sigma dot of the last string (t = 21) almost crosses the origin (σ = 0.5). A slight increase and it will do so.

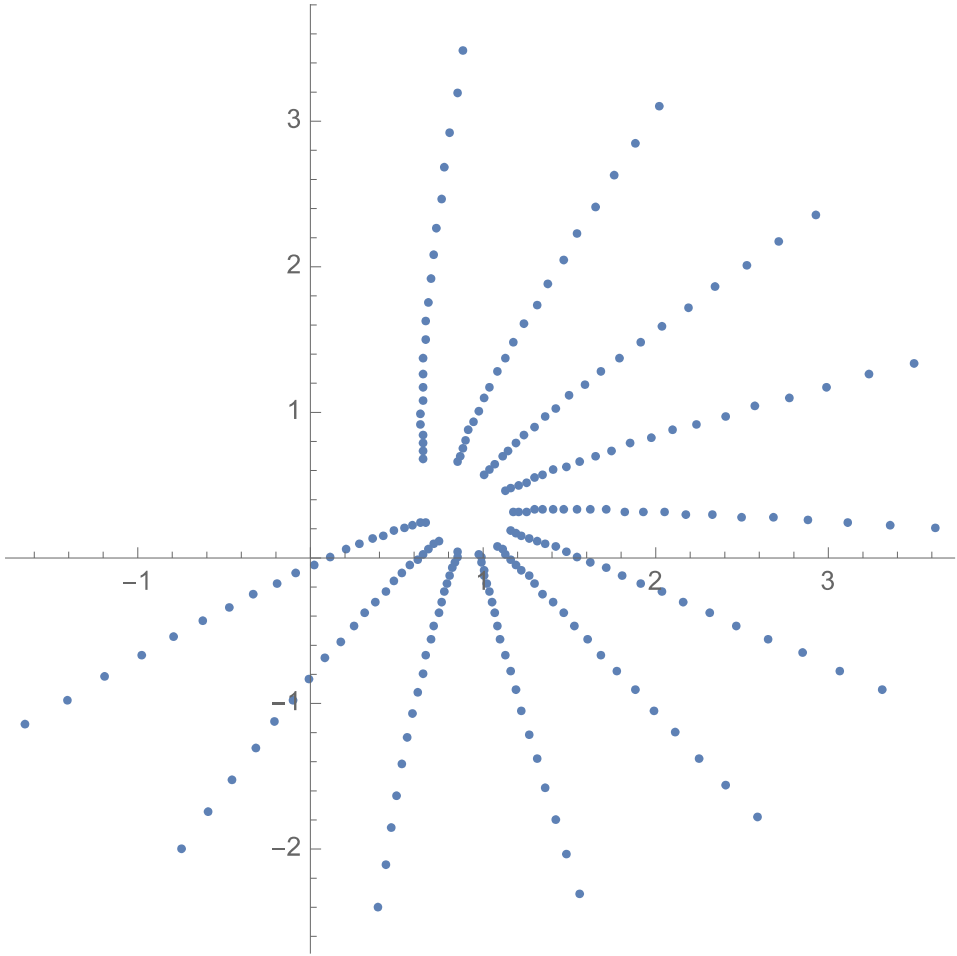

[{σ, 0, 1, .05}, {t, 21, 23, .2}]
There is the same zero as above near $t = 21$. It appears evident that the string's $\sigma = 1$ ends are focused on a center of rotation although this center appears to move around as well. Why all strings appear as nearly straight rays for $\sigma < 0.5$ remains to be understood. The $t$ spacing is now $0.2$.

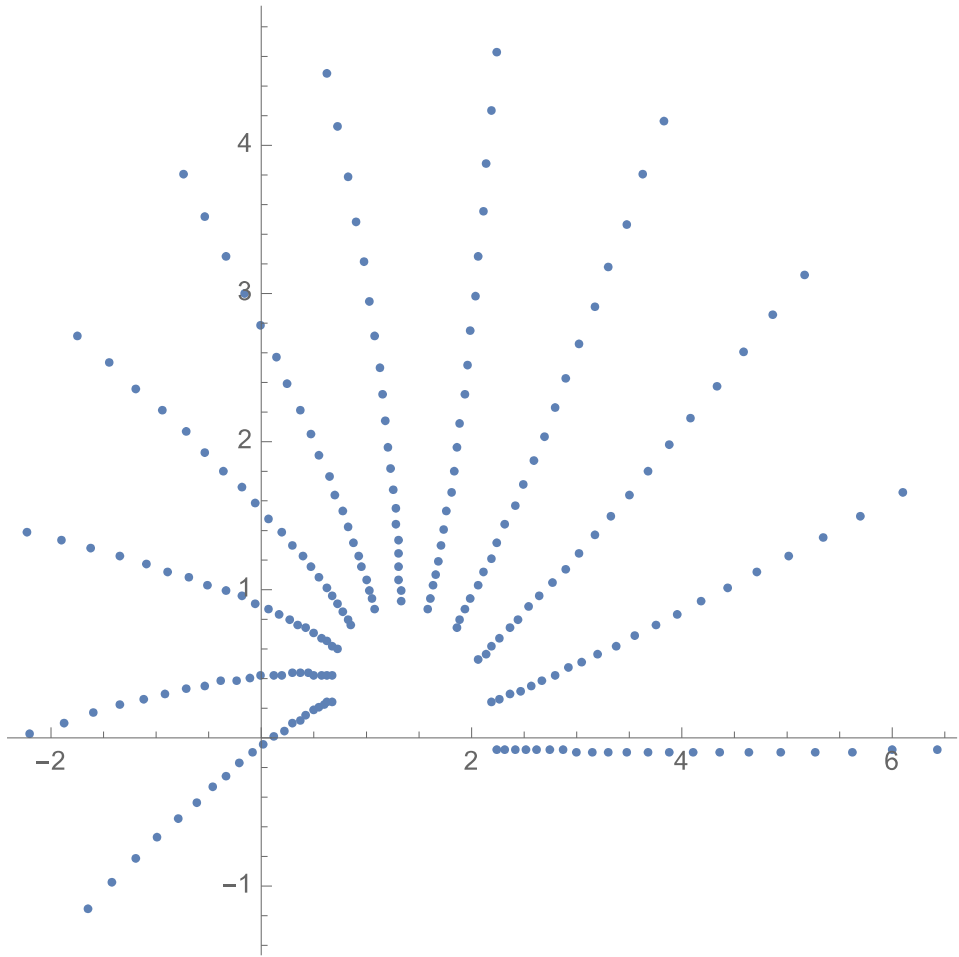

[{σ, 0, 1, .05}, {t, 21, 24, .2}]

There are no zeros in this t range. Smaller increments in t fill up one rotation of a circle with a given number of strings. Mysterious centers of rotation are evident. For example, in the first figure, for t running from 1 to 11 separated by ones, 1.33 rotations of a circle were produced, whereas in the figure above 11 strings separated by 0.2 produced 5/8 of a circle. The area of the first figure is roughly 4 × 5 and the area of this figure is roughly 6 × 8. This increase reflects the increase in string lengths.

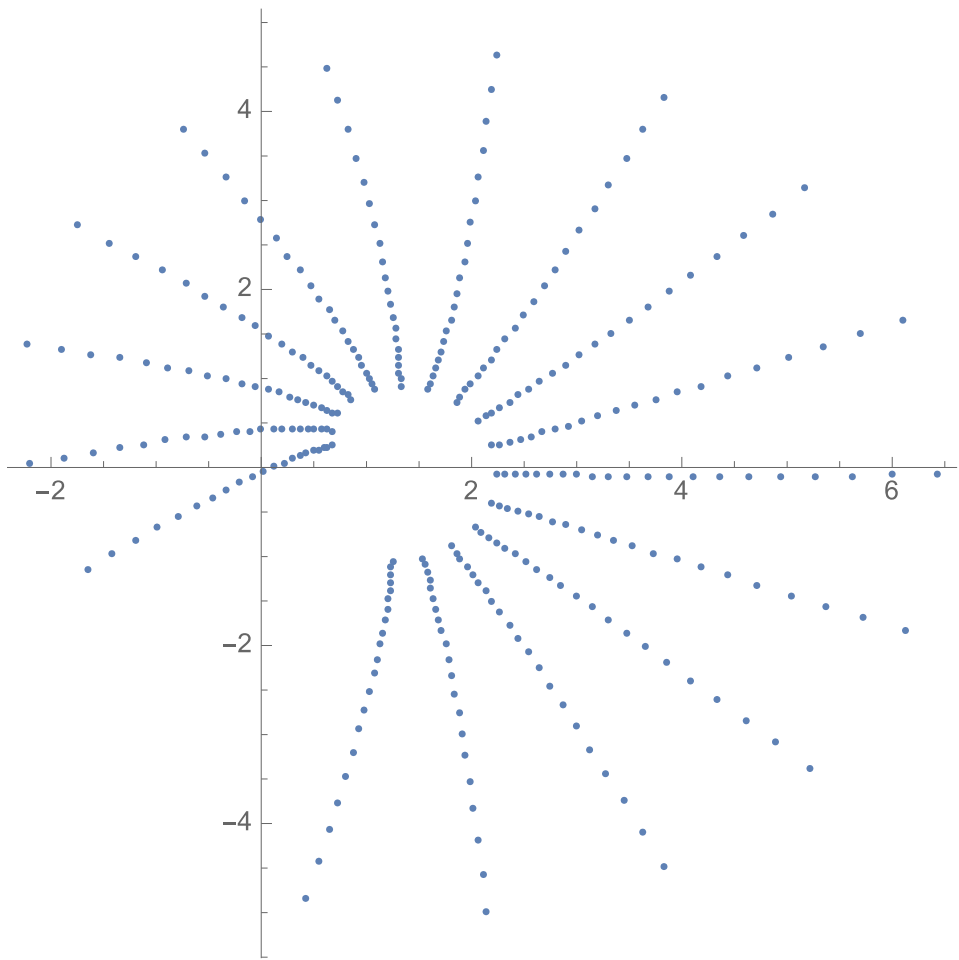

[{σ, 0, 1, .05}, {t, 24, 26, .2}]
There is a nontrivial zero for $t = 25.010$ ... The eleventh dot of the center string ($t = 25$) is almost right on the origin. The center of rotation is contracting and moving. The $t = 26$ string is shorter than the $t = 24$ string. In the long run length increases with increasing $t$.

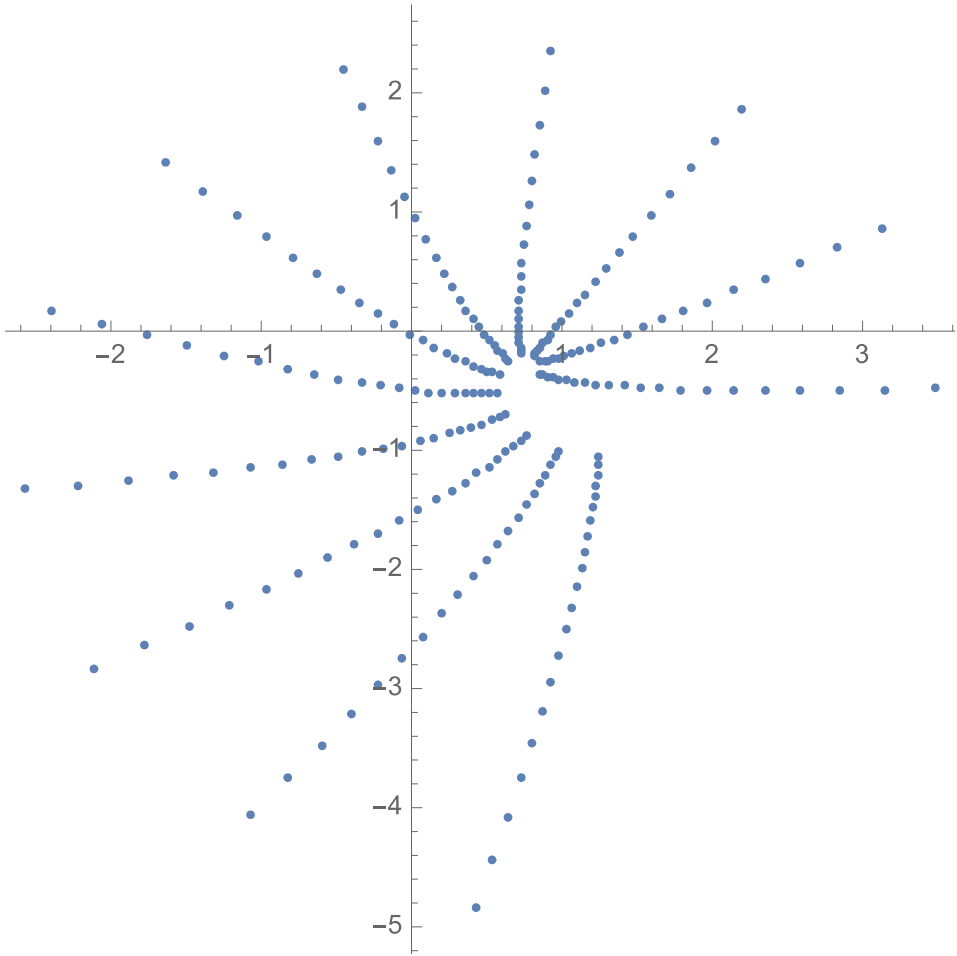

[{σ, 0, 1, .05}, {t, 26, 28, .2}]
A trivial zero for $t = 27.1941 \ldots$ ($k = 3$) is apparent as the $\sigma = 1$ dot of the seventh string ($t = 27.2$) from the lower right crosses the origin. The strings are getting longer again.

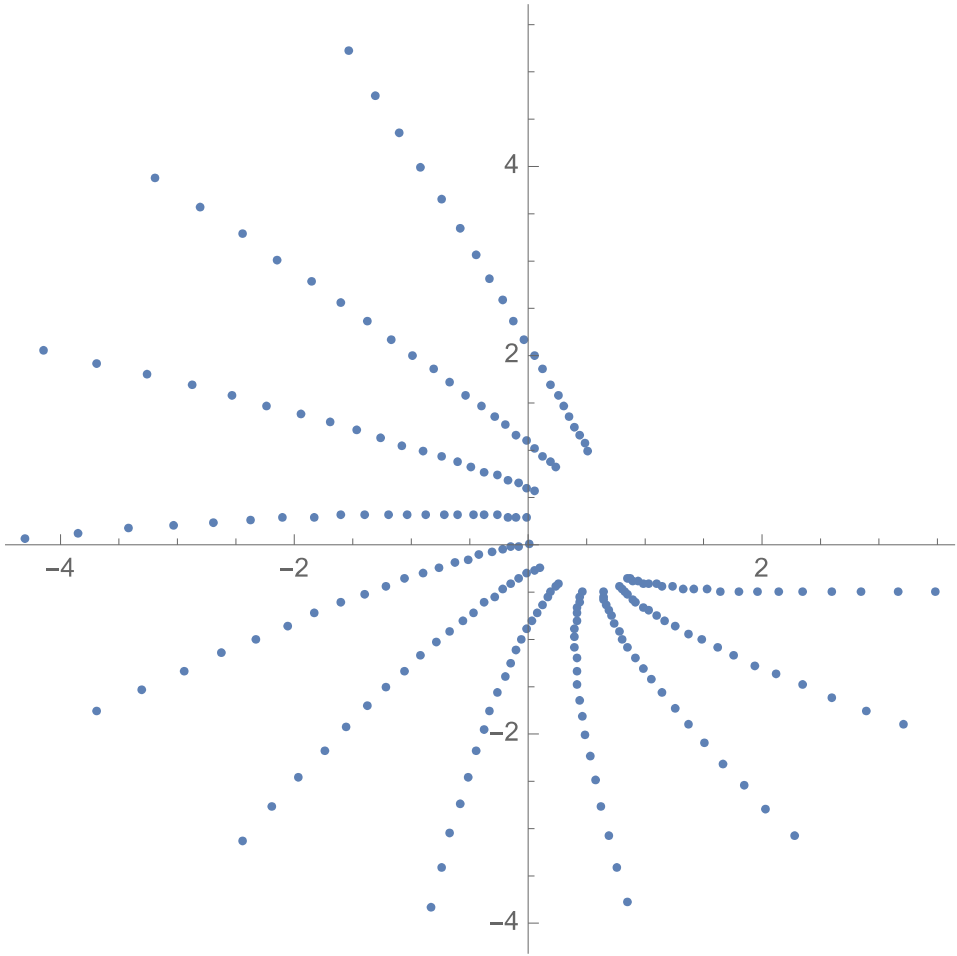

[{σ, 0, 1, .05}, {t, 29, 31, .2}]
There is a nontrivial zero at $t = 30.424$ ... just above the eleventh dot ($\sigma = 0.5$) of the eighth string, counting clockwise from 3 o'clock. The nonuniformity of the distribution of dots caused by the action of Eta is increasing. The dense end is always the $\sigma = 1$ end.

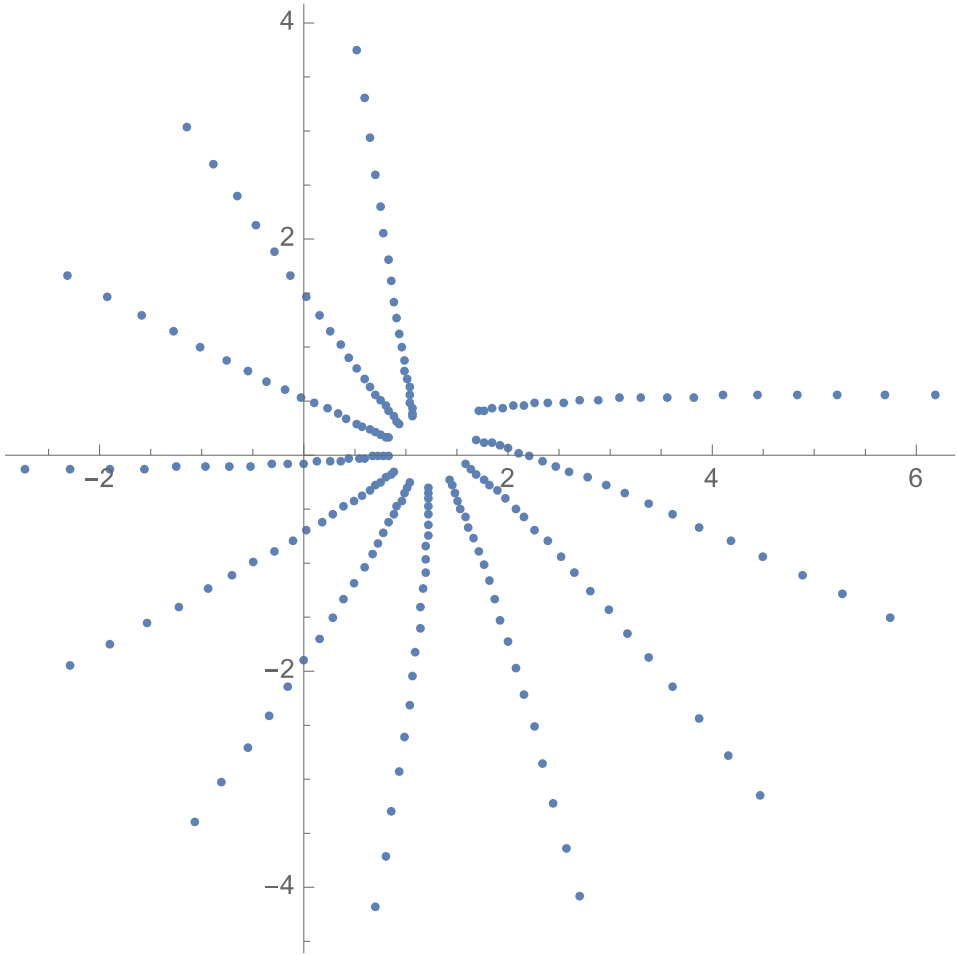

[{σ, 0, 1, .05}, {t, 31, 33, .2}]
There is a nontrivial zero at $t = 32.935$ ... This is just below the eleventh dot of the $t = 33$ string (at 9 o'clock).

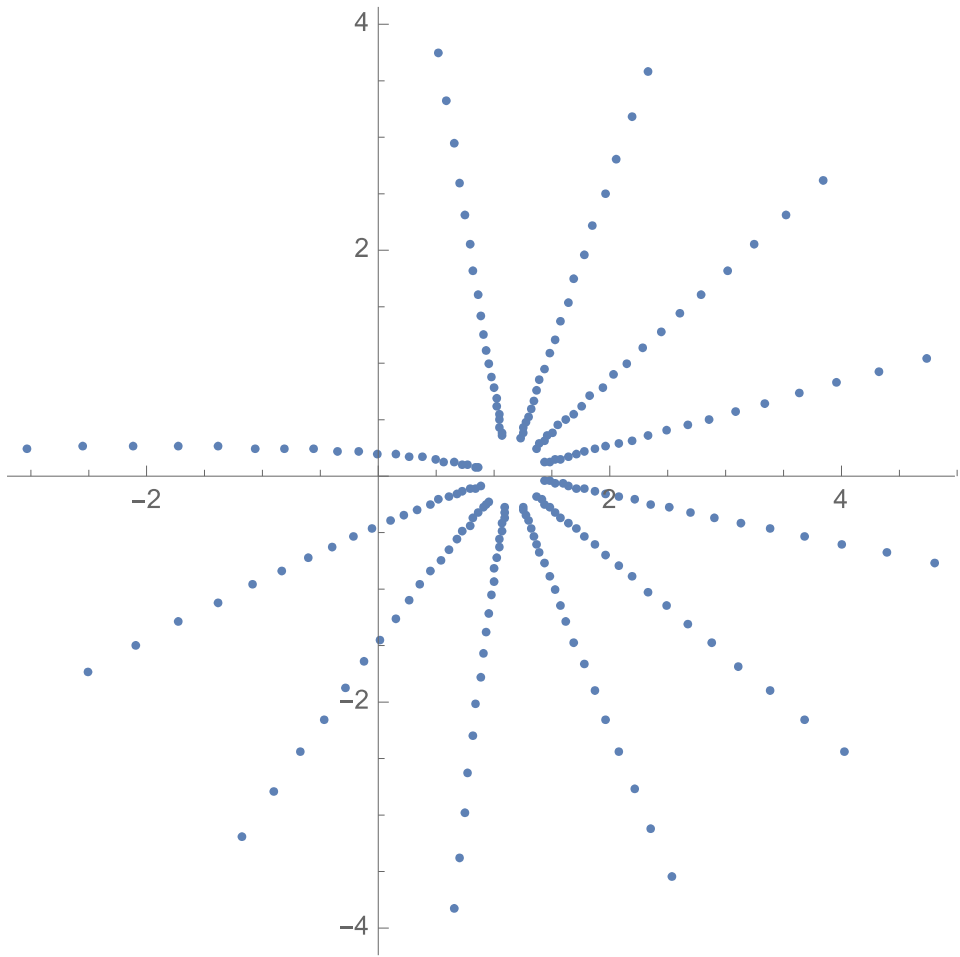

[{σ, 0, 1, .05}, {t, 35, 37, .2}]
There is a trivial zero at t = 36.2588 ... (k = 4). This is very near the σ =1 end of the seventh string counting clockwise from 5 o'clock. Notice that the curvature of the σ = 1 ends of strings 36.8 and 37 are greater than in any earlier figure.

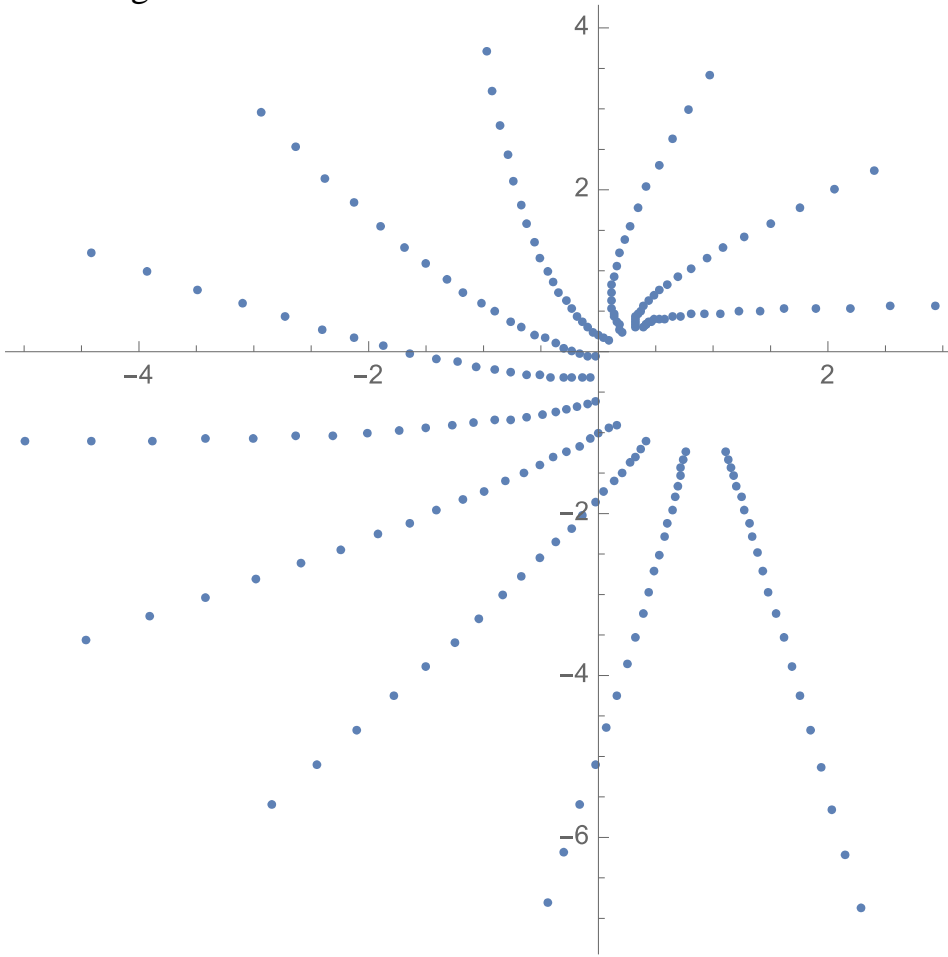

[{σ, .6, 1, .05}, {t, 35, 37, .2}]]
This is a magnification of the σ = 1 end of the strings in the figure above for σ greater or equal to 0.6.

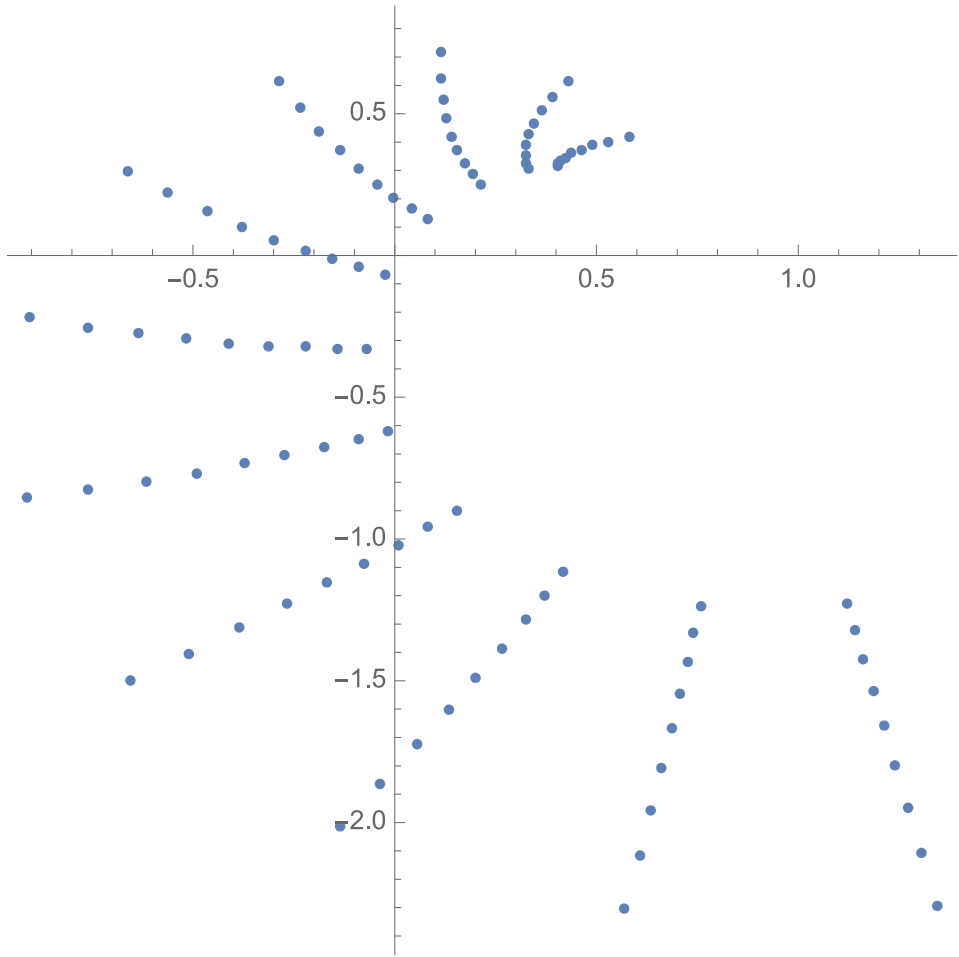

[{σ, 0, 1, .05}, {t, 37, 39, .2}]
There is a nontrivial zero at $t = 37.586$ ... Counting clockwise from the string at 3 o'clock, the fourth string is for $t = 37.6$ and the eleventh dot of this string is right at the origin. Something dramatic seems to be happening at the $\sigma = 1$ ends in this region of the strings.

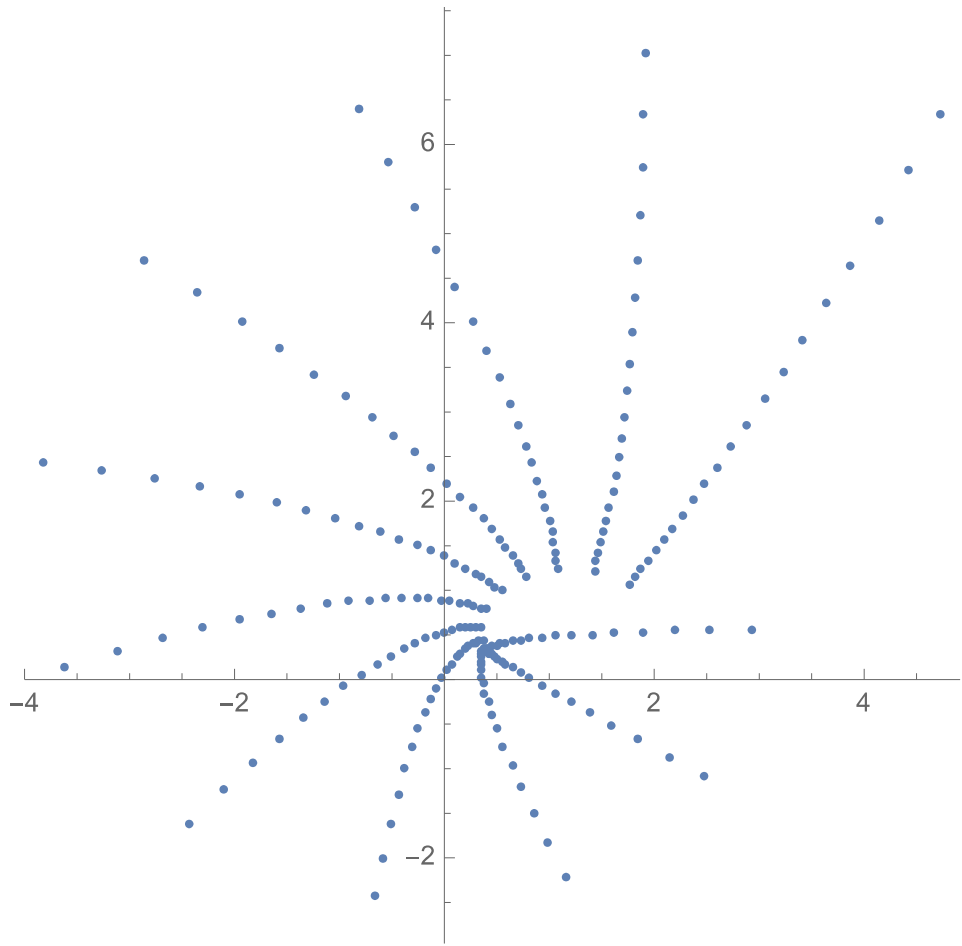

[{σ, .6, 1, .05}, {t, 37, 39, .2}]]
Cutting out the small values of σ enlarges this region.

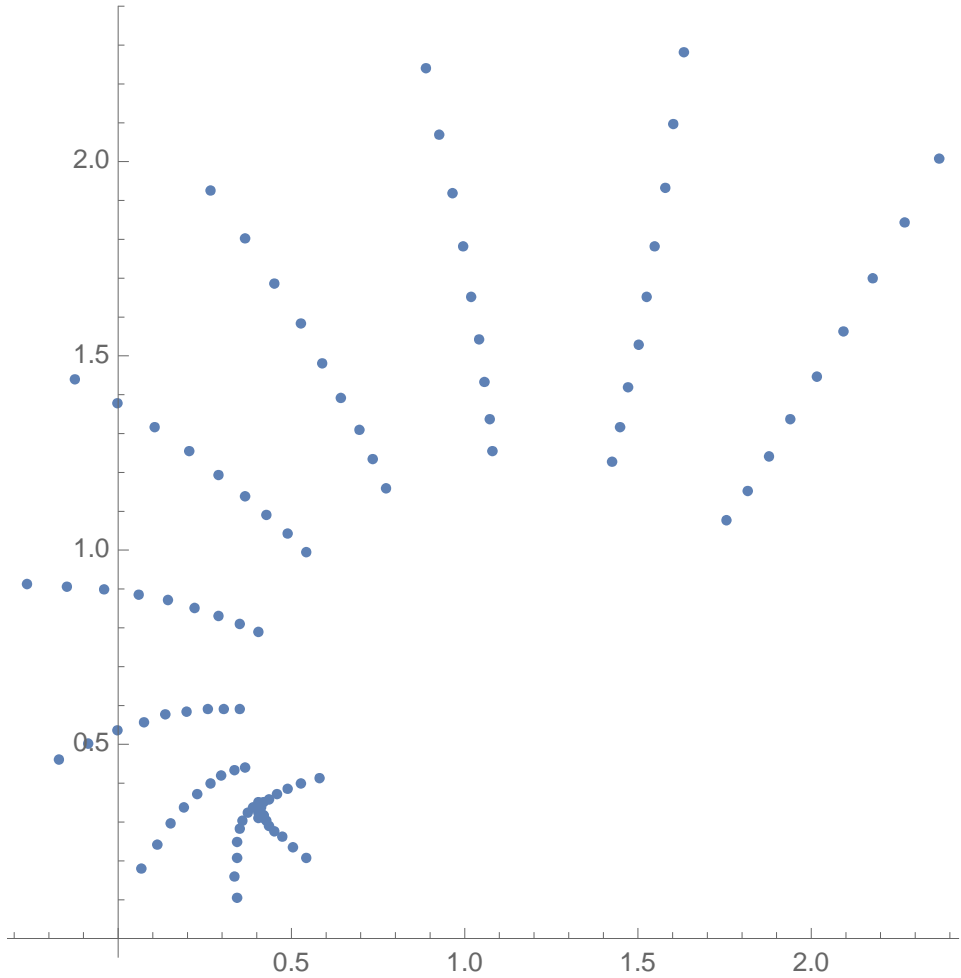

[{σ, .6, 1, .05}, {t, 37, 38, .2}]

A strange configuration appears at this magnification. By reducing the maximum $t$ from 39 to 38 more magnification is achieved. Three string ends are seen together. Their curvatures are a little bit bigger than for other strings and their ends are very close together. The first two strings have opposite curvatures.

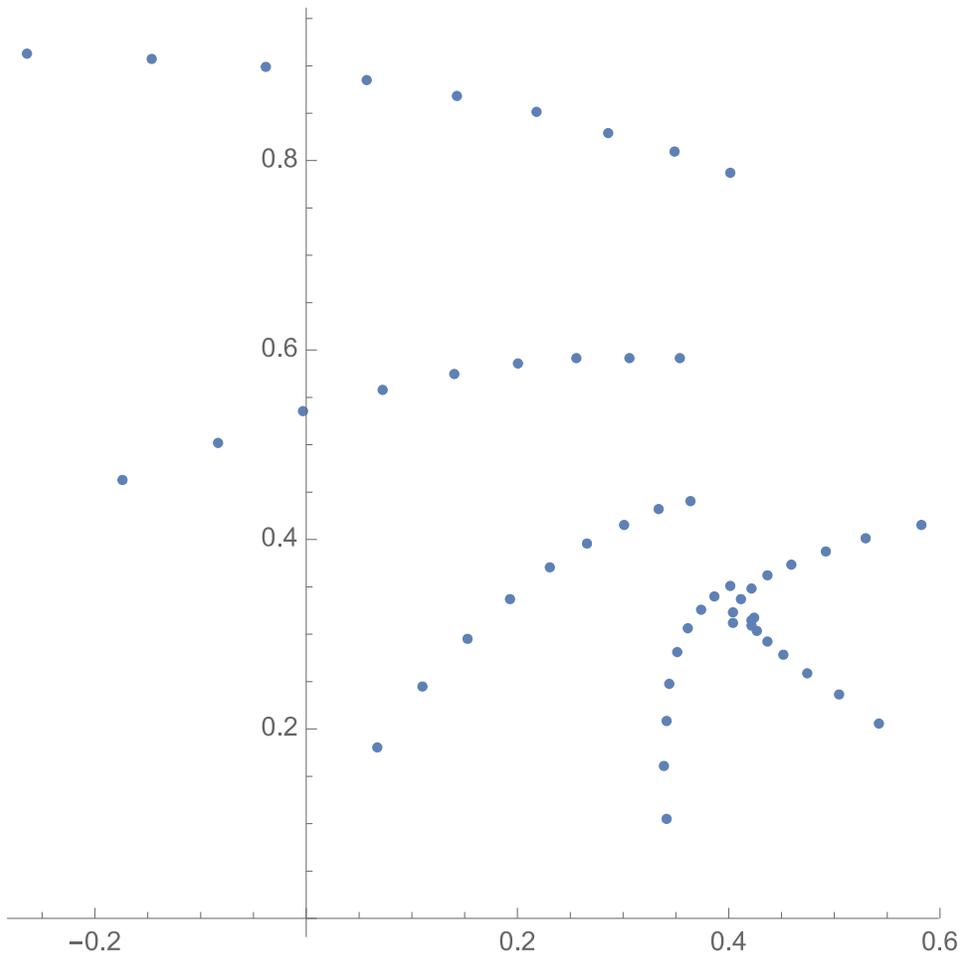

[{σ, .5, 1, .05}, {t, 37, 38, .2}]
There is a nontrivial zero in this range with $t = 37.586$ ... Having cut $\sigma$ off at 0.6 above eliminated the origin crossing for this zero. By adding $\sigma = 0.5$ back in the crossing can be seen close to the origin after Eta acts.

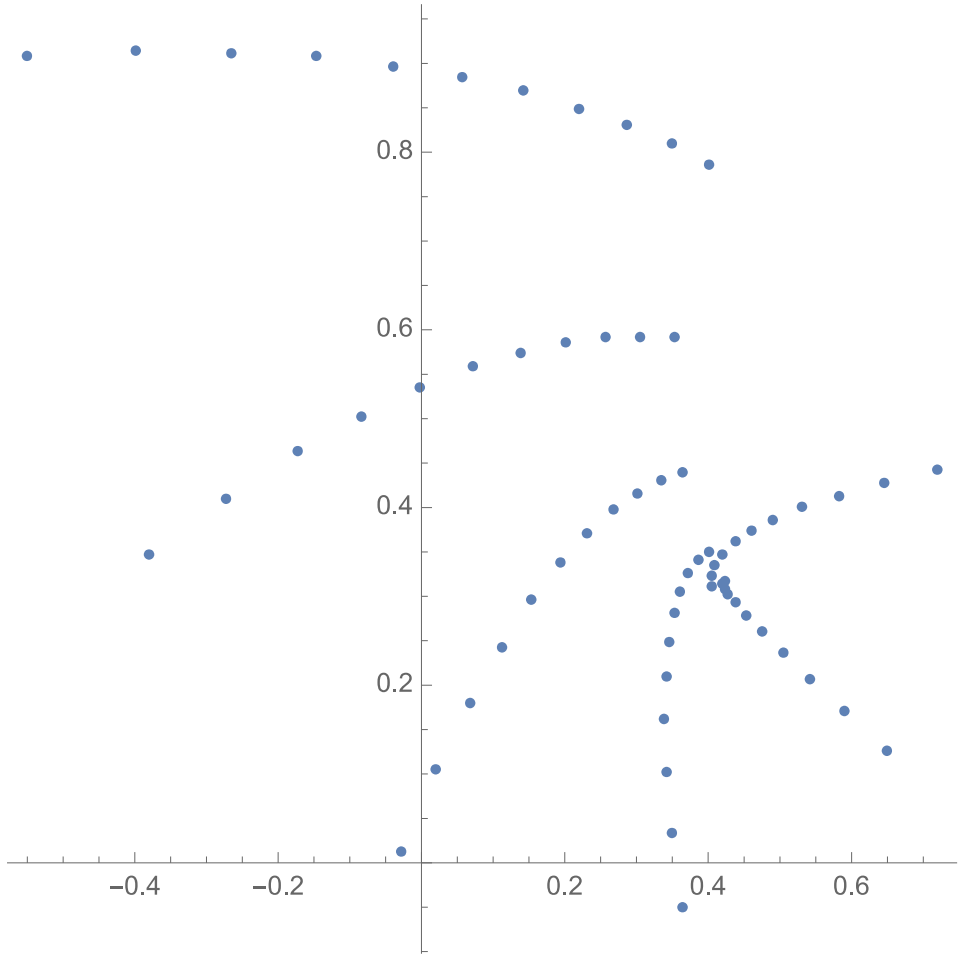

[{σ, 0, 1, .05}, {t, 39, 41, .2}]

In this range there is the nontrivial zero with $t = 40.918$ ... This is about halfway between the $t$ equals 40.8 and 41 strings. Indeed, the eleventh dot ($\sigma = 0.5$) for $t = 40.8$ is just to the left of the ordinate and the eleventh dot for $t = 41$ is just to the right of the ordinate.

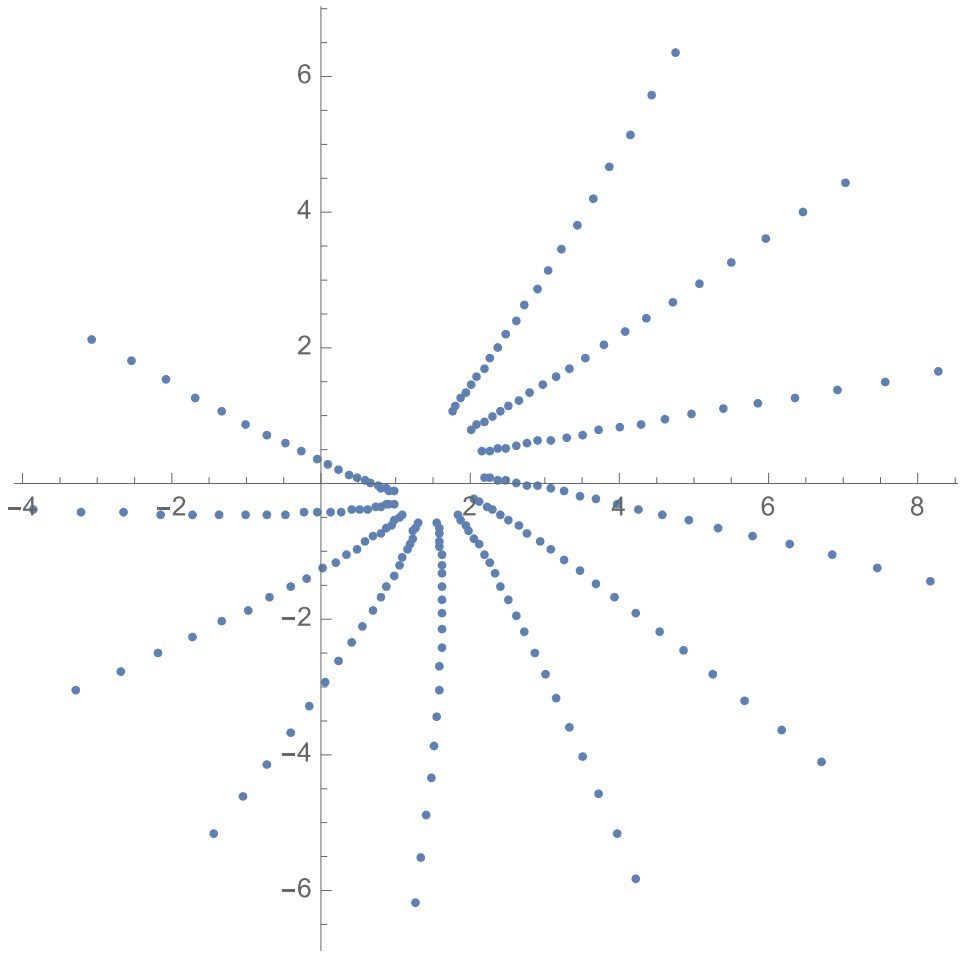

[{σ, 0, 1, .05}, {t, 43, 44.4, .2}]

There is a nontrivial zero for $t = 43.327$ ... The tenth dot from the left end of string 43.4 is on the ordinate just above the origin. If Eta were evaluated for 43.327 ... instead of for 43.4 then the eleventh dot would be at the origin as required. Note that $t$ is increasing less and less as the strings go around a circle. The strings are longer and the lack of uniformity in dot density is becoming more extreme. The centers of rotation continue to move around.

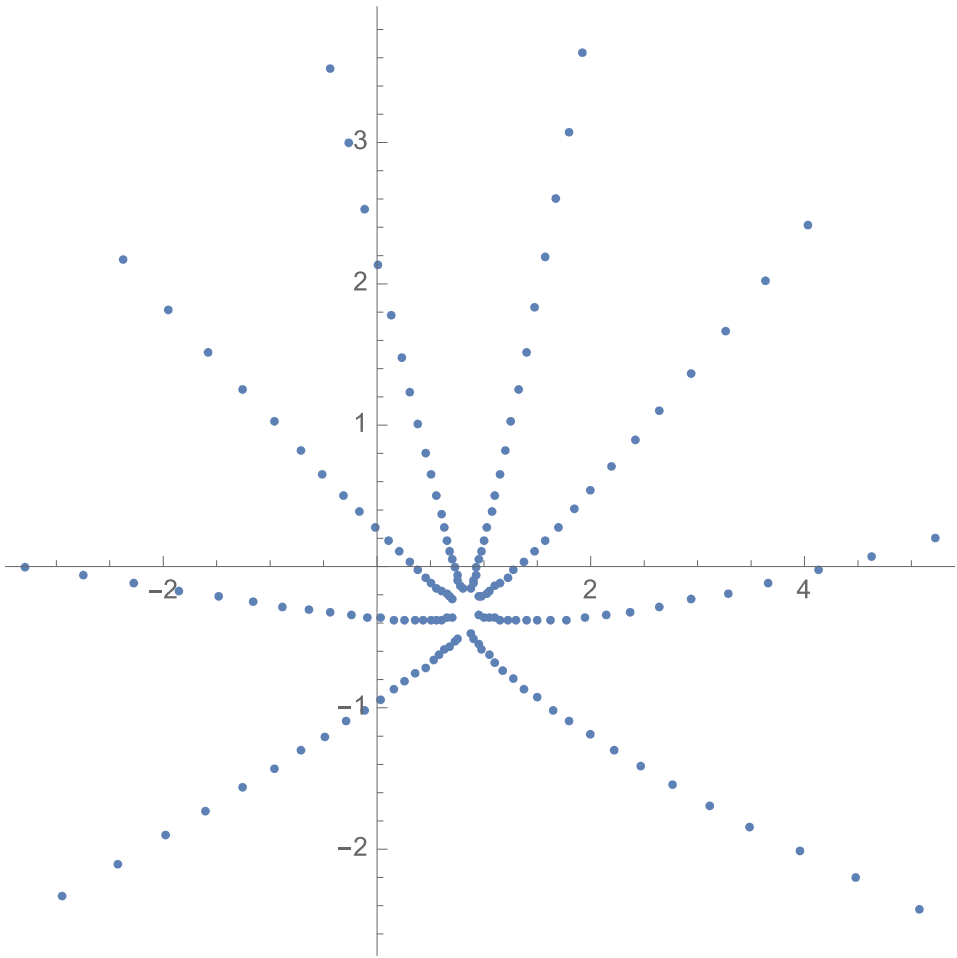

[{σ, 0, 1, .05}, {t, 44.4, 46, .2}]
Only a single trivial zero is in this range. It has $t = 45.3235 \ldots$ (k = 5). Interpolating this value into the figure between 45.2 and 45.4 it is clear that the $\sigma = 1$ end of the string will be at the origin.

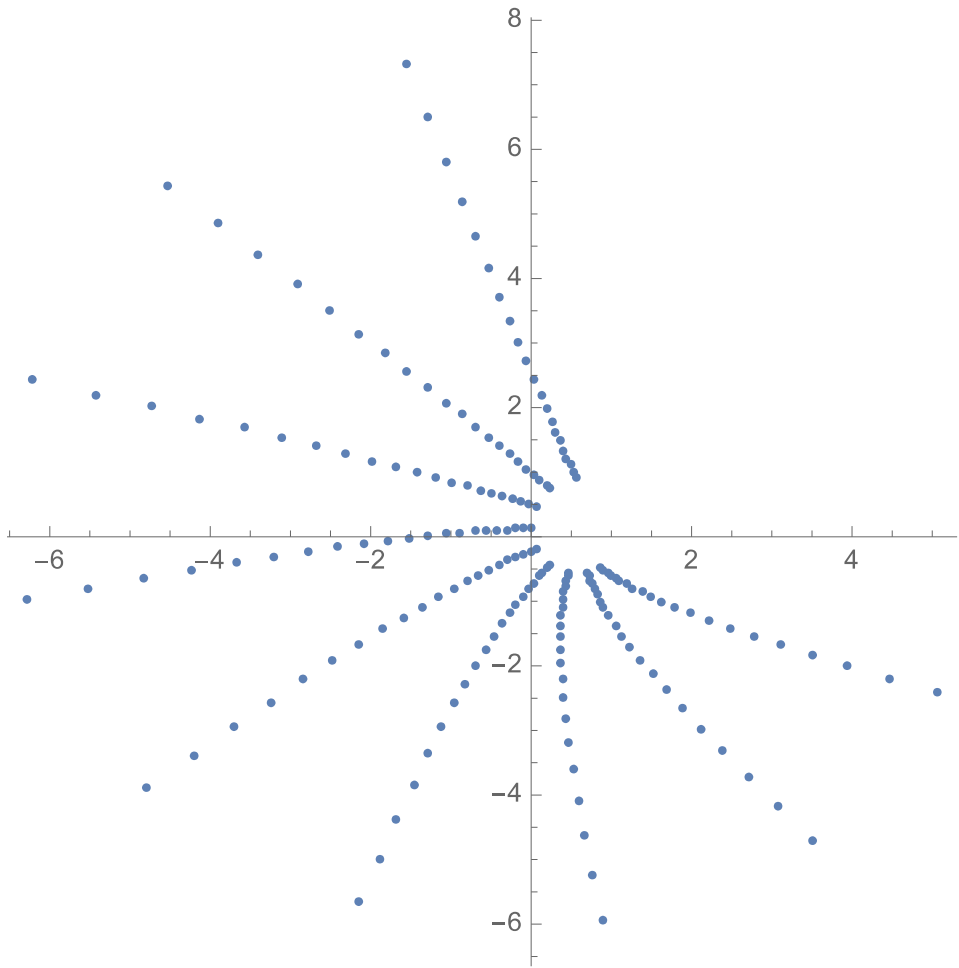

[{σ, 0, 1, .05}, {t, 47.6, 49, .2}]
There is a nontrivial zero with $t = 48.005$ ... The eleventh dot (counting from the left) of the third string counting clockwise is right on the origin. It is continuing to appear that the organizing centers of groups of the strings are always associated with $\sigma = 1$ ends of strings.

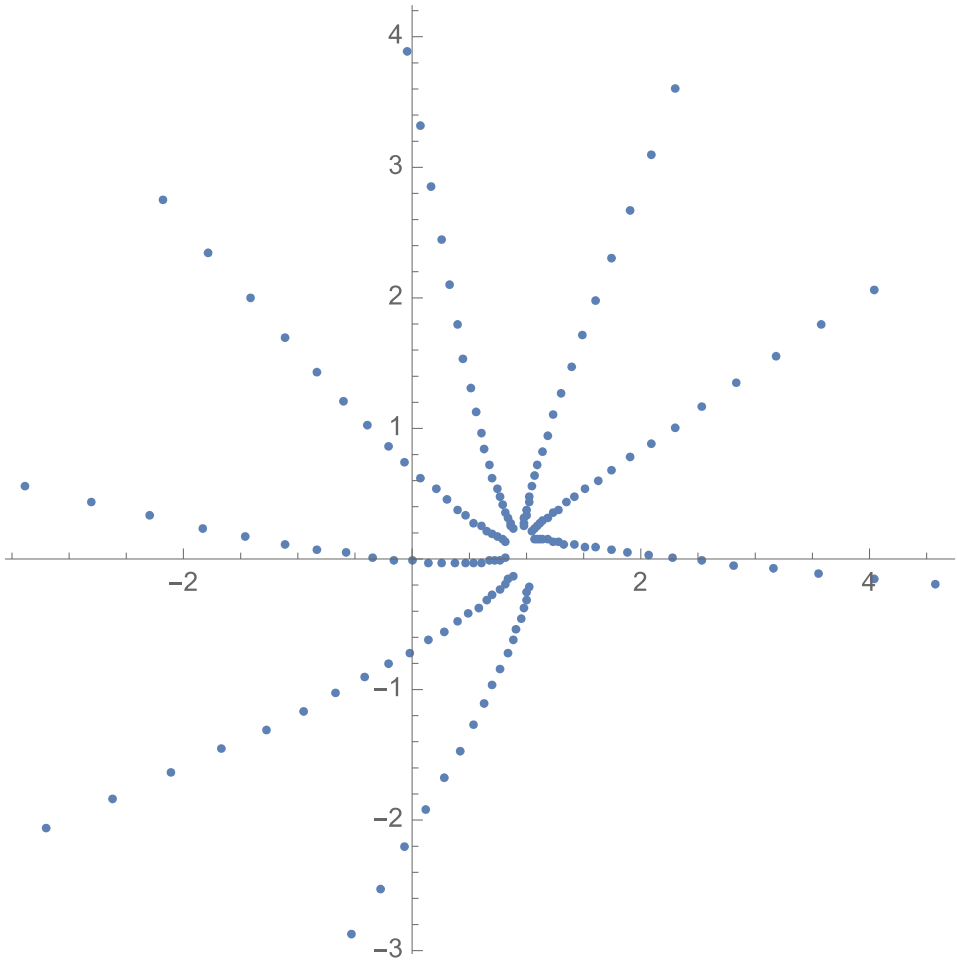

[{σ, 0, 1, .05}, {t, 49, 50.4, .2}]
There is a nontrivial zero with $t = 49.773\ldots$ in this range. The eleventh dot of the fifth string ($t = 49.8$) counting clockwise is almost on the origin.

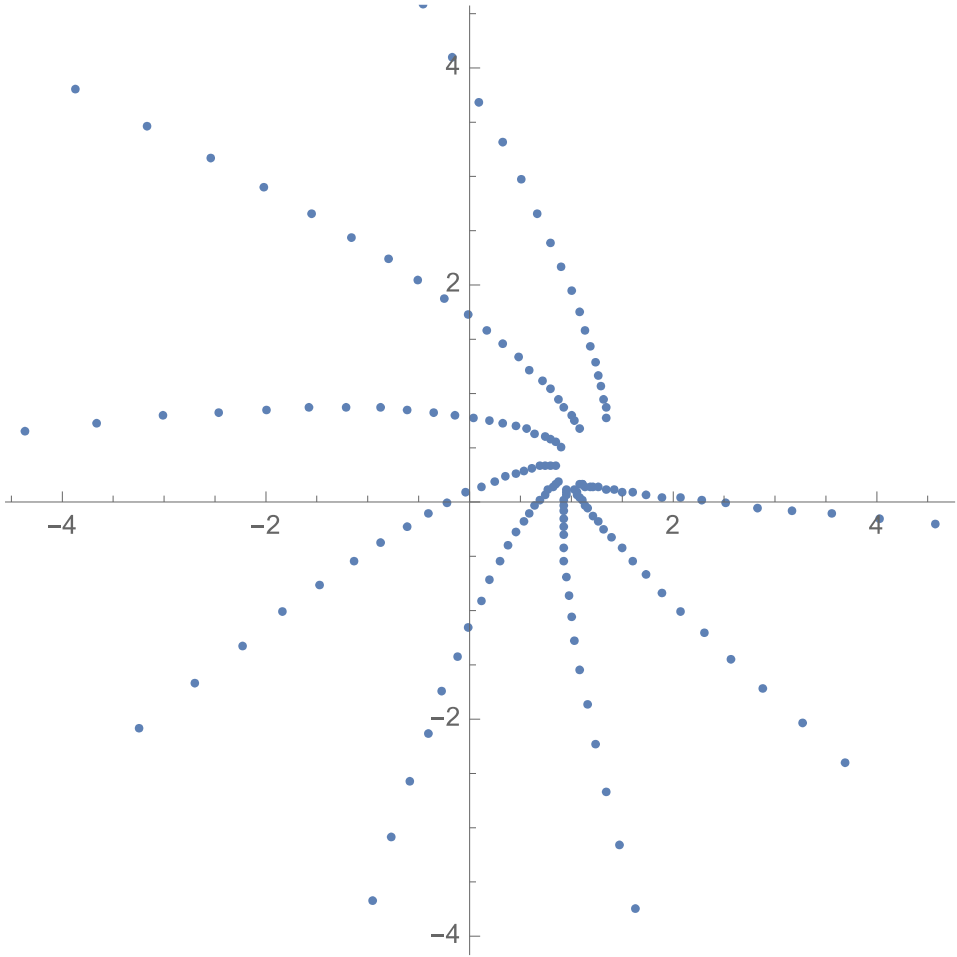

[{σ, 0, 1, .05}, {t, 50.4, 52, .2}]
In this range there are no zeros.

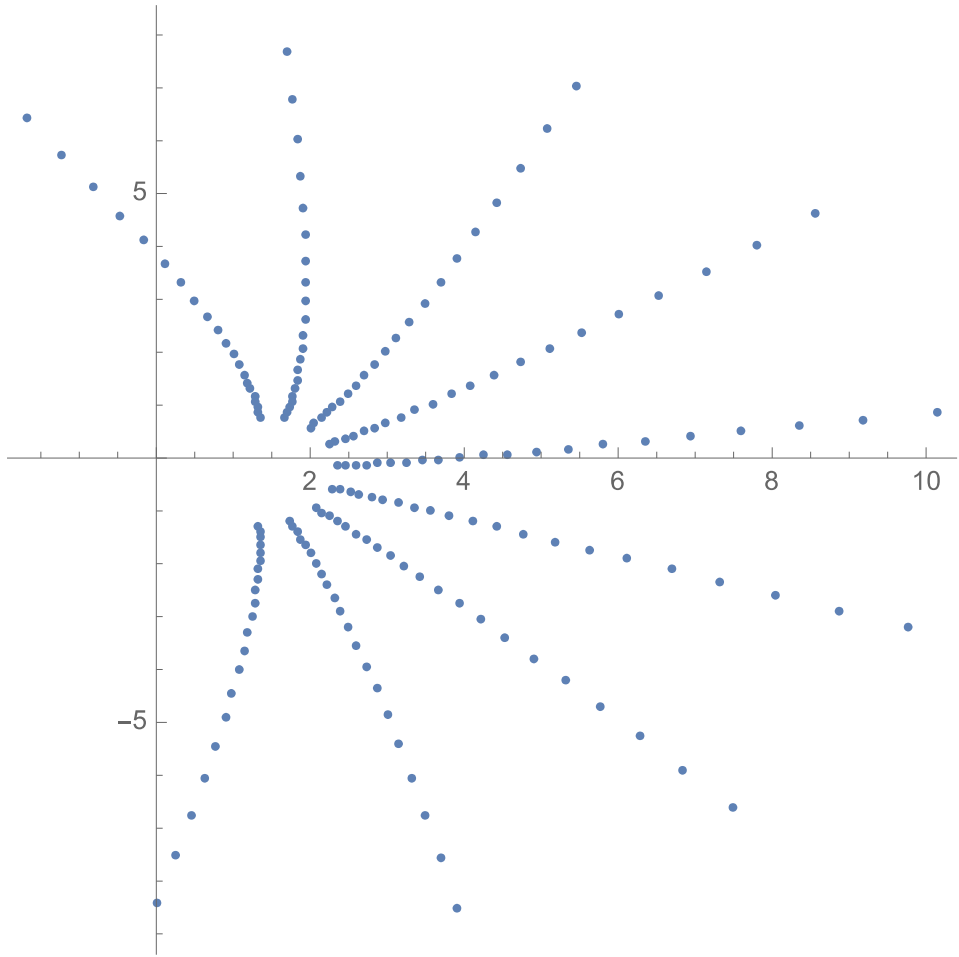

[{σ, 0, 1, .05}, {t, 52, 53.4, .2}]
There is a zero with t = 52.970 ... Counting clockwise the sixth string (t = 53) at 12 o'clock has its eleventh dot (counting down) just to the right of the origin. Distinct dots appear in these plots for σ less than or equal to 0.5 and dots are forming an indistinct smear for σ greater than 0.5.

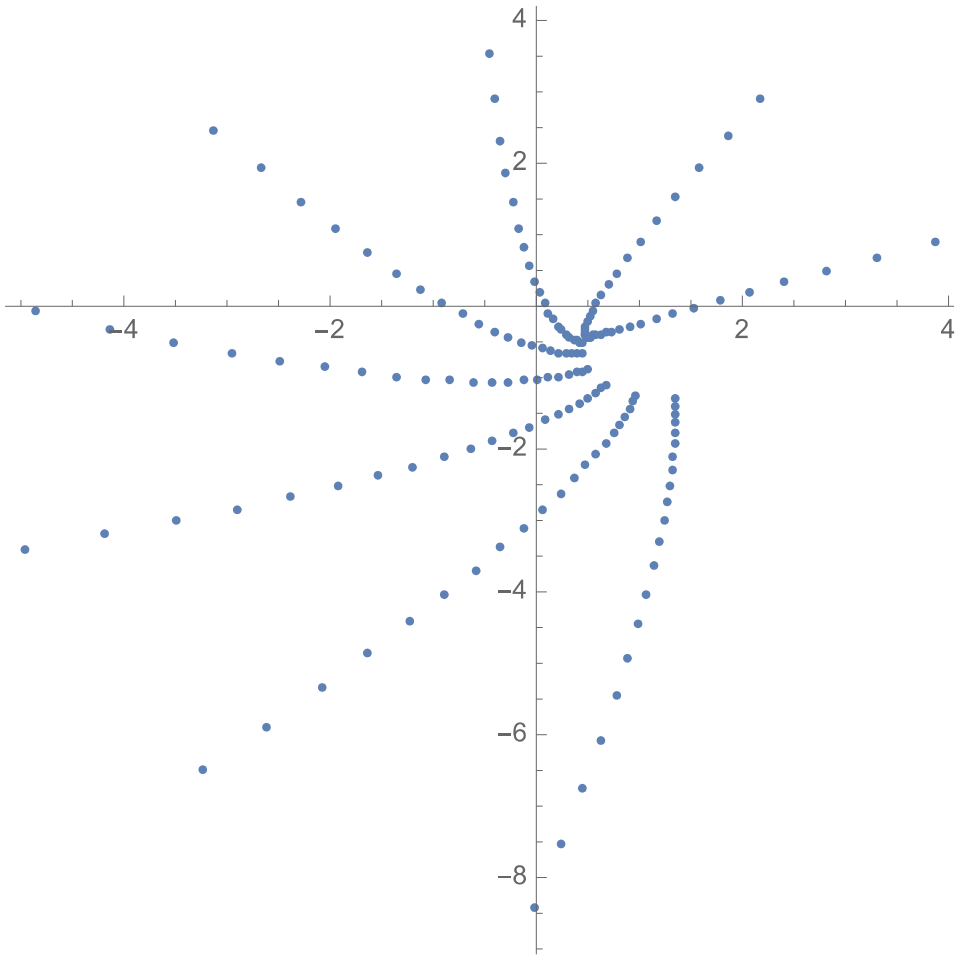

[{σ, 0, 1, .05}, {t, 53.4, 55, .2}]
There is a trivial zero with $t = 54.3882$ ... ($k = 6$). Counting clockwise the sixth string with $t = 54.4$ has the Eta value for $\sigma = 1$ right at the origin.

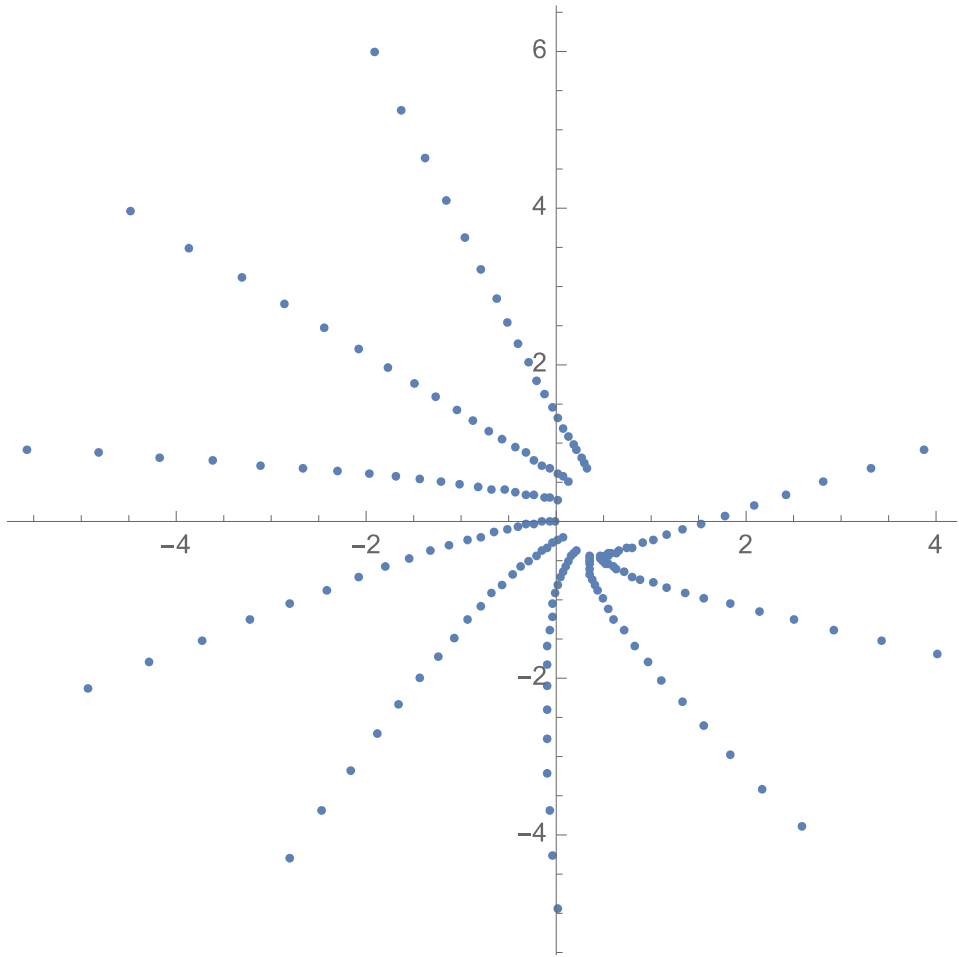

[{σ, 0, 1, .05}, {t, 55, 56.4, .2}]

There is a nontrivial zero with $t = 56.446\ldots$ This is just beyond the range for $t$ in this figure and will be discussed in the next figure. Notice here how tight the center of rotation has become. This is partly because the $t$ string lengths are increasing and partly because it is novel behavior.

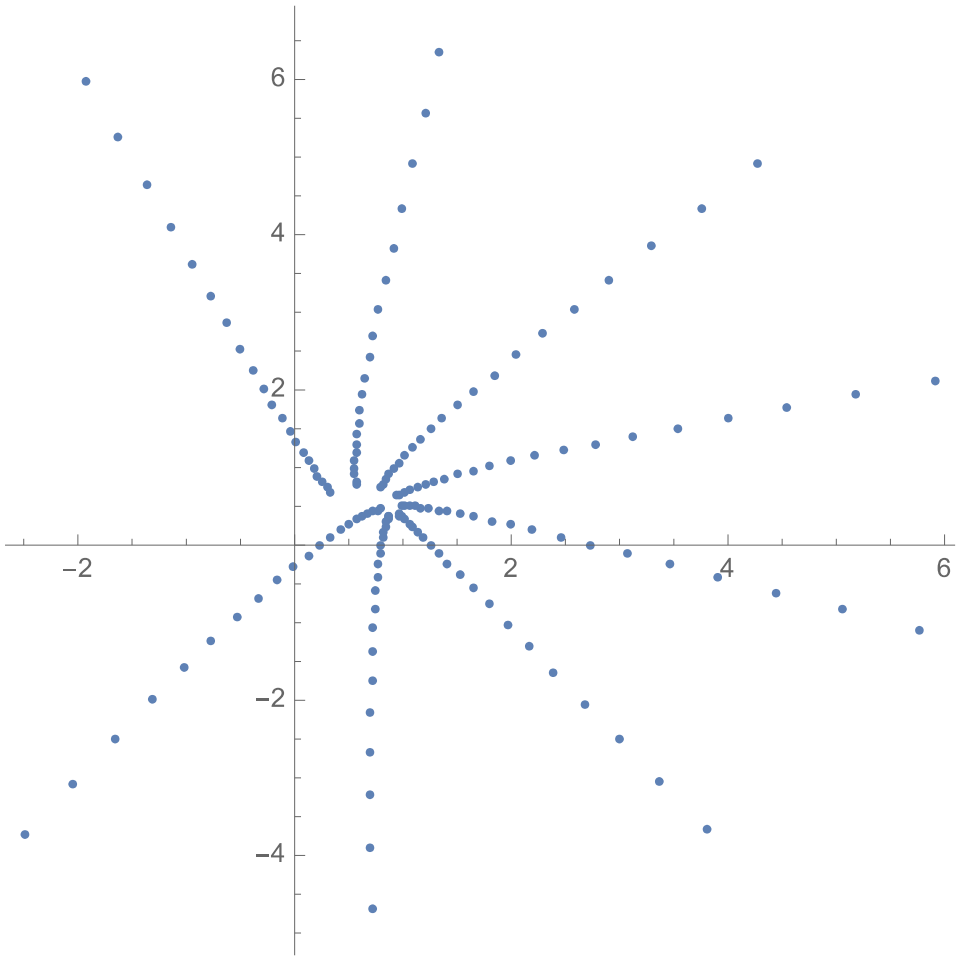

[{σ, 0, 1, .05}, {t, 56.4, 58, .2}]
Now $t = 56.446 \ldots$ is in the exhibited range. It is close to the first string in the figure at about 7 o'clock. The eleventh dot is closest to the origin. The difference in dot density between $\sigma < 0.5$ and $\sigma > 0.5$ is growing more dramatic.

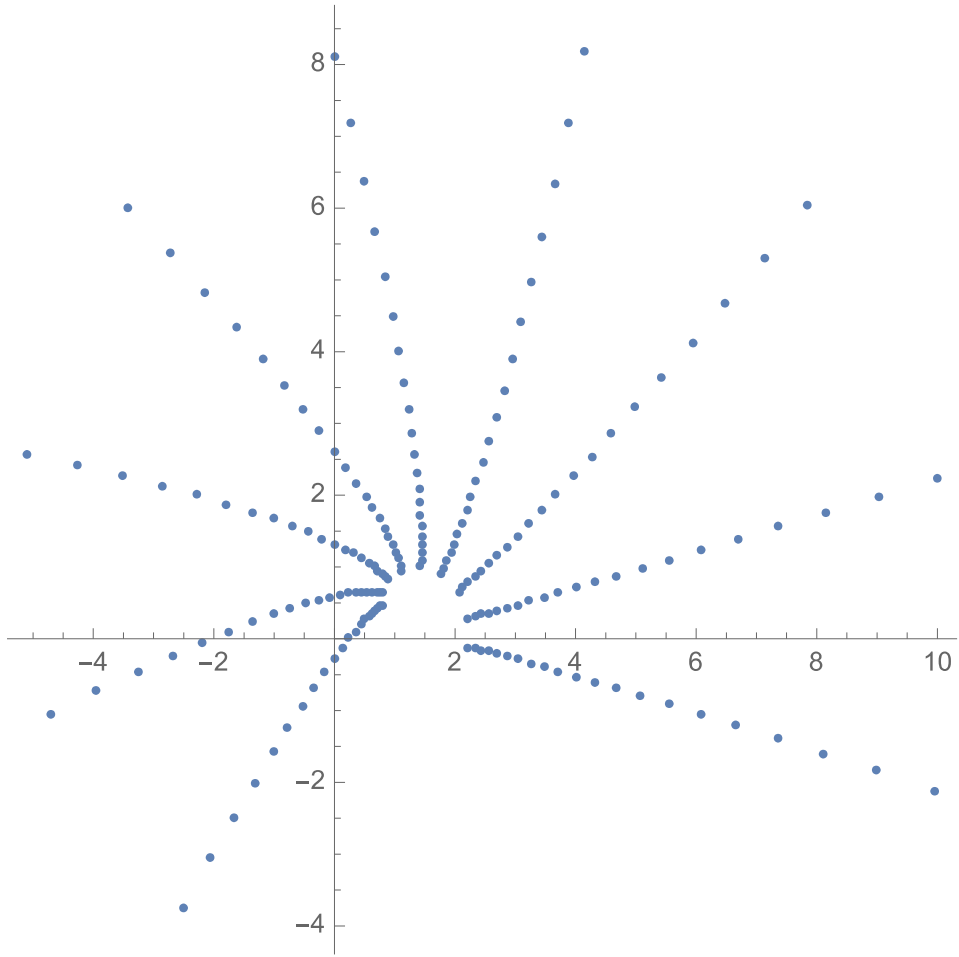

[{σ, 0, 1, .05}, {t, 58, 59.6, .2}]
There is a nontrivial zero at $t = 59.347$ ... This is near the eighth string at $t = 59.4$. Counting down the dots the eleventh dot is just above and to the right of the origin.

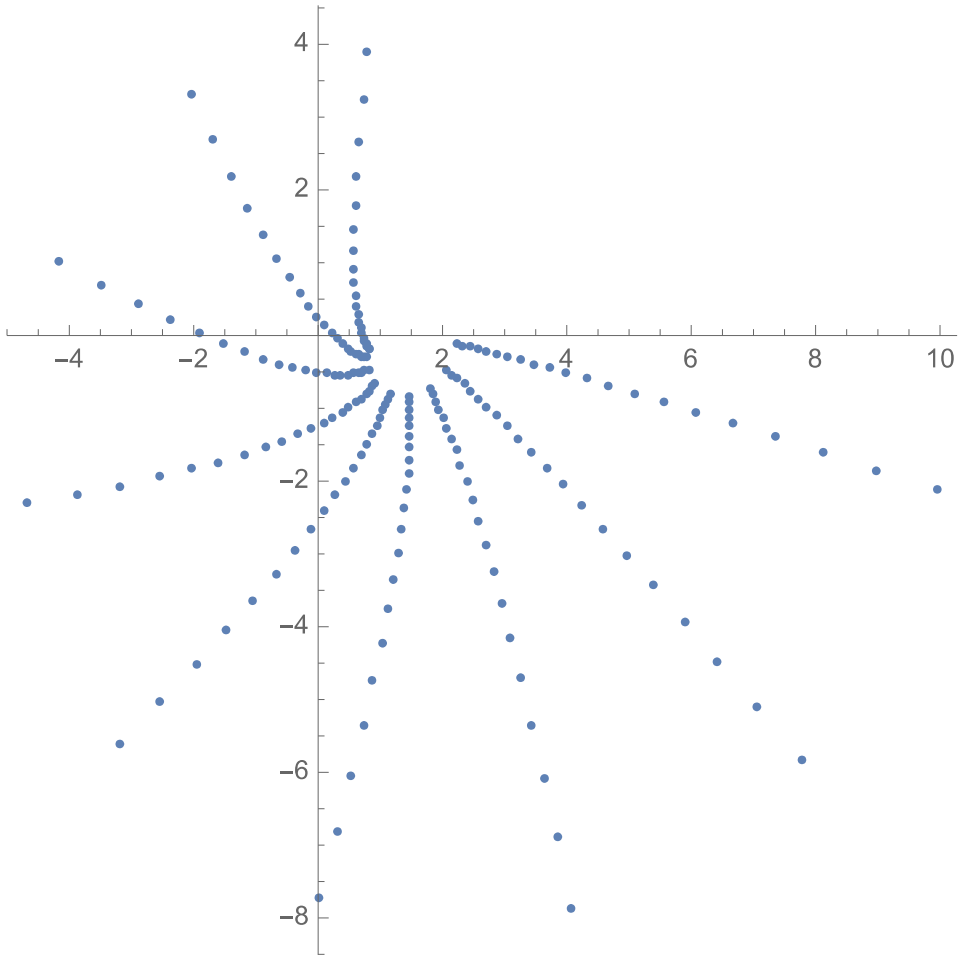

[{σ, 0, 1, .05}, {t, 59.6, 61, .2}]
There is another nontrivial zero at $t = 60.831$ ... This time there is a complicated morass and the plot has been cut off slightly to the right of $-2$ on the abscissa. From the figure above it is seen that the $t = 59.6$ string is the most vertical one. Counting around clockwise to the seventh string ($t = 60.8$) shows it coming in from the left a bit below the origin. The dot just below the origin is the eighth dot from the left but three more dots of this string have been cut off. So, it is the eleventh dot after all.

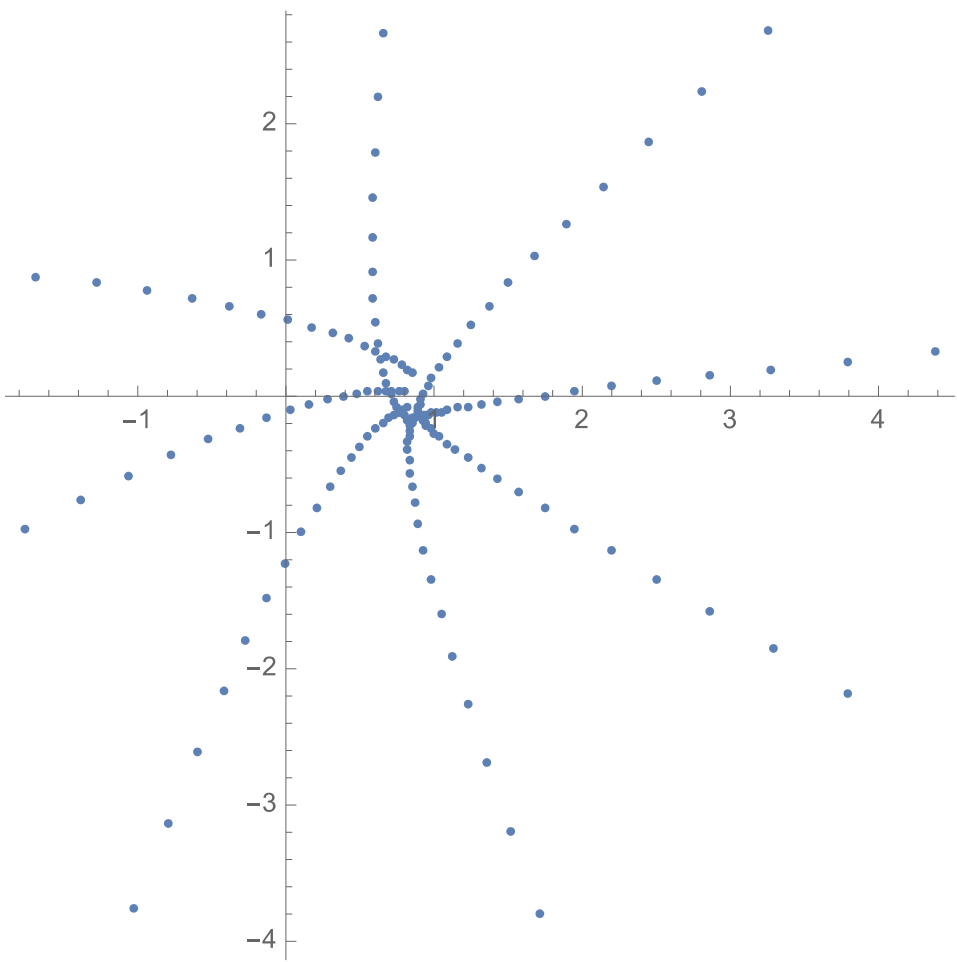

[{σ, .5, 1, .05}, {t, 59.6, 61, .2}]

This figure is a magnification of the region nearer to the origin created by cutting off σ < 0.5. This makes it clearer where the string ends are and how the curvature is changing from that of the first three strings to that of the last five.

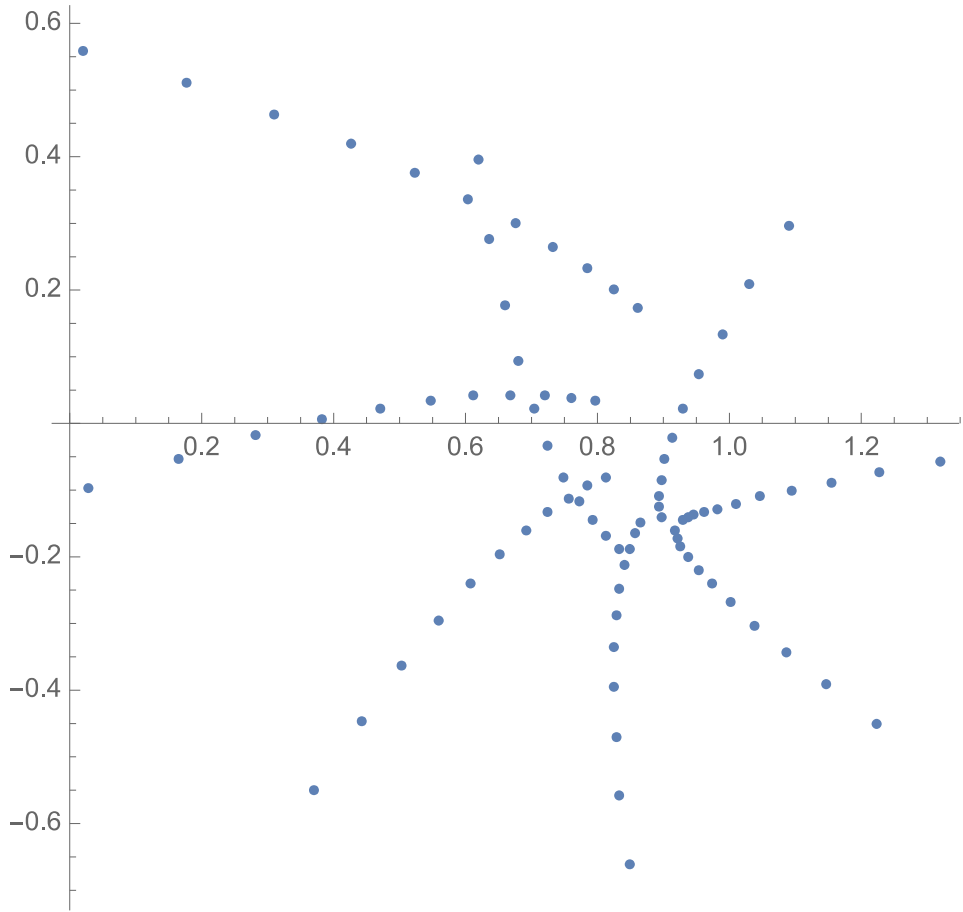

[{σ, .5, 1, .05}, {t, 61, 62.4, .2}]
There are no zeros in this range. The $\sigma < 0.5$ dots have been cut off in order to magnify the behavior near the origin. Why these strings are so nearly straight has not been explained.

[{σ, .5, 1, .05}, {t, 62.4, 64, .2}]

There is a trivial zero for $t = 63.4529$ ... (k = 7). The $t = 63.4$ string is close to this value of $t$ and is the sixth string counting clockwise. The $\sigma = 1$ dot for this string is just below the origin. Again the $\sigma < 0.5$ dots have been omitted.

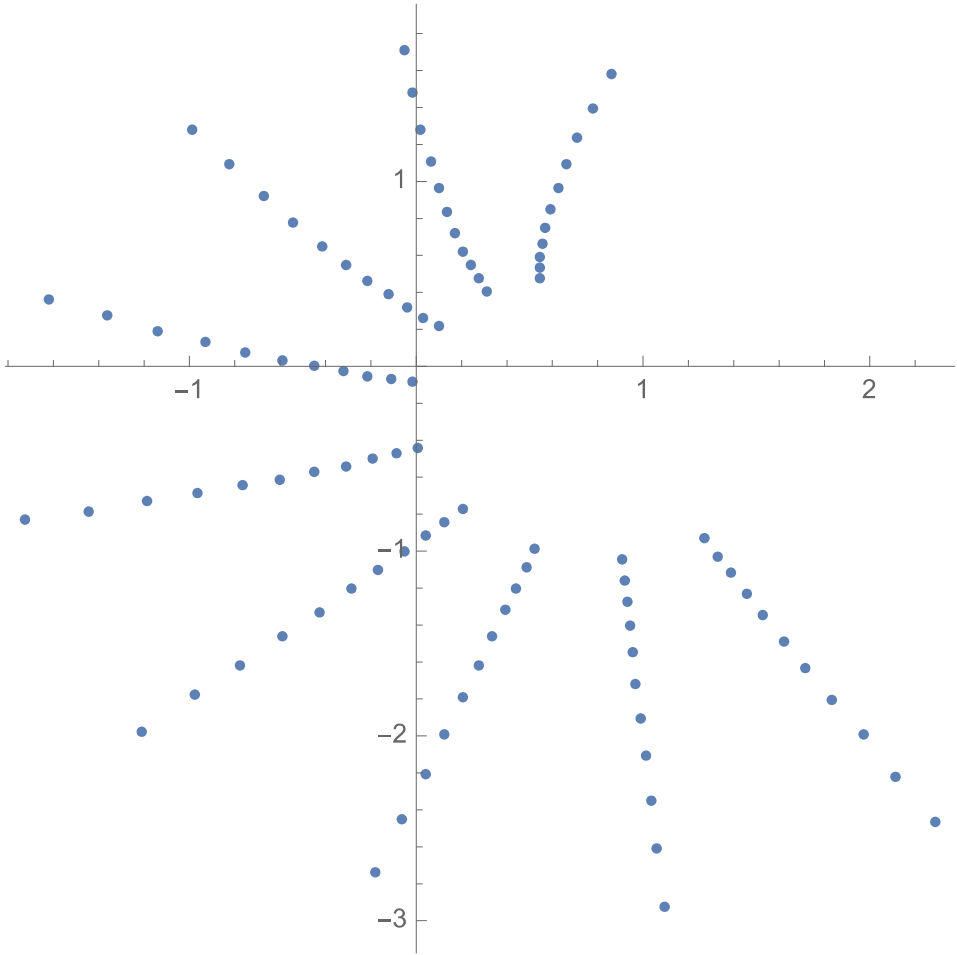

[{σ, .5, 1, .05}, {t, 64, 65.6, .2}]

There is a nontrivial zero in this range with $t = 65.112\ldots$ This is almost halfway between the $t = 65$ and the $t = 65.2$ strings in the figure. It is clear that a string there would cross the origin although perhaps one more dot is needed to reach that far. That means using a sigma range of [{σ, .45, 1, .05}]. In addition, the strings are getting so far apart that there is overlap of strings.

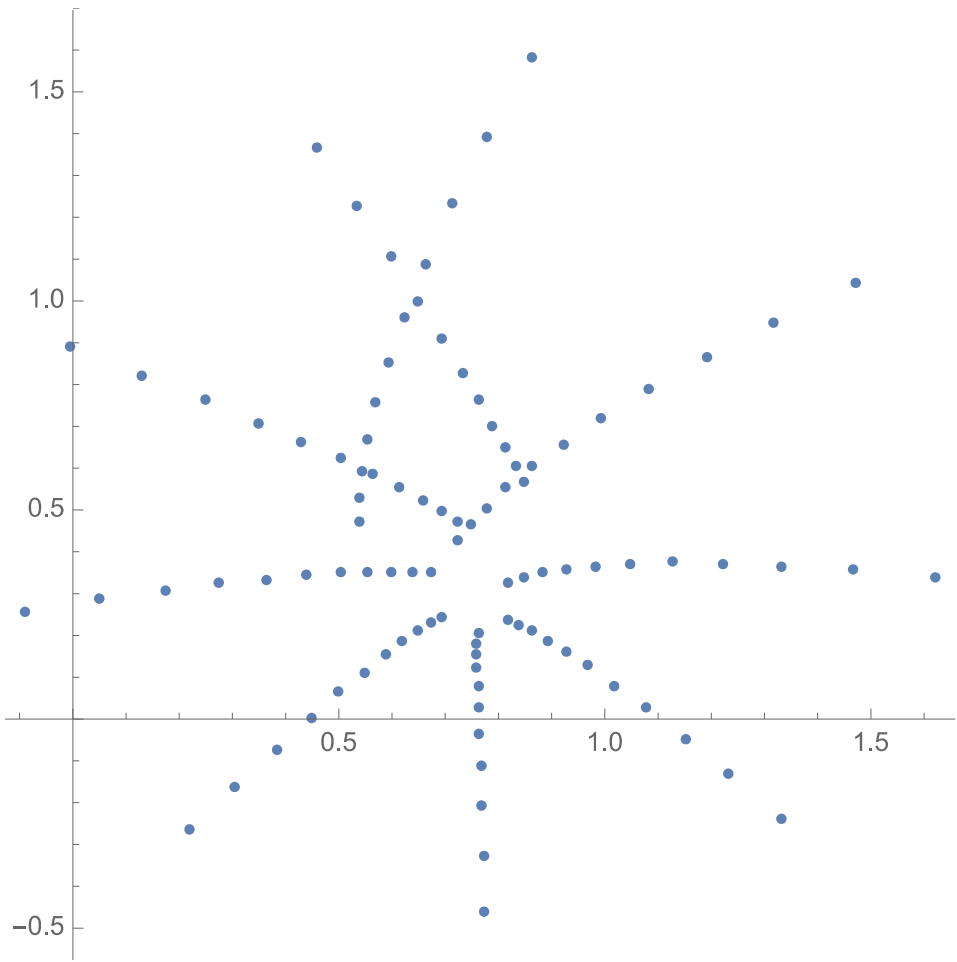

[{σ .5, 1, .05}, {t, 64, 65.2, .1}]
To eliminate the overlap in the previous figure the range for $t$ is shortened and the spacing between $t$ values is halved. A twelfth dot proves to be unnecessary and the eleventh string counting clockwise has its eleventh dot just below the origin.

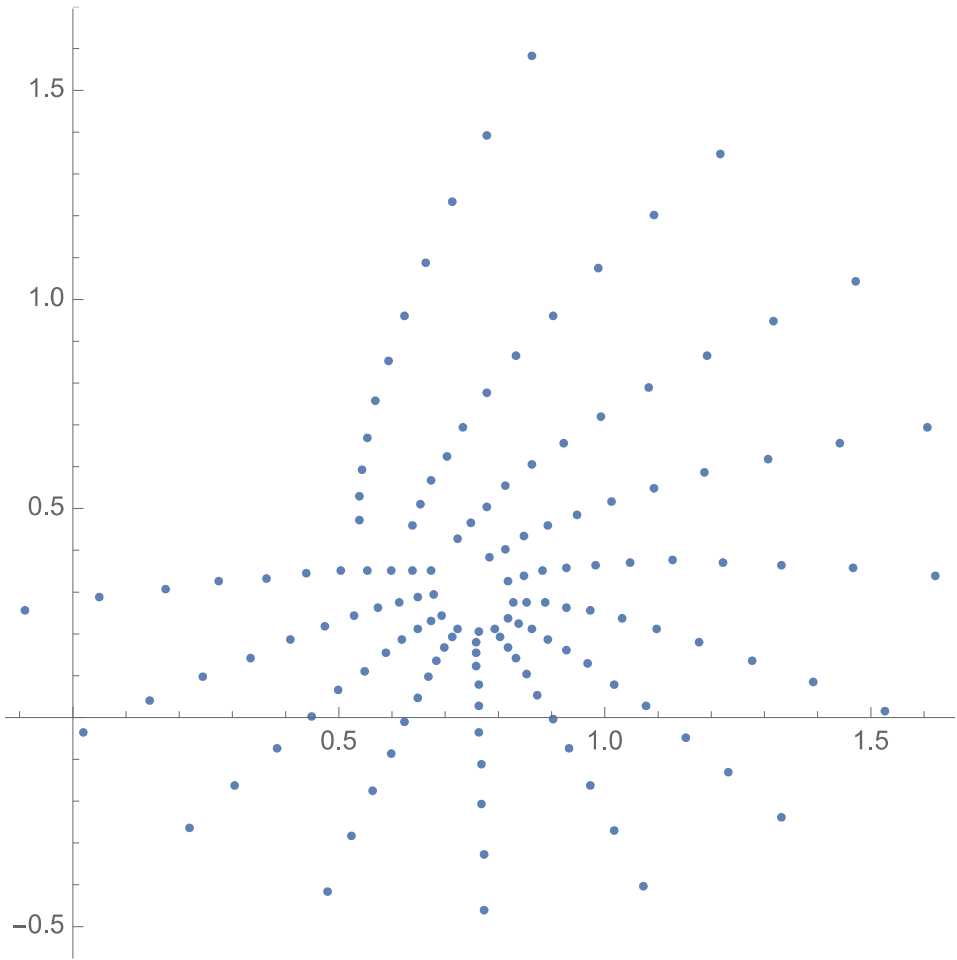

[{σ, .5, 1, .05}, {t, 66.4, 67.6, .1}]
There is a nontrivial zero at $t = 67.079$ ... The eighth string, counting clockwise, is the $t = 67.1$ string and its eleventh dot is just above the origin.

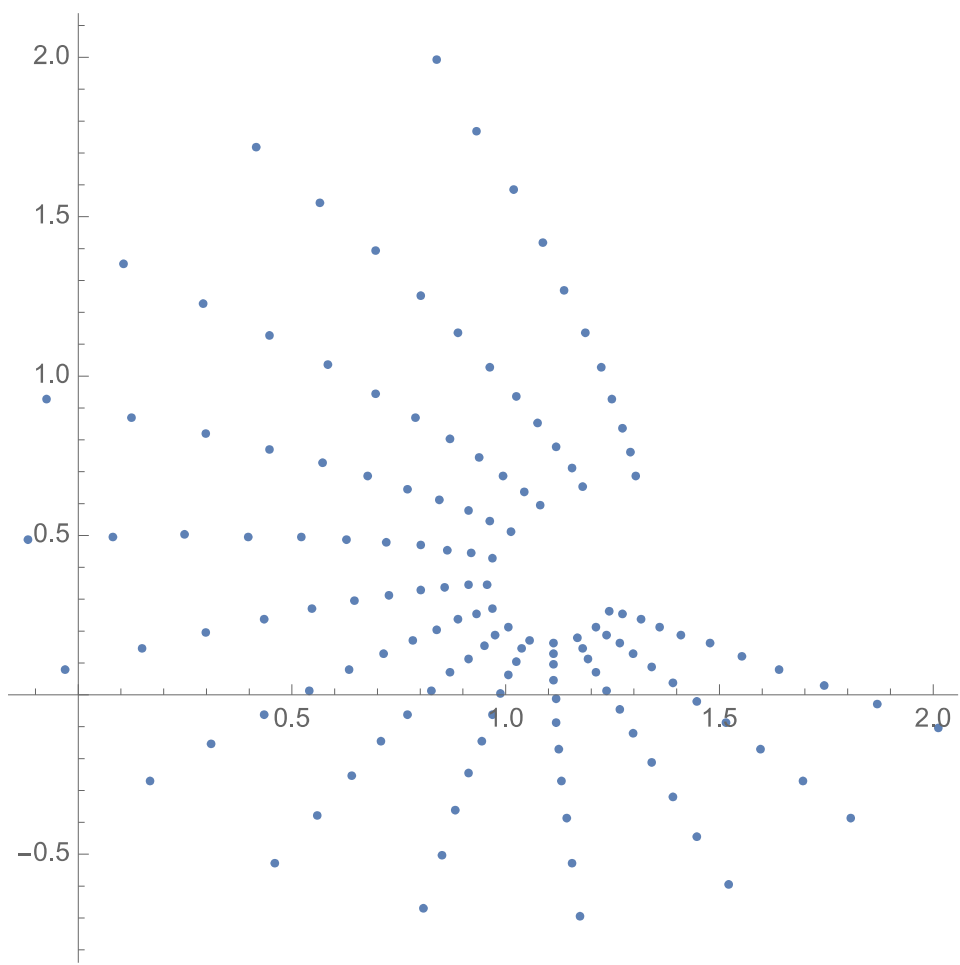

[{σ, 0, 1, .05}, {t, 66.4, 67.6, .1}]

The full range of σ values is restored. Because the small σ 's create long tails on the strings it is difficult to see the details near the $\sigma = 1$ ends of the strings. While the figure above has an area of roughly $2 \times 2.5$ the figure here is roughly $12 \times 12$. For larger $t$ values, such as 7000, the difference in scale between the $\sigma > 0.5$ and $\sigma < 0.5$ regions of the figures is much greater, as was shown in [1].

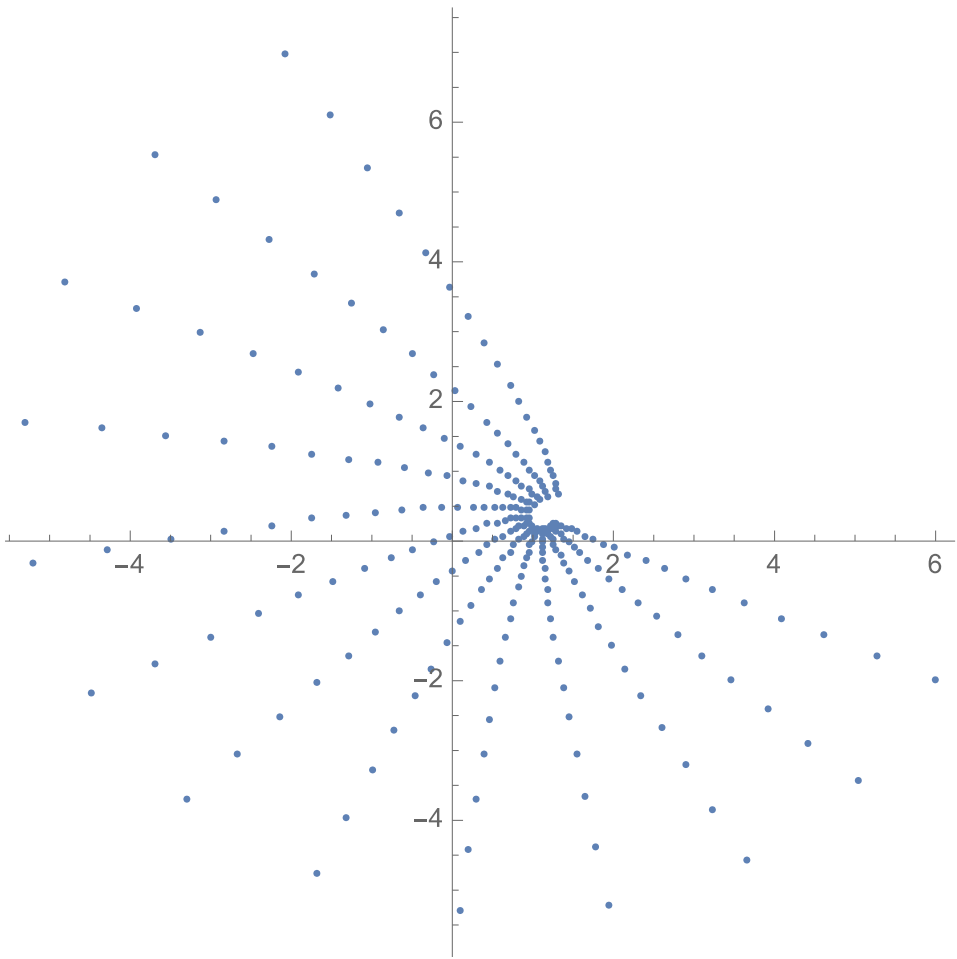

# References

[1] Riemann, G. F. B., *Über die Anzahl der Primzahlen unter einer gegebenen Grösse.* Monatsberichte der Berliner Akademie, November, 1859.

[2] Euler, L., *Variae observationes circa series infinitas.* Thesis, St. Petersburg Academy, 1737.

[3] Gourdon, X., "The $10^{13}$ first zeros of the Riemann Zeta function, and zeros computation at very large height", version October 24, 2004.

[4] Odlyzko, A. M., The $10^{22}$-th zero of the Riemann zeta function. In M. van Frankenhuysen and M. L. Lapidus, editors, Dynamical, Spectral, and Arithmetic Zeta Functions, number 290 in Contemporary Math. series, pages 139–144. Amer. Math. Soc., 2001.

[5] Perez-Marco, R., "Notes on the Riemann Hypothesis.", *Jornadas Sobre Los Problemas Del Milenio*, Barcelona, junio 2011. (see p.72)

[6] Odlyzko, A. M., http://www.dtc.umn.edu/~odlyzko/zeta_tables/zeros3